\documentclass[preprint,12pt]{elsarticle}

\usepackage[utf8]{inputenc}
\usepackage[T1]{fontenc}
\usepackage[margin=1in]{geometry}
\usepackage{hyperref}
\usepackage{url}
\usepackage{booktabs}
\usepackage{amsmath,amsfonts}
\usepackage{nicefrac}
\usepackage{microtype}
\usepackage{graphicx}
\usepackage{subcaption}
\usepackage{xcolor}
\usepackage{enumitem}
\usepackage{bm}
\usepackage{comment} %用于多行注释的包
\usepackage[ruled,vlined]{algorithm2e}

\usepackage {booktabs}
\usepackage{subcaption}
\usepackage{caption}
\usepackage{booktabs} % for professional tables
\usepackage{graphicx}
\usepackage{float}
\usepackage{amsmath}      % 支持数学公式
\usepackage{amssymb}
\usepackage{mathtools}
\usepackage{amsthm}
\usepackage{multirow}
\usepackage{makecell}   % 用于单元格内换行
\usepackage{tcolorbox}
\usepackage{adjustbox}
\usepackage{enumitem}
\usepackage{tabularx}
\usepackage{tikz}
\usetikzlibrary{positioning,arrows.meta,shapes.geometric,calc,backgrounds}
\usetikzlibrary{matrix}
\usetikzlibrary{decorations.pathreplacing}
\usetikzlibrary{calc}

\newlength{\subfigwidth}
\newlength{\insetfigwidth}
\graphicspath{{./images/}}

\usepackage[numbers]{natbib}

\begin{document}

\begin{frontmatter}

\title{Robust Physics-Guided Diffusion for Full-Waveform Inversion}
% ---------- Authors ----------
\author[aff1]{Jishen Peng\fnref{fn1}}
\author[aff1]{Enze Jiang\fnref{fn1}}
\author[aff1,aff2,aff3,aff4]{Zheng Ma}
\author[aff5]{Xiong-Bin Yan\corref{cor1}}

% ---------- Footnotes ----------
\fntext[fn1]{Jishen Peng and Enze Jiang contributed equally to this work.}
\cortext[cor1]{Corresponding author.}

% ---------- Email ----------
\ead{yanxb2015@163.com}

% ---------- Affiliations ----------
\address[aff1]{School of Mathematical Sciences, Shanghai Jiao Tong University, Shanghai, China}
\address[aff2]{Institute of Natural Sciences, MOE-LSC, Shanghai Jiao Tong University, Shanghai, China}
\address[aff3]{Qing Yuan Research Institute, Shanghai Jiao Tong University, Shanghai, China}
\address[aff4]{CMA-Shanghai, Shanghai Jiao Tong University, Shanghai, China}
\address[aff5]{School of Mathematics and Statistics, Lanzhou University, Lanzhou, China.}

\begin{abstract}
We develop a robust physics-guided diffusion framework for full-waveform inversion that combines a score-based generative prior with likelihood guidance computed through wave-equation simulations. We adopt a transport-based data-consistency potential (Wasserstein-2), incorporating wavefield enhancement via bounded weighting and observation-dependent normalization, thereby improving robustness to amplitude imbalance and time/phase misalignment. On the inference side, we introduce a preconditioned guided reverse-diffusion scheme that adapts the guidance strength and spatial scaling throughout the reverse-time dynamics, yielding a more stable and effective data-consistency guidance step than standard diffusion posterior sampling (DPS). Numerical experiments on OpenFWI datasets demonstrate improved reconstruction quality over deterministic optimization baselines and standard DPS under comparable computational budgets.
\end{abstract}

\end{frontmatter}

% keywords can be removed
%\keywords{First keyword \and Second keyword \and More}

\section{Introduction}
In geophysics, seismic inversion is widely used to infer quantitative subsurface models that match observed seismic recordings.
From a mathematical viewpoint, recovering medium parameters from time-dependent wavefield measurements constitutes a nonlinear and ill-posed inverse problem,
primarily due to limited acquisition geometry and band-limited observations. Several classes of methods have been developed, including velocity analysis based on stacked traces \cite{Berkhout+1997},
migration and traveltime-based methods \cite{zelt+1992,clement+2001}, Born-approximation-based linearized inversion \cite{hudson1981use,muhumuza2018seismic},
and full-waveform inversion (FWI) \cite{tarantola1984inversion,warner2013anisotropic,jakobsen2015full}.
Among these, FWI is distinguished by its potential to recover high-resolution subsurface structure by exploiting both phase and amplitude information through the solution of a PDE-constrained optimization problem that minimizes a chosen data-misfit functional between observed and synthetic seismograms.

Full-waveform inversion (FWI) is a PDE-constrained inverse problem that aims to reconstruct subsurface parameters (e.g., wavespeed or velocity) from time-dependent seismic recordings \cite{virieux2009overview,fichtner2010full,tarantola1984inversion}.
Given an acquisition geometry and a forward wave-propagation operator $\mathcal{F}$, the classical formulation seeks a velocity model $v$ such that $\mathcal{F}(v)$ matches the observed data $d_{\mathrm{obs}}$, where $v\in\mathbb{R}^m$ denotes the discretized velocity field and
$\mathcal{F}:\mathbb{R}^m\to\mathbb{R}^n$ maps model parameters to recorded traces.
Despite its ability to recover high-resolution structures by exploiting the full wavefield content, FWI is widely recognized as a challenging nonconvex optimization problem: the forward map is nonlinear, the inverse problem is ill-conditioned,
and the misfit landscape typically contains many local minima due to phase ambiguity and limited illumination \cite{virieux2009overview,yao2019tackling,yang2018application,pratt1999seismic}.
These difficulties are especially pronounced when the initial model is inaccurate and when the observations are contaminated by noise (and other data errors).
In addition, strong multi-scale components with markedly different amplitudes can further bias residual-based objectives and degrade stability
\cite{bunks1995multiscale,virieux2009overview}.

Since pointwise $\ell_2$ objectives penalize the pointwise residual in time, they are highly sensitive to small time/phase shifts between observed and synthetic traces.
This sensitivity leads to a strongly nonconvex objective landscape and the well-known cycle-skipping phenomenon \cite{gauthier1986two,metivier2016measuring}, in which iterative
methods can become trapped in local minima associated with incorrect phase (kinematic) alignment.
Moreover, seismic traces often exhibit pronounced dynamic-range imbalance:
early-arriving phases can have substantially larger amplitudes than later reflections, causing both the misfit and its gradient
to be dominated by a small portion of the data \cite{liu2017robust}. As a consequence, the resulting gradients can be dominated by high-amplitude components, leading updates to underfit weaker later phases.
A large body of work aims to mitigate these issues through multi-scale continuation \cite{bunks1995multiscale}, alternative misfit functions
based on phase/envelope attributes \cite{bozdaug2011misfit,chi2014full}, adaptive waveform inversion \cite{warner2016adaptive},
and geometry-aware metrics such as optimal-transport/Wasserstein distances \cite{yang2018application,qiu2017full,li2022quadratic};
see also \cite{gao2023review} for a recent review of misfit functions and adjoint sources.

In parallel with advances in data misfit design, recent years have witnessed rapid progress in learning expressive priors for
inverse problems using generative modeling. In seismic imaging, a range of generative models---including GAN-based geological priors, latent-variable models such as VAEs, and normalizing flows---have been explored
to mitigate ill-posedness and to facilitate uncertainty quantification \cite{mosser2020stochastic,lopez2022geophysical,sun2024enabling}.
Score-based diffusion models (also known as score-based generative models) provide a particularly attractive mechanism: they learn the
score field $\nabla_x \log p_t(x)$ of progressively noised marginals and enable conditional generation by reversing the associated diffusion
dynamics \cite{ho2020denoising,song2020score}.
From the perspective of inverse problems, diffusion models can be interpreted as providing a learned prior over plausible velocity fields, trained solely from samples of the velocity model (i.e., parameter-field samples), without requiring paired data of the form (wavefield, velocity)
and without embedding a forward solver in the prior-training loop; conditional inference is then performed at test time by incorporating data-consistency
(likelihood) information during reverse-time sampling or related posterior-sampling schemes \cite{song2021solving,chung2022diffusion}.
This decoupling is especially appealing in seismic imaging, where repeated wavefield simulation is expensive and acquisition geometries may
vary across applications; accordingly, diffusion-based Bayesian FWI and diffusion-driven velocity modeling have started to be investigated in
recent work \cite{li2025diffusioninv,wang2024controllable}.

To solve the inverse problem given observations, guided reverse diffusion methods combine the learned prior score with a data-consistency term induced by a potential $\Phi$ (i.e., a negative log-likelihood, or more generally a data-misfit functional that measures agreement between simulated data and the observed data). Diffusion posterior sampling (DPS) is a representative method: at reverse step $i$ it forms a clean-state estimate $\hat v_0^{(i)}$, which is a denoised estimate of the underlying clean variable associated with the current noisy state. This estimate is computed using the learned prior score $s_{\theta^\ast}(\cdot,i)$. A likelihood-guidance step is then applied using the gradient $\nabla \Phi(\hat v_0^{(i)})$. However, a direct use of the baseline DPS update with a fixed global guidance step size is poorly matched to the characteristics of full-waveform inversion (FWI) in two fundamental respects.

First, the effectiveness of DPS critically depends on the choice of the potential $\Phi$.
In FWI, conventional pointwise $\ell_2$ (residual-based) potentials inherit two structural deficiencies.
On the one hand, the quadratic residual makes the induced descent direction disproportionately influenced by large-amplitude arrivals, resulting in markedly unbalanced updates across traces and time windows.
On the other hand, the pointwise $\ell_2$ geometry is highly sensitive to small time/phase shifts, which typically produces a highly nonconvex misfit landscape
with spurious local minima (cycle-skipping) associated with incorrect phase alignment.
When such a potential is used for likelihood guidance in DPS, these effects can yield correction terms whose scale and direction are poorly matched to the current reverse-time state,
thereby degrading the stability and effectiveness of the reverse-time dynamics.

Second, even with a robust potential, a single scalar guidance strength in baseline DPS cannot capture the spatially varying sensitivity of the misfit to the velocity field $v$. Early reverse-time steps often yield a rough $\hat v_0^{(i)}$, so
overly aggressive guidance can be unreliable, whereas later steps typically benefit from stronger data-consistency enforcement. Moreover,
illumination and acquisition geometry induce strong spatial variability in sensitivity, so a single global step size may overshoot in
well-illuminated regions while remaining ineffective in poorly illuminated ones. These mismatches can reduce stability and efficiency along the reverse-time iterations and ultimately degrade reconstruction quality.

\textit{Contributions.}
This work develops a robust physics-guided diffusion framework for full-waveform inversion (FWI), combining a score-based generative prior with an OT-based data-consistency potential and a stabilized, variable-metric guidance mechanism in guided reverse diffusion.
Our main contributions are:
\begin{itemize}
  \item \textbf{An OT-based data-consistency potential for physics-guided diffusion.}
We design a robust likelihood potential for guided reverse diffusion by combining bounded amplitude-adaptive weighting of seismic traces
with a one-dimensional Wasserstein discrepancy computed via quantile functions. When used as the data-consistency term in DPS-type guidance, the resulting guidance signal is less dominated by high-amplitude arrivals and less sensitive to small time/phase misalignment, thereby alleviating cycle-skipping and yielding a more balanced and stable data-consistency mechanism for FWI.

  \item \textbf{Adaptive, variable-metric guidance for guided reverse diffusion.}
  We introduce a variable-metric (diagonal) preconditioner for the guidance step, replacing the scalar DPS step size by
  a matrix $P_i=\rho_i D_i$.
  Here $\rho_i$ adapts the guidance strength along the reverse trajectory to suppress unreliable early updates,
  while the diagonal scaling $D_i$ balances the spatially heterogeneous sensitivity of the misfit with respect to the velocity field in FWI. This yields a more stable and effective guidance mechanism than standard DPS in seismic settings.

  \item \textbf{Numerical validation.}
  We evaluate the proposed method on OpenFWI benchmarks and compare against deterministic optimization-based inversion baselines
  and the original DPS under matched computational budgets, demonstrating improved reconstruction quality and stability.
  We further assess generalization on other standard benchmarks beyond OpenFWI, where the method maintains strong performance.
\end{itemize}

%Although our formulation admits a Bayesian interpretation through a %learned prior $p(v)$ (via diffusion) and a potential
%$\Phi$ (via OT-based data consistency), the emphasis of this paper is on robust reconstruction and stabilized conditional sampling,
%rather than comprehensive uncertainty quantification.

The remainder of the paper is organized as follows.
In Section~\ref{sec:setting} we state the FWI problem and cast it in Bayesian form.
Section~\ref{sec:diffusion} recalls score-based diffusion models and the learned prior score needed for reverse-time dynamics.
Section~\ref{sec:baseline_dps} derives a conditional reverse-diffusion formulation for Bayesian inversion and presents a baseline discretization in the spirit of diffusion posterior sampling (DPS).
Section~\ref{sec:method} introduces our robust physics-guided diffusion scheme, including an OT-based data-consistency potential and an adaptive, variable-metric guidance mechanism.
Section~\ref{sec:exp} provides numerical results, comparative studies on OpenFWI benchmarks, and generalization tests on additional standard benchmark datasets.
Section~\ref{sec:conclusion} concludes with a summary of the main results.

\section{Problem formulation}\label{sec:setting}

\subsection{Full-waveform inversion (FWI)}\label{subsec:fwi}
Full-waveform inversion (FWI) aims to recover a spatially varying wave speed $v(\mathbf{x})$ from seismic observations recorded at receivers $\{\mathbf{x}_r\}_{r=1}^{N_r}$. Given observed seismograms $d_r(t)=d(\mathbf{x}_r,t)$, the estimate of $v$ is typically obtained by minimizing a waveform-misfit functional between observed and simulated data, subject to the governing wave equation. By exploiting the full waveform information---including phase, amplitude, and multiple arrivals---FWI can in principle achieve higher resolution than kinematic approaches such as traveltime tomography \cite{phillips1991traveltime,chen2023adjoint,hao2024topography}. At the same time, the resulting PDE-constrained optimization problem is strongly nonlinear and generally nonconvex; robust performance therefore requires careful treatment of noise and modeling errors, together with appropriate regularization or prior information.

In what follows, the simulated data are generated by a forward operator $\mathcal{F}(v)$ defined through the solution of an acoustic wave equation. Specifically, we consider seismic wave propagation under the constant-density acoustic approximation in an isotropic medium. The wavefield $u(\mathbf{x},t)$ denotes the acoustic pressure, and the unknown parameter is the spatially varying wave speed $v(\mathbf{x})>0$. Neglecting variable-density and anisotropic effects leads to the standard Laplacian as the spatial operator. Let $\Omega\subset\mathbb{R}^d$ $(d=2,3)$ be the domain and $T>0$ the final time. For a given source term $s(\mathbf{x},t)$, the forward problem reads
\begin{equation}\label{eq:solving}
\begin{cases}
\dfrac{1}{v^2(\mathbf{x})}\partial_{tt}u(\mathbf{x},t)-\Delta u(\mathbf{x},t)=s(\mathbf{x},t),
& (\mathbf{x},t)\in \Omega\times(0,T],\\[4pt]
u(\mathbf{x},0)=0,\qquad \partial_t u(\mathbf{x},0)=0,
& \mathbf{x}\in\Omega,\\[4pt]
\partial_n u(\mathbf{x},t)=0,
& (\mathbf{x},t)\in \Gamma_{\mathrm{ref}}\times(0,T],
\end{cases}
\end{equation}
where $\Gamma_{\mathrm{ref}}\subset\partial\Omega$ denotes the (top) reflecting boundary and $\partial_n$ is the outward normal derivative. On the remaining boundary $\Gamma_{\mathrm{abs}}:=\partial\Omega\setminus\Gamma_{\mathrm{ref}}$, we impose absorbing boundary conditions, implemented in practice via perfectly matched layers (PML), in order to minimize artificial reflections from the truncation of the computational domain.
 
We take the domain to be the rectangle 
$\Omega=[0,L_x]\times[0,L_z]\subset\mathbb{R}^2$ with $\mathbf{x}=(x,z)$, where the reflecting boundary corresponds to the top side $\Gamma_{\mathrm{ref}}=\{(x,z)\in\partial\Omega:\ z=L_z\}$. We discretize $\Omega$ by a tensor-product grid with nodes $\{(x_i,z_j)\}_{i=1,\dots,m_x;\, j=1,\dots,m_z}$. The wave-speed field is represented by its nodal values $v_{ij}:=v(x_i,z_j)$, collected in a matrix $V=\{v_{ij}\}\in\mathbb{R}^{m_x\times m_z}$, or equivalently by the vectorization $v=\mathrm{vec}(V)\in\mathbb{R}^{m}$, where $m=m_x m_z$.

Receivers are deployed exclusively on the top boundary $\Gamma_{\mathrm{ref}}\subset\partial\Omega$.
We denote the discrete receiver set by
$\Gamma_r := \{\mathbf{x}_r\}_{r=1}^{N_r}\subset \Gamma_{\mathrm{ref}}$. For the $j$-th experiment with source term $s(\mathbf{x},t;\boldsymbol{\xi}_j)$, we solve \eqref{eq:solving} (with this source) to obtain the wavefield $u_j(\mathbf{x},t;v,\boldsymbol{\xi}_j)$. Let $R$ denote the receiver sampling (restriction) operator, defined by
\begin{equation*}
(Ru_j)_r(t):=u_j(\mathbf{x}_r,t;v,\boldsymbol{\xi}_j),\qquad r=1,\dots,N_r,
\end{equation*}
so that the corresponding noise-free synthetic data for experiment $j$ are
$d^{(j)}_r(t)=(Ru_j)_r(t)$. Stacking the traces across all experiments and receivers defines the parameter-to-observable operator
\begin{equation*}
\mathcal{F}:\ v\longmapsto \{(Ru_j)_r(t)\}_{j=1,\dots,N_s;\, r=1,\dots,N_r},
\end{equation*}
which maps the velocity model to the predicted receiver recordings. After discretizing time at $\{t_k\}_{k=1}^{N_T}$, the data can be viewed as an array
\begin{equation*}
d=\{d_{j,r}(t_k)\}\in\mathbb{R}^{N_s\times N_r\times N_T},\quad
(\mathcal{F}(v))_{j,r,k}=(Ru_j)_r(t_k).
\end{equation*}
The observation model is then
\begin{equation}\label{eq:data_model}
d_{\mathrm{obs}}=\mathcal{F}(v)+\eta,
\end{equation}
where $\eta$ denotes the measurement noise. We adopt an additive Gaussian noise $\eta\sim\mathcal{N}(0,\Sigma)$, where $\Sigma$ is the noise covariance in data space. In the simplest case of independent and identically distributed noise, $\Sigma=\sigma^2 I$, so that $\eta\sim\mathcal{N}(0,\sigma^2 I)$ with $\sigma>0$ controlling the noise level.

Classically, FWI is formulated as the minimization of a data-misfit functional (possibly augmented by regularization).
In this work, we instead adopt a Bayesian formulation and infer $v$ through the posterior, where data consistency is encoded by the likelihood potential $\Phi$ (i.e., the negative log-likelihood), defined below.

\subsection{Bayesian formulation of FWI}\label{subsec:bayes_fwi}
Given the observation model \eqref{eq:data_model}, we cast FWI as a Bayesian inverse problem
and seek the subsurface velocity $v$ through the posterior distribution $p(v | d_{\mathrm{obs}})$.
In this work, we use the posterior formulation primarily as a computational target for a diffusion-based inversion algorithm,
rather than pursuing full posterior uncertainty quantification.

\textit{Likelihood as a data-consistency potential.}
Rather than fixing the likelihood to a pointwise $\ell_2$ discrepancy a priori, we express it through a general
\emph{data-consistency potential} $\Phi(\,\cdot\,;d_{\mathrm{obs}})$(i.e., a negative log-likelihood up to an additive constant), following the standard Bayesian inverse problem formulation; see, e.g., \cite{stuart2010inverse,dashti2013bayesian}:
\begin{equation}\label{eq:likelihood_general}
p(d_{\mathrm{obs}}|v)\ \propto\ \exp\!\big(-\Phi(v;d_{\mathrm{obs}})\big).
\end{equation}
This formulation decouples the data-error model \eqref{eq:data_model} from the choice of data-misfit functional used to quantify agreement between simulated and observed data.
As a classical special case, the additive Gaussian noise model in \eqref{eq:data_model} yields
\begin{equation}\label{eq:likelihood}
p(d_{\mathrm{obs}}| v)\ \propto\
\exp\left(-\tfrac12\big\|\mathcal{F}(v)-d_{\mathrm{obs}}\big\|_{\Sigma^{-1}}^{2}\right),\quad
\|r\|_{\Sigma^{-1}}^{2}:=r^{\top}\Sigma^{-1}r.
\end{equation}
corresponding to $\Phi(v;d_{\mathrm{obs}})=\tfrac12\|\mathcal{F}(v)-d_{\mathrm{obs}}\|_{\Sigma^{-1}}^{2}$.

Throughout the paper, the observation $d_{\mathrm{obs}}$ is fixed. We often suppress this dependence and write $\Phi(v)\equiv \Phi(v;d_{\mathrm{obs}})$.
In what follows, we keep $\Phi$ general until Section~\ref{sec:method}, where we instantiate it by a robust OT-based waveform misfit tailored to amplitude dominance, cycle-skipping, and objective scaling.

\textit{Prior.}
To regularize the inverse problem and encode admissible structural assumptions on the subsurface,
we endow the wave speed $v$ with a prior distribution $\pi_0(\mathrm{d}v)$ (with density $p(v)$ in the finite-dimensional discretization).
Such prior information is crucial in FWI, which is typically ill-posed due to limited acquisition geometry and band-limited data,
and hence may admit non-unique and unstable reconstructions under purely data-driven fitting.

A common choice is a Gibbs-type prior of the form
\begin{equation*}
p(v) \propto \exp\big(-\mu\mathcal{R}(v)\big),
\end{equation*}
where $\mathcal{R}(v)$ is a regularization functional and $\mu>0$ controls its strength.
In this case, the negative log-prior contributes the term $\mu\,\mathcal{R}(v)$ to the posterior objective,
so that classical deterministic regularization is recovered at the level of maximum a posteriori (MAP) estimation.

\textit{Posterior.}
Combining the likelihood~\eqref{eq:likelihood_general} with the prior $p(v)$ via Bayes' rule gives the posterior distribution
\begin{equation}\label{eq:posterior}
p(v | d_{\mathrm{obs}})\propto p(d_{\mathrm{obs}}| v) p(v)
\propto \exp\big(-\Phi(v;d_{\mathrm{obs}})\big) p(v).
\end{equation}
Once $\Phi$ is specified, the posterior combines the data-fit term encoded by the likelihood with the prior information encoded by $p(v)$.

\textit{Connection to deterministic FWI (MAP).}
Classical deterministic formulations of FWI are often posed as minimizing a data-misfit functional, possibly augmented by regularization. Within the Bayesian formulation, such approaches can be interpreted as computing a MAP estimate. Taking the negative logarithm of~\eqref{eq:posterior} yields, up to an additive constant,
\begin{equation*}
-\log p(v | d_{\mathrm{obs}}) =\Phi(v;d_{\mathrm{obs}})-\log p(v)+C,
\end{equation*}
where $C$ is independent of $v$. Consequently,
\begin{equation}\label{eq:map_estimator}
v_{\mathrm{MAP}}= \arg\min_{v\in\mathcal{V}} \Phi(v;d_{\mathrm{obs}})-\log p(v),
\end{equation}
where $\mathcal{V}$ denotes the admissible set (e.g., enforcing $v>0$).
If $\Phi$ is chosen as a quadratic (Gaussian) data-misfit and $p(v)\propto\exp(-\mu\mathcal{R}(v))$,
then~\eqref{eq:map_estimator} reduces to a classical regularized PDE-constrained optimization problem for FWI.
In the remainder of the paper, we compute posterior-driven reconstructions using a diffusion-based solver with a learned prior score and a data-consistency potential $\Phi$.

% ============================================================
\section{Score-Based Generative Model}\label{sec:diffusion}
Diffusion models, also known as \emph{score-based generative models} \cite{song2020score,ho2020denoising,yang2023diffusion}, define a generative process by gradually perturbing samples with Gaussian noise and then learning to reverse this corruption process. From the score-based viewpoint, these models are grounded in the estimation of the \emph{score function}, i.e., the gradient of the log-density with respect to the variable of interest.

\subsection{Setup and notation}\label{subsec:diff_setup}
We identify the discretized velocity field $v\in\mathbb{R}^m$ ($m=m_x m_z$) with the data vector $\mathbf{x}_0$ and write $\mathbf{x}_0:=v$. This is purely a notational convention to align with diffusion-model formulations, where $\mathbf{x}_0$ denotes an uncorrupted sample at diffusion time $t=0$. We regard $\mathbf{x}_0$ as a random vector following an unknown distribution $p_{\mathrm{data}}$ on $\mathbb{R}^m$, representing the population of plausible subsurface models. In our Bayesian setting, $p_{\mathrm{data}}$ plays the role of the prior $p(v)$, from which we assume access to samples.

\subsection{Forward diffusion process}\label{subsec:forward}
A diffusion model begins by defining a continuum of progressively noised variables $\{\mathbf{x}(t)\}_{t\in[0,T_{\rm{diff}}]}$ through an It\^{o} stochastic differential equation (SDE) \cite{song2021scorebased}
\begin{equation*}
\mathrm{d}\mathbf{x}(t) = f(\mathbf{x}(t),t)\,\mathrm{d}t + g(t)\,\mathrm{d}\mathbf{w}_t,
\qquad \mathbf{x}(0)=\mathbf{x}_0\sim p_{\mathrm{data}},
\end{equation*}
where $\mathbf{w}_t$ is a standard $m$-dimensional Wiener process, $f$ is the drift, and $g$ is the diffusion coefficient. Under mild regularity assumptions \cite{oksendal2003stochastic}, the law of $\mathbf{x}(t)$ admits a density $p_t$, and $\{p_t\}_{t\in[0,T_{\rm{diff}}]}$ evolves according to the Kolmogorov forward (Fokker--Planck) equation \cite{song2021scorebased}
\begin{equation*}
\partial_t p_t(\mathbf{x}) = -\nabla \cdot \big(f(\mathbf{x},t) p_t(\mathbf{x})\big) + \tfrac12 g^2(t)\Delta p_t(\mathbf{x}),
\quad p_{t=0}=p_{\mathrm{data}}.
\end{equation*}
In particular, the forward SDE induces a family of distributions $\{p_t\}$ obtained by Gaussian perturbations of $p_{\mathrm{data}}$.

\subsection{Variance-preserving diffusion and transition kernel}\label{subsec:vp}
We adopt the \emph{variance-preserving} (VP) SDE \cite{song2020score}, defined by the linear It\^{o} SDE
\begin{equation}\label{eq:vp_sde}
\mathrm{d}\mathbf{x}(t) =-\frac{1}{2} \beta(t)\,\mathbf{x}(t)\,\mathrm{d}t +\sqrt{\beta(t)}\,\mathrm{d}\mathbf{w}_t,
\quad t\in[0,T_{\rm{diff}}],
\end{equation}
where $\beta(t)>0$ is a prescribed noise schedule. In our setting, $\mathbf{x}(0)=\mathbf{x}_0$ represents a clean (vectorized) velocity model, while $\mathbf{x}(t)$ denotes its progressively perturbed version.

Since \eqref{eq:vp_sde} is linear with time-dependent coefficients, it admits an explicit mild solution. Define the integrating factor
\begin{equation*}
\alpha(t):= \exp\!\Big(-\tfrac12\int_0^t \beta(\tau)\,\mathrm{d}\tau\Big).
\end{equation*}
Applying It\^{o}'s formula to $\alpha(t)^{-1}\mathbf{x}(t)$ yields
\begin{equation}\label{eq:vp_solution}
\begin{aligned}
\mathbf{x}(t)
&= \alpha(t)\mathbf{x}_0
  + \alpha(t)\int_0^t \alpha(s)^{-1}\sqrt{\beta(s)}\,\mathrm{d}\mathbf{w}_s \\
&= \alpha(t)\mathbf{x}_0
  + \int_0^t
    \exp\!\Big(-\tfrac12\int_s^t \beta(\tau)\,\mathrm{d}\tau\Big)
    \sqrt{\beta(s)}\,\mathrm{d}\mathbf{w}_s .
\end{aligned}
\end{equation}
The stochastic integral in \eqref{eq:vp_solution} is Gaussian with mean zero. Setting
\begin{equation*}
\hat \sigma^2(t) := 1-\exp\!\Big(-\int_0^t \beta(\tau)\,\mathrm{d}\tau\Big)
= 1-\alpha^2(t),
\end{equation*}
we obtain the Gaussian transition kernel
\begin{equation}\label{eq:vp_transition}
\mathbf{x}(t) | \mathbf{x}_0\sim\mathcal{N}\!\big(\alpha(t)\,\mathbf{x}_0,\ \hat \sigma^2(t)\,I\big),
\end{equation}
or, equivalently, the reparameterization
\begin{equation*}
\mathbf{x}(t)=\alpha(t)\,\mathbf{x}_0+\hat \sigma(t)\,\boldsymbol{z},
\quad \boldsymbol{z}\sim\mathcal{N}(0,I).
\end{equation*}
If the schedule is chosen so that $\int_0^{T_{\rm{diff}}} \beta(\tau) \mathrm{d}\tau$ is sufficiently large, then $\alpha(T_{\rm{diff}})\approx 0$ and $\hat \sigma^2(T_{\rm{diff}})\approx 1$, implying that $\mathbf{x}(T_{\rm{diff}})$ is approximately standard Gaussian and essentially independent of $\mathbf{x}_0$. In practice we sample $\mathbf{x}(T_{\rm{diff}})\sim \mathcal{N}(0,I)$ to initialize the reverse-time procedure.

\subsection{Reverse-time dynamics and score learning}\label{subsec:reverse}
Let $q_t(\cdot|\mathbf{x}_0)$ denote the transition kernel of the Markov diffusion \eqref{eq:vp_sde}, and define, for each $t\in[0,T_{\rm{diff}}]$, the time-$t$ density
\begin{equation*}
p_t(\mathbf{x}):=\int_{\mathbb{R}^m}q_t(\mathbf{x} | \mathbf{x}_0)p_{\mathrm{data}}(\mathbf{x}_0)\mathrm{d}\mathbf{x}_0,
\end{equation*}
whenever the density exists. The main object of interest is the score field
\begin{equation*}
s^\star(\mathbf{x},t):=\nabla_{\mathbf{x}} \log p_t(\mathbf{x}),
\end{equation*}
i.e., the logarithmic derivative of the marginal density induced by the forward diffusion at time $t$ (equivalently, at the corresponding noise level).

Under mild regularity assumptions, a time-reversal formula for diffusion processes \cite{anderson1982reverse} implies that the reverse-time dynamics associated with \eqref{eq:vp_sde} can be expressed in terms of the score. In the VP case, where $f(\mathbf{x},t)=-\tfrac12\beta(t)\mathbf{x}$ and $g^2(t)=\beta(t)$, the reverse-time SDE can be written as
\begin{equation}\label{diff:backward}
\mathrm{d}\mathbf{x}(t)=\Big(-\tfrac12\beta(t)\mathbf{x}(t)-\beta(t)s^\star(\mathbf{x}(t),t)\Big)\mathrm{d}t
+\sqrt{\beta(t)}\,\mathrm{d}\bar{\mathbf{w}}_t,
\qquad t\in[0,T_{\rm diff}],
\end{equation}
interpreted in reverse time, i.e., the process is run backward from $t=T_{\rm diff}$ to $t=0$, where $\bar{\mathbf{w}}_t$ is a Wiener process in the reverse-time parameterization.

We approximate the true score $s^\star(\mathbf{x},t)=\nabla_{\mathbf{x}}\log p_t(\mathbf{x})$ by a neural network
$s_\theta(\mathbf{x},t)$ (the \emph{score network}), parameterized by $\theta$.
A natural training objective is the mean-squared error regression:
\begin{equation}\label{eq:L1_rewrite}
\mathcal{L}_1(\theta):=\mathbb{E}_{t\sim \mathcal{U}(0,T_{\rm{diff}})}\mathbb{E}_{\mathbf{x}\sim p_t}\Big[
\kappa(t)\|s_\theta(\mathbf{x},t)-\nabla_{\mathbf{x}}\log p_t(\mathbf{x})\|_2^2
\Big],
\end{equation}
where $\kappa(t)>0$ reweights the contributions of different noise levels.
However, \eqref{eq:L1_rewrite} is not directly computable in practice, since the marginal density $p_t$ is induced by the unknown data distribution $p_{\mathrm{data}}$.

Since the marginal score $\nabla_{\mathbf{x}}\log p_t(\mathbf{x})$ is intractable, we adopt denoising score matching (DSM) \cite{vincent2011connection}.
DSM exploits the fact that the forward transition kernel $q_t(\cdot\mid \mathbf{x}_0)$ is known, and trains the score network to match the corresponding conditional score. Concretely, it minimizes
\begin{equation}\label{eq:L2_rewrite}
\mathcal{L}_2(\theta):=\mathbb{E}_{t\sim \mathcal{U}(0,T_{\rm{diff}})}\mathbb{E}_{\mathbf{x}_0\sim p_{\mathrm{data}}}
\mathbb{E}_{\mathbf{x}\sim q_t(\cdot\mid\mathbf{x}_0)}\Big[
\kappa(t)\,\|s_\theta(\mathbf{x},t)-\nabla_{\mathbf{x}} \log q_t(\mathbf{x} | \mathbf{x}_0)\|_2^2\Big].
\end{equation}
The mathematical rationale is the identity
\begin{equation*}
\nabla_{\mathbf{x}} \log p_t(\mathbf{x})
=\mathbb{E}_{\mathbf{x}_0\sim p(\cdot\,\mid\mathbf{x}(t)=\mathbf{x})}
\left[\nabla_{\mathbf{x}} \log q_t(\mathbf{x} | \mathbf{x}_0)\right],
\end{equation*}
so the marginal score is the conditional expectation of the conditional score.
By the $L^2$ projection property of conditional expectation, $\mathcal{L}_2(\theta)=\mathcal{L}_1(\theta)+C$ with a constant $C\ge 0$ independent of $\theta$. Hence $\mathcal{L}_1$ and $\mathcal{L}_2$ have the same minimizers and $\mathcal{L}_2(\theta)\ge \mathcal{L}_1(\theta)$.
For the VP diffusion, $q_t(\cdot | \mathbf{x}_0)$ is Gaussian (\eqref{eq:vp_transition}), hence
\begin{equation*}
\nabla_{\mathbf{x}}\log q_t(\mathbf{x} | \mathbf{x}_0)= -\frac{\mathbf{x}-\alpha(t)\mathbf{x}_0}{\hat \sigma^2(t)}.
\end{equation*}

After training, let $\theta^\ast$ denote the learned parameters (e.g., an empirical minimizer of the training loss \eqref{eq:L2_rewrite}).
Replacing the unknown score in the reverse-time dynamics by $s_{\theta^\ast}$ yields the approximate reverse SDE
\begin{equation*}
\mathrm{d}\mathbf{x}(t)=\Big(f(\mathbf{x}(t),t)-g^2(t)\,s_{\theta^\ast}(\mathbf{x}(t),t)\Big)\mathrm{d}t
+g(t)\,\mathrm{d}\bar{\mathbf{w}}_t,\qquad t\in[0,T_{\rm diff}],
\end{equation*}
interpreted in reverse time (run backward from $t=T_{\rm diff}$ to $t=0$), with $f(\mathbf{x},t)=-\tfrac12\beta(t)\mathbf{x}$ and $g(t)=\sqrt{\beta(t)}$ for the VP diffusion.

% ============================================================
\section{Conditional Diffusion for Bayesian Inversion}\label{sec:baseline_dps}

\subsection{Diffusion-based posterior solver}\label{subsec:dps_solver}
We now describe a diffusion-based solver for computing posterior-driven reconstructions of $v$ from the Bayesian model \eqref{eq:posterior}.
Following the diffusion posterior sampling (DPS) framework, we incorporate data consistency through a likelihood potential and use a learned score model to represent prior information.
For notational consistency with the DPS literature, we denote the observation by $y^\delta:=d_{\mathrm{obs}}$ and recall the reverse-time SDE \eqref{diff:backward} for the VP diffusion.
\emph{Throughout this section, the observation $y^\delta$ is fixed; hence we suppress this dependence in the notation and write $\Phi(\cdot)\equiv \Phi(\cdot; y^\delta)$.}
To draw samples from the Bayesian posterior distribution $p(v | y^\delta)$ in \eqref{eq:posterior},
we replace the marginal score $s^\star(\mathbf{x},t)=\nabla_{\mathbf{x}}\log p_t(\mathbf{x})$
by the conditional score $\nabla_{\mathbf{x}}\log p_t(\mathbf{x} | y^\delta)$, leading to
\begin{equation*}
\mathrm{d}\mathbf{x}(t)=\Big(-\tfrac{1}{2}\beta(t)\mathbf{x}(t)-\beta(t)\nabla_{\mathbf{x}}\log p_t(\mathbf{x}(t)| y^\delta)\Big)\mathrm{d}t
+\sqrt{\beta(t)}\,\mathrm{d}\bar{\mathbf{w}}_t,
\qquad t\in[0,T_{\rm diff}],
\end{equation*}
interpreted in reverse time (run backward from $t=T_{\rm diff}$ to $t=0$).
Using Bayes' rule, the conditional score decomposes as
\begin{equation}\label{eq:condi_score}
\nabla_{\mathbf{x}(t)}\log p_t(\mathbf{x}(t) | y^\delta)
=\nabla_{\mathbf{x}(t)}\log p_t(\mathbf{x}(t))+\nabla_{\mathbf{x}(t)}\log p_t(y^\delta | \mathbf{x}(t)),
\end{equation}
so that, given an approximation of the marginal score $\nabla\log p_t(\mathbf{x}(t))\approx s_{\theta^*}(\mathbf{x}(t),t)$,
it remains to approximate the term $\nabla_{\mathbf{x}(t)}\log p_t(y^\delta|\mathbf{x}(t))$.

\textit{From likelihood to a general potential.}
Under a Gaussian likelihood \eqref{eq:likelihood}, a common approximation is
\begin{equation*}
\nabla_{\mathbf{x}(t)}\log p_t(y^{\delta} | \mathbf{x}(t))
\simeq -\frac{1}{2}\nabla_{\mathbf{x}(t)}\|y^{\delta}-\mathcal{F}(\hat{\mathbf{x}}_0(t))\|_{\Sigma^{-1}}^2,
% \label{eq:appr1}
\end{equation*}
where $\hat{\mathbf{x}}_0(t)$ denotes the denoised (posterior-mean) estimate of the clean variable associated with $\mathbf{x}(t)$.
Following Eq.~(10) in \cite{chung2022diffusion}, we set
\begin{equation*}
\hat{\mathbf{x}}_0(t)
:= \frac{1}{\sqrt{\bar{\alpha}(t)}}\Big(\mathbf{x}(t) + \big(1-\bar{\alpha}(t)\big)s_{\theta^*}(\mathbf{x}(t),t)\Big),
\end{equation*}
where $s_{\theta^*}(\mathbf{x},t)$ is the trained score network and, for the VP diffusion,
$\bar{\alpha}(t):=\exp\big(-\int_0^t \beta(\tau)d\tau\big)$.

In our formulation \eqref{eq:likelihood_general}, data consistency is encoded through a potential $\Phi$.
Accordingly, we approximate the conditioning term by
\begin{equation}\label{eq:appr1_general}
\nabla_{\mathbf{x}(t)}\log p_t(y^{\delta}|\mathbf{x}(t)) \simeq
-\nabla_{\mathbf{x}(t)}\Phi(\hat{\mathbf{x}}_0(t)),
\end{equation}
so that the observational constraint enters the reverse-time dynamics via the chosen potential~$\Phi$.

Combining~\eqref{eq:condi_score} with~\eqref{eq:appr1_general} and replacing the prior score by the trained score network yields
\begin{equation}\label{eq:grad_appr}
\nabla_{\mathbf{x}(t)}\log p_t(\mathbf{x}(t)|y^{\delta})
\simeq
s_{\theta^*}(\mathbf{x}(t),t)-\zeta\,\nabla_{\mathbf{x}(t)}\Phi(\hat{\mathbf{x}}_0(t)),
\end{equation}
where $\zeta>0$ is a user-chosen guidance-strength parameter that rescales the data-consistency term to account for objective scaling.

\subsection{Baseline DPS discretization for FWI}\label{subsec:baseline_ddpm}
Let $\{\beta_i\}_{i=1}^N$ be the discrete noise schedule, $\alpha_i:=1-\beta_i$, and $\bar\alpha_i:=\prod_{j=1}^i \alpha_j$.
Denote by $v_i$ the reverse-time diffusion state at step $i$ (DDPM indexing \cite{ho2020denoising}), so that $v_i\approx \mathbf{x}(t_i)$.
Here $i=N$ corresponds to the highest noise level (an approximately Gaussian state), and $i=1$ corresponds to the lowest noise level; thus
$\{v_i\}_{i=1}^N$ traces a discrete reverse-time trajectory that progressively denoises the state from $v_N$ toward a clean reconstruction.
Denote by $s_{\theta^\ast}(\cdot,i)$ the trained score network.
The standard clean estimator is
\begin{equation}\label{eq:v0hat}
\hat v_0^{(i)}=\hat v_0(v_i,i)
:=\frac{1}{\sqrt{\bar\alpha_i}}
\Big(v_i + (1-\bar\alpha_i)s_{\theta^\ast}(v_i,i)\Big).
\end{equation}

\textit{Baseline likelihood potential (least squares).}
In the original DPS formulation \cite{chung2022diffusion}, data consistency is enforced through a Gaussian likelihood,
leading to the least-squares potential
\begin{equation*}
\Phi_{\ell_2}(v):=\tfrac12\|\mathcal{F}(v)-y^\delta\|_{\Sigma^{-1}}^2.
\end{equation*}

Given the prior-driven proposal $v_{i-1}'$, baseline DPS applies a scalar guidance correction evaluated at the current state $v_i$:
\begin{align}\label{eq:dps_baseline_step}
 v_{i-1}&=v_{i-1}'-\zeta\nabla_{v_i}\Phi_{\ell_2}(\hat v_0^{(i)}),\\ \nonumber
 &=v_{i-1}'-\zeta\big(\mathcal{F}^{'}(\hat v_0^{(i)})\big)^\ast
\Sigma^{-1}\big(\mathcal{F}(\hat v_0^{(i)})-y^\delta\big),\qquad \zeta>0.
\end{align}

Based on the above formulation, the DPS-based procedure for solving the full waveform inversion problem is summarized in Algorithm~\ref{al:base_dps}.
\begin{algorithm}[t]
\textbf{Require}: $N$, schedule $\{\beta_i\}_{i=1}^N$, score network $s_{\theta^\ast}$, variance $\{\hat \sigma_i\}$ fixed observation $d_{\mathrm{obs}}$,
hyperparameters $\{\zeta_i\}$ (and optional clipping bounds), objective function $\Phi_{l_2}$.\\
\textbf{Result}: obtain the inversion velocity model $v=v_0$.\\
Initialize $v_N\sim\mathcal{N}(0,I)$; set $\alpha_i:=1-\beta_i$ and $\bar\alpha_i:=\prod_{j=1}^i\alpha_j$.\\
\For{$i=N$ \KwTo $1$}{
  $\hat v_0^{(i)} \leftarrow \frac{1}{\sqrt{\bar\alpha_i}}\big(v_i+(1-\bar\alpha_i)s_{\theta^\ast}(v_i,i)\big)$;\\
  $\mathbf{z}\sim \mathcal{N}(0,\mathbf{I})$;\\
  $v'_{i-1}\leftarrow\frac{\sqrt{\alpha_i}(1-\bar \alpha_{i-1})}{1-\bar \alpha_i}v_i+\frac{\sqrt{\bar \alpha_{i-1}}\beta_i}{1-\bar \alpha_i}\hat v_0^{(i)}+\hat \sigma_i\mathbf{z}$;\\
  $v_{i-1}\leftarrow v_{i-1}'-\zeta\,\nabla_{v_i} \Phi_{l_2}(\hat{v}_0)$;\\
}
\caption{Baseline Diffusion Posterior Sampling (DPS)}
\label{al:base_dps}
\end{algorithm}

\subsection{Limitations of $\ell_2$-guided DPS in FWI}\label{limits_dps}
Although the Algorithm~\ref{al:base_dps} is consistent with a Gaussian noise model, its direct application to FWI presents several well-known and practically important difficulties. Specifically, the pointwise $\ell_2$ potential may induce cycle-skipping and objective imbalance, while the use of a single global scalar guidance parameter is inadequate for the strongly space-dependent sensitivities characteristic of FWI.

\textit{(i) Cycle-skipping and nonconvexity.}
The $\ell_2$ misfit compares waveforms pointwise in time; small phase shifts can yield large residuals and highly oscillatory gradients.
As a consequence, $\Phi_{\ell_2}$ typically exhibits a strongly nonconvex landscape with many spurious local minima, which can lead to
cycle-skipping in solving the present inverse problem. Within DPS, the guidance term is applied to the denoised estimate $\hat v_0^{(i)}$, and early
reverse-time iterates often produce rough $\hat v_0^{(i)}$; in this regime, the $\ell_2$ gradient can be particularly unreliable,
causing unstable or ineffective guidance.

\textit{(ii) Amplitude dominance and objective scaling.}
Seismic data often exhibit a pronounced dynamic-range imbalance across arrivals. Since $\Phi_{\ell_2}$ is quadratic in the residual,
the gradient $\nabla\Phi_{\ell_2}$ is dominated by the largest residual components, which are typically associated with high-amplitude early arrivals.
As a result, the baseline correction step tends to prioritize fitting the strongest events, while weaker phases that are informative for deeper or poorly
illuminated regions are comparatively down-weighted.
Moreover, $\Phi_{\ell_2}(v)=\tfrac12\|\mathcal{F}(v)-y^\delta\|_{\Sigma^{-1}}^2$ is used without normalization by an observation-dependent scale; when the data (or residuals
after preprocessing) have small magnitude, both $\Phi_{\ell_2}$ and $\nabla\Phi_{\ell_2}$ can become numerically small, making the guidance strength $\zeta$ in
\eqref{eq:dps_baseline_step} delicate to tune and potentially leading to overly weak data-consistency enforcement.

\textit{(iii) Heterogeneous sensitivity and scalar step size.}
FWI sensitivities are strongly space-dependent (e.g., due to limited illumination), so a single global scalar guidance parameter $\zeta$ cannot balance updates across
well-illuminated and poorly-illuminated regions. This mismatch can slow down progress in weak-sensitivity areas and may yield unstable behavior in sensitive regions.

\smallskip
\noindent These observations motivate two modifications: we (a) replace the pointwise $\ell_2$ potential by a robust OT-based
data-consistency potential to mitigate cycle-skipping and amplitude imbalance, and (b) introduce a variable-metric preconditioned
guidance scheme to stabilize and balance spatial updates. We present these components next.
% ============================================================
\section{Our method: OT-guided preconditioned DPS for FWI}\label{sec:method}

This section presents two modifications motivated by the limitations of $\ell_2$-guided DPS in FWI discussed in the above section: (a) an OT-based data-consistency potential $\Phi=\mathcal{J}$ that mitigates amplitude dominance, reduces phase-sensitivity (cycle-skipping), and stabilizes objective scaling; and (b) a variable-metric preconditioned guidance scheme that adapts the guidance strength across noise levels and balances spatial updates under heterogeneous sensitivities.

\subsection{OT-based data-consistency potential}\label{subsec:ot_misfit}
As discussed in Section~\ref{limits_dps}, pointwise least-squares waveform fitting can lead to amplitude-dominated updates,
strong phase sensitivity (cycle-skipping), and ill-conditioned objective scaling. We therefore construct a robust OT-based misfit
$\mathcal{J}$ in three steps, each targeting one of these difficulties.

\textit{(I) Amplitude-adaptive data weighting.}
Let $d_{\mathrm{obs}}(s,r,t)$ denote the observed trace associated with the source--receiver pair $(s,r)$, and let
$d_{\mathrm{syn}}(s,r,t;v)$ be its synthetic counterpart generated by the velocity model $v$. We introduce a globally normalized
instantaneous amplitude
\begin{equation*}
a(s,r,t)=\frac{|d_{\mathrm{obs}}(s,r,t)|}{\max_{i,j,t}|d_{\mathrm{obs}}(s_i,r_j,t)|}\in[0,1],
\end{equation*}
and define the bounded weight
\begin{equation}\label{eq:omega_def}
\omega(s,r,t)=\frac{1}{1+k\,a(s,r,t)}, \quad k\ge 0.
\end{equation}
By construction, $\omega(s,r,t)\in\big[1/(1+k),\,1\big]$ and is monotone decreasing in $a(s,r,t)$, so large-amplitude arrivals are
down-weighted while low-amplitude portions are largely preserved. We apply the same pointwise weighting to both observed and synthetic
signals,
\begin{equation}\label{eq:weighted_traces}
\tilde d_{\mathrm{obs}}(s,r,t)=\omega(s,r,t)\,d_{\mathrm{obs}}(s,r,t),
\quad\tilde d_{\mathrm{syn}}(s,r,t;v)=\omega(s,r,t)\,d_{\mathrm{syn}}(s,r,t;v),
\end{equation}
so that the comparison is performed under a fixed, observation-dependent weighting. Equivalently,
$\tilde d_{\mathrm{syn}}-\tilde d_{\mathrm{obs}}=\omega\,(d_{\mathrm{syn}}-d_{\mathrm{obs}})$, and thus the contribution of
high-amplitude events to the misfit and its gradient is reduced. Importantly, $\omega$ is computed from the observations and kept fixed during inversion, so the weighting does not introduce additional model-dependent nonlinearity. This amplitude-adaptive weighting alleviates scale disparities induced by raw amplitudes and prevents the descent direction from being dominated by the high-amplitude portions of the data.

\begin{figure}[h]
    \centering
    \begin{subfigure}[b]{0.45\textwidth}
        \centering
        \includegraphics[width=\textwidth]{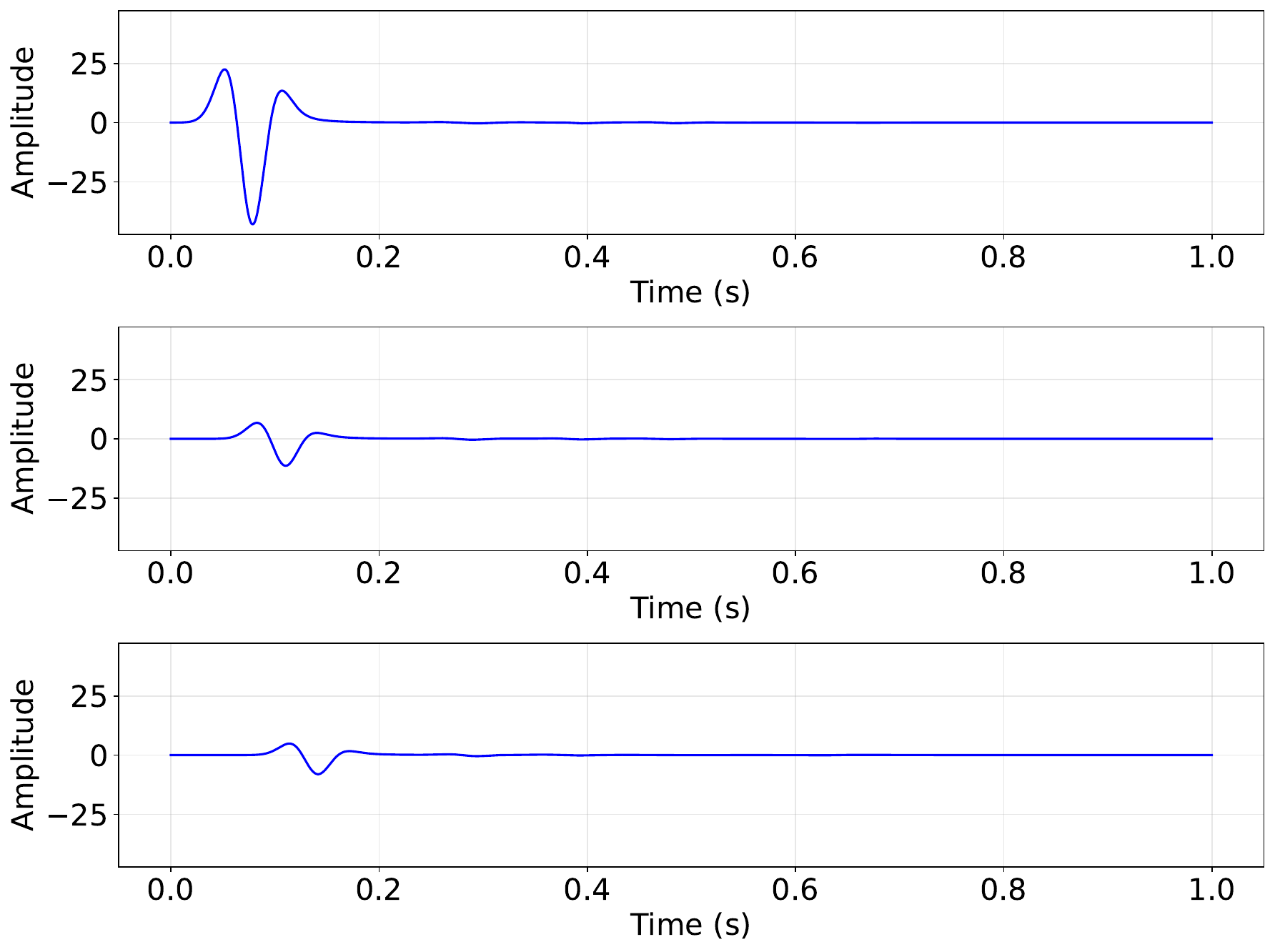}
        \caption{The original observed wavefield.}
        \label{fig:origin_wave}
    \end{subfigure}
    \hfill
    \begin{subfigure}[b]{0.45\textwidth}
        \centering
        \includegraphics[width=\textwidth]{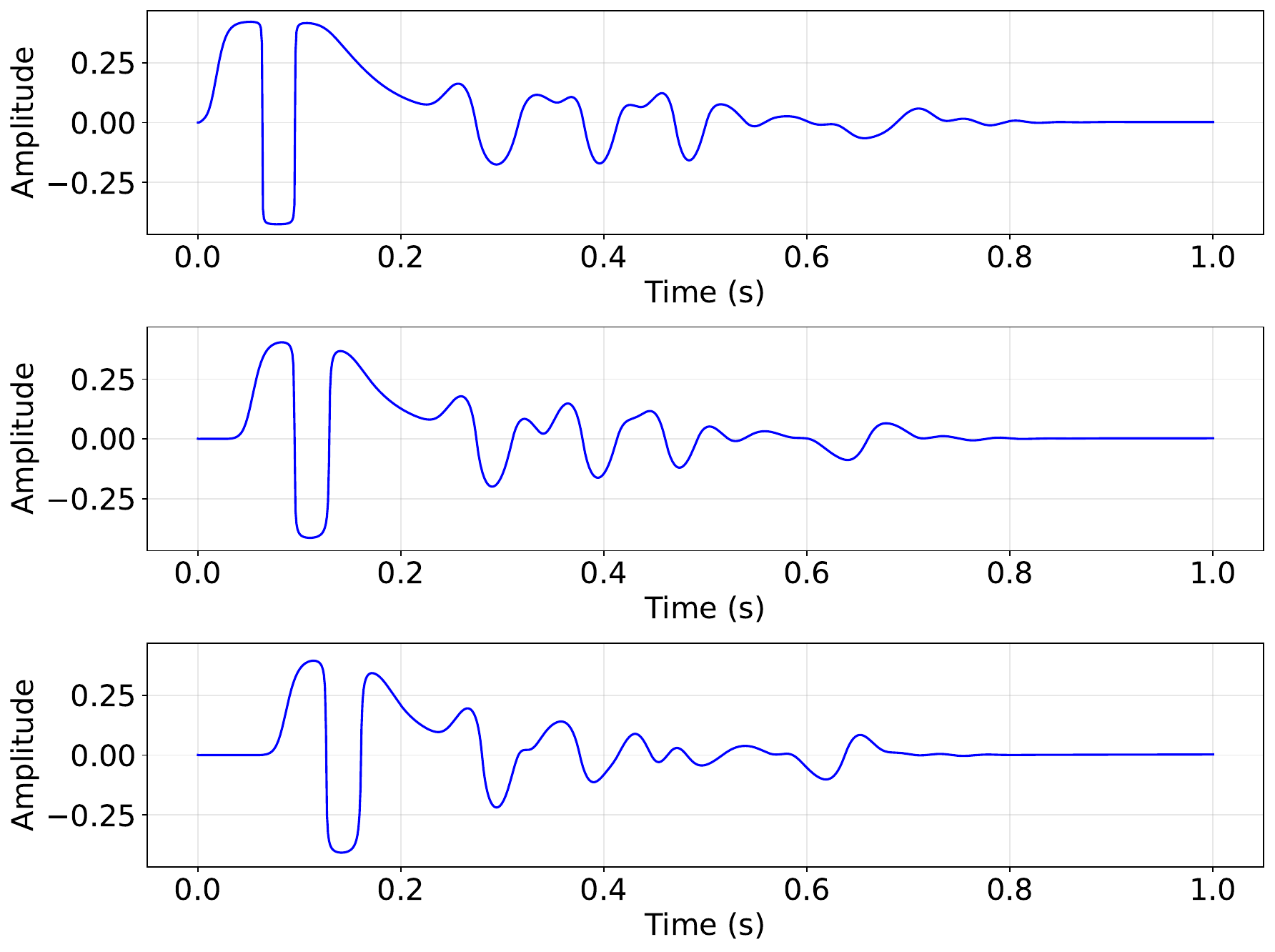}
        \caption{The weighted wavefield with $k=100$.}
        \label{fig:weighted_wave}
    \end{subfigure}
    \caption{Comparison of the original and weighted wavefield traces at Receiver~1 for the first three source--receiver pairs. 
    The original recordings are dominated by early large-amplitude arrivals, whereas the proposed bounded weighting compresses the dynamic range and renders later weak-amplitude events more comparable in scale.}
    \label{fig:wave}
\end{figure}

As an empirical illustration, we take CurveFault-A dataset \cite{deng2022openfwi} as an example and compare the original traces and their weighted counterparts at Receiver~1 for the first three source--receiver pairs in Figure~\ref{fig:wave}. 
In the original recordings (Figure~\ref{fig:origin_wave}), the early-arriving portions exhibit substantially larger amplitudes than later-arriving events, thereby overwhelming the contribution of weaker phases in the least-squares misfit.
After applying the bounded weighting \eqref{eq:omega_def} with $k=100$ (Figure~\ref{fig:weighted_wave}), the dynamic range is effectively compressed: high-amplitude arrivals are suppressed while low-amplitude events are largely preserved. This amplitude balancing reduces the tendency of large-amplitude portions of the data to dominate the gradient and yields a more balanced weighted comparison for subsequent inversion.

Nevertheless, even with amplitude balancing, pointwise $\ell_2$ misfits can remain highly sensitive to small time shifts and phase
mismatch. To alleviate the resulting nonlinearity and reduce cycle-skipping, we therefore adopt the Wasserstein-$p$ distance as the data-misfit measure (we use $p=2$ in all computations).

\textit{(II) Transport-based misfit (Wasserstein-$p$).}
The least-squares misfit induces a pointwise $\ell_2$ geometry that is overly sensitive to temporal misalignment, which can lead to a highly nonconvex objective landscape and promote cycle-skipping. We therefore adopt the Wasserstein-$p$ distance, which equips the data space with an optimal-transport geometry and measures discrepancy via mass displacement (after suitable normalization), rather than by pointwise differences.

Since the Wasserstein distance is formally defined on nonnegative measures with equal total mass, we follow the standard OT-FWI
construction~\cite{yang2018application} and transform each enhanced 1D trace into a probability density function (PDF) via a fixed
shift--rescale normalization. Specifically, for a trace $\tilde d(t)$ (either $\tilde d_{\mathrm{obs}}$ or $\tilde d_{\mathrm{syn}}$),
we introduce a constant $c'>0$ to ensure nonnegativity. To guarantee that $\tilde d(t)+c'\ge 0$ on $[0,T]$ for all traces, we fix $(s,r)$ and take
\begin{equation*}
c' = 1.1\Bigl|\min_{t}\tilde d_{\mathrm{obs}}(s,r,t)\Bigr|,
\end{equation*}
and keep $c'$ constant throughout the inversion, so that the mapping does not introduce additional model-dependent nonlinearity in the
gradient. In practice, $c'$ is chosen sufficiently large so that $\tilde d_{\mathrm{syn}}+c'\ge 0$ throughout the inversion. We then define
the density mapping
\begin{equation}
P(\tilde d)(t)=
\frac{\tilde d(t)+c'}{\int_0^T \bigl(\tilde d(\tau)+c'\bigr)d\tau},
\label{eq:pdf_map}
\end{equation}
which ensures $P(\tilde d)(t)\ge 0$ and $\int_0^T P(\tilde d)(t)dt=1$.
Accordingly, for each source--receiver pair $(s,r)$, the associated probability densities are
\begin{equation*}
% \label{eq:rho_def}
\tilde\rho_{\mathrm{obs}}(\cdot;s,r)=P\!\left(\tilde d_{\mathrm{obs}}(\cdot;s,r)\right),\quad
\tilde\rho_{\mathrm{syn}}(\cdot;s,r,v)=P\!\left(\tilde d_{\mathrm{syn}}(\cdot;s,r,v)\right).
\end{equation*}

With the signals transformed into densities, the Wasserstein-$p$ distance is defined by
\begin{equation*}
W_p\left(\tilde\rho_{\mathrm{syn}},\tilde\rho_{\mathrm{obs}}\right)=
\left(\inf_{\gamma\in\tilde{\Gamma}(\tilde\rho_{\mathrm{syn}},\tilde\rho_{\mathrm{obs}})}
\int_{[0,T]\times[0,T]} |t-\tau|^p\, d\gamma(t,\tau)\right)^{1/p},\quad p\ge 1,
\end{equation*}
where $\tilde{\Gamma}(\tilde\rho_{\mathrm{syn}},\tilde\rho_{\mathrm{obs}})$ denotes the set of admissible transport plans with marginals
$\tilde\rho_{\mathrm{syn}}$ and $\tilde\rho_{\mathrm{obs}}$. In the one-dimensional setting, $W_p$ admits an explicit representation
in terms of quantile functions. Let $\tilde F_{\mathrm{syn}}$ and $\tilde F_{\mathrm{obs}}$ be the cumulative distribution functions (CDFs) of
$\tilde\rho_{\mathrm{syn}}$ and $\tilde\rho_{\mathrm{obs}}$, respectively; then
\begin{equation}\label{eq:Wp_1d}
W_p\left(\tilde\rho_{\mathrm{syn}},\tilde\rho_{\mathrm{obs}}\right)=
\left(\int_0^1\left|\tilde F_{\mathrm{syn}}^{-1}(\xi)-\tilde F_{\mathrm{obs}}^{-1}(\xi)\right|^p\, d\xi\right)^{1/p},
\end{equation}
where $\tilde F^{-1}$ denotes the quantile function, defined as the generalized inverse
$\tilde F^{-1}(\xi):=\inf\{t\in[0,T]: \tilde F(t)\ge \xi\}$. While we present the definition for general $p\ge 1$, all subsequent results and
computations in this paper use $p=2$. Numerically, for each trace we compute the PDF using \eqref{eq:pdf_map}, construct the CDF via
cumulative summation, and evaluate the quantile difference in \eqref{eq:Wp_1d} by interpolation.

\textit{(III) Misfit scale normalization.}
To stabilize the numerical scale of the OT guidance term, we rescale the aggregated discrepancy by an observation-dependent
characteristic magnitude. Specifically, let $\tilde F_{\mathrm{obs},s,r}^{-1}$ denote the quantile function associated with the observed PDF
for the pair $(s,r)$, and define the scale factor
\begin{equation*}
S_{\mathrm{obs}}=\sum_{s}\sum_{r}\left(\int_0^1\lvert \tilde F_{\mathrm{obs},s,r}^{-1}(\xi)\rvert^p\, d\xi\right)^{1/p},
% \label{eq:scale_factor}
\end{equation*}
which depends only on observations and is independent of the velocity model $v$.
Therefore, this standardization does not alter the minimizer of the objective; it merely rescales the functional and its gradient by a constant factor, improving numerical conditioning, and reducing sensitivity to step-size tuning in the guidance updates.

\textit{Weighted and normalized Wasserstein misfit.}
Combining the above ingredients, we define the normalized Wasserstein-based objective functional as

\begin{equation}\label{eq:final_objective}
\mathcal J(v)=
\frac{\displaystyle\sum_{s}\sum_{r}
W_p\left(\tilde{\rho}_{\mathrm{syn}}(\cdot;s,r;v),\tilde{\rho}_{\mathrm{obs}}(\cdot;s,r)\right)}
{\displaystyle \sum_{s}\sum_{r}\left(\int_0^1\lvert \tilde F_{\mathrm{obs},s,r}^{-1}(\xi)\rvert^p\, d\xi\right)^{1/p}},
\end{equation}
where $\tilde{\rho}_{\mathrm{syn}}$ and $\tilde{\rho}_{\mathrm{obs}}$ are the PDFs obtained from the enhanced wavefields
$\tilde d_{\mathrm{syn}}$ and $\tilde d_{\mathrm{obs}}$, respectively. 

This construction of $\mathcal{J}$ is designed to (i) reduce the dominance of high-amplitude arrivals through bounded weighting,
(ii) alleviate phase sensitivity via the OT geometry, and (iii) stabilize the numerical scale through observation-dependent normalization.

\textit{Choice of likelihood potential.}
With the above construction, we instantiate the data-consistency potential in the Bayesian likelihood \eqref{eq:likelihood_general} as
\begin{equation*}
\Phi(v;d_{\mathrm{obs}}) := \mathcal{J}(v),
\end{equation*}
and we often suppress the dependence on $d_{\mathrm{obs}}$ by writing $\mathcal{J}(v)\equiv \mathcal{J}(v;d_{\mathrm{obs}})$ and
$\Phi(v)\equiv \Phi(v;d_{\mathrm{obs}})$.
\subsection{Preconditioned Guidance in DPS}\label{subsec:precond_dps}
Building on the conditional-score approximation \eqref{eq:grad_appr},
we now introduce a variable-metric (diagonal) preconditioned guidance scheme tailored to FWI.

\textit{Our modification: variable-metric preconditioned guidance.}
We replace the uniform scalar step size by a positive diagonal preconditioner,
\begin{equation}\label{eq:precond_step}
 v_{i-1}=v_{i-1}'-P_i\nabla_{v_i}\Phi(\hat v_0^{(i)}), \quad P_i=\rho_i D_i,
\end{equation}
where $\rho_i>0$ is a scalar schedule and $D_i$ is a diagonal positive matrix.
When $P_i=\zeta I$, \eqref{eq:precond_step} reduces to the baseline DPS step.
The guidance gradient in \eqref{eq:precond_step} is computed with respect to $v_i$ through the clean estimator \eqref{eq:v0hat}:
\begin{equation}
\label{eq:chain_rule}
\nabla_{v_i}\Phi(\hat v_0^{(i)})
=\Big(\nabla_{v_i}\hat v_0(v_i,i)\Big)^{\top}\nabla_{\hat v_0}\Phi(\hat v_0^{(i)}).
\end{equation}
We evaluate \eqref{eq:chain_rule} in practice by automatic differentiation.

\textit{Stabilized guidance scheduling.}
Early reverse-time iterates typically produce $\hat v_0^{(i)}$ with spurious high-frequency artifacts. In this regime, the guidance direction can be unreliable, and overly aggressive updates may destabilize the reverse trajectory. We therefore adopt a conservative guidance strength when $\hat v_0^{(i)}$ is still rough, and gradually increase it as the estimate becomes smoother and more structured.

To quantify the roughness of the current denoised estimate, we introduce a TV-based indicator \cite{rudin1992nonlinear}. Specifically,
let $\mathrm{TV}(\cdot)$ denote the discrete total-variation semi-norm on the 2D grid:
\begin{equation}\label{eq:tv}
\mathrm{TV}(v)=\frac{1}{N_g}\sum_{\ell}\big(\vert(\nabla_x v)_\ell\vert+\vert(\nabla_z v)_\ell\vert\big),
\end{equation}
where $\ell$ indexes grid points, $N_g$ is the number of grid points, and $\nabla_x,\nabla_z$ are forward finite differences.
We use $\mathrm{TV}(\hat v_0^{(i)})$ as a roughness proxy because spurious oscillatory artifacts in early iterates typically manifest
themselves through large local finite differences and hence elevated TV values. At the same time, TV is less sensitive than quadratic
smoothness indicators to genuine sharp interfaces, which are common in subsurface models. Therefore, it provides a simple and
computationally efficient criterion for distinguishing unstable, artifact-dominated iterates from structurally meaningful ones.

Based on this proxy, we define the step-dependent scalar factor
\begin{equation}\label{eq:rho}
\rho_i=\rho_0\exp\left(-\frac{(\mathrm{TV}(\hat v_0^{(i)})-\tau)_+}{c}\right),
\qquad (\cdot)_+:=\max\{\cdot,0\},
\end{equation}
with hyperparameters $\rho_0>0$, $c>0$, and an optional threshold $\tau\ge 0$.
By construction, $\rho_i\in(0,\rho_0]$. Moreover, $\rho_i=\rho_0$ whenever $\mathrm{TV}(\hat v_0^{(i)})\le\tau$, and it decreases
monotonically with $\mathrm{TV}(\hat v_0^{(i)})$ once $\mathrm{TV}(\hat v_0^{(i)})>\tau$.
The parameter $c$ sets the decay scale: larger $c$ leads to slower attenuation, that is, less sensitivity to roughness.
The threshold $\tau$ defines a tolerance level for total variation and prevents excessive suppression of guidance in the presence of
genuine interfaces, such as sharp geological contrasts.

\textit{Spatial diagonal preconditioning.}
FWI gradients typically decay with depth due to limited illumination and attenuation, leading to highly nonuniform update magnitudes.
To mitigate this imbalance, we introduce a diagonal rescaling based on the current misfit gradient with respect to the clean estimate.
Let
\begin{equation*}
g^{(i)}:=\nabla_{\hat v_0}\Phi(\hat v_0^{(i)})\in\mathbb{R}^n.
\end{equation*}
We define the diagonal weights
\begin{equation}\label{eq:kappa}
\kappa^{(i)}_j=\left(\frac{\|g^{(i)}\|_\infty+\varepsilon}{|g^{(i)}_j|+\varepsilon}\right)^{\gamma},\quad
D_i:=\mathrm{diag}(\kappa^{(i)}_1,\dots,\kappa^{(i)}_n),
\end{equation}
where $\gamma\ge 0$ and $\varepsilon>0$. We build $D_i$ from
$g^{(i)}=\nabla_{\hat v_0}\Phi$ because it directly reflects the physical FWI sensitivity in the model domain; the chain rule in
\eqref{eq:chain_rule} then maps this sensitivity back to the current diffusion state.
%For numerical stability one may clip $\kappa^{(i)}_j$ to a prescribed interval $[\kappa_{\min},\kappa_{\max}]$.

\textit{Interpretation via a descent inequality.}
For each reverse-time index $i$, consider the composite functional $\tilde\Phi_i(v):=\Phi(\hat v_0^{(i)}(v))$.
Assume that $\tilde\Phi_i$ is differentiable and has a locally $L_i$-Lipschitz gradient on a neighborhood containing the relevant iterates. Then the descent lemma \cite{nesterov2013introductory} yields that, for any symmetric positive semidefinite matrix $P$ satisfying $\|P\|_2\le 1/L_i$,
\begin{equation}\label{eq:descent_lemma}
\tilde\Phi_i(v-P\nabla\tilde\Phi_i(v))
\le \tilde\Phi_i(v)-\frac{1}{2}\|\nabla\tilde\Phi_i(v)\|_{P}^2,
\quad\| a\|_{P}^2:=a^\top P a.
\end{equation}
Applying \eqref{eq:descent_lemma} pointwise with $P$ at the current iterate, this estimate provides a heuristic justification for choosing $P_i=\rho_iD_i$. Indeed, by tuning $\rho_i$ and clipping the diagonal entries of $D_i$, we control the effective spectral norm of $P_i$ and obtain a stabilized variable-metric guidance step, analogous to diagonal preconditioning in ill-conditioned inverse problems.

Combining the above ingredients, we arrive at the following final Algorithm~\ref{al:precond_dps}.
\begin{algorithm}[t]
\textbf{Require}: $N$, schedule $\{\beta_i\}_{i=1}^N$, score network $s_{\theta^\ast}$, variance $\{\hat \sigma_i\}$, fixed observation $d_{\mathrm{obs}}$,
hyperparameters $\rho_0,c,\tau,p,\varepsilon$ (and optional clipping bounds), objective function $\mathcal J$ as \eqref{eq:final_objective}.\\
\textbf{Result}: approximate posterior sample $v_0$.\\
Initialize $v_N\sim\mathcal{N}(0,I)$; set $\alpha_i:=1-\beta_i$ and $\bar\alpha_i:=\prod_{j=1}^i\alpha_j$.\\
\For{$i=N$ \KwTo $1$}{
  $\hat v_0^{(i)} \leftarrow \frac{1}{\sqrt{\bar\alpha_i}}\big(v_i+(1-\bar\alpha_i)s_{\theta^\ast}(v_i,i)\big)$;\\
  $\mathbf{z}\sim \mathcal{N}(0,\mathbf{I})$;\\
  $v'_{i-1}\leftarrow\frac{\sqrt{\alpha_i}(1-\bar \alpha_{i-1})}{1-\bar \alpha_i}v_i+\frac{\sqrt{\bar \alpha_{i-1}}\beta_i}{1-\bar \alpha_i}\hat v_0^{(i)}+\hat \sigma_i\mathbf{z}$;\\
  $g^{(i)}\leftarrow \nabla_{\hat v_0^{(i)}}\mathcal J(\hat v_0^{(i)})$;\\
  $\kappa^{(i)}_j\leftarrow ((\|g^{(i)}\|_\infty+\varepsilon)/( |g^{(i)}_j|+\varepsilon))^{p}$; form $D_i=\mathrm{diag}(\kappa^{(i)}_1,\dots,\kappa^{(i)}_d)$;\\
  $\rho_i \leftarrow \rho_0\exp\!\big(-(\mathrm{TV}(\hat v_0^{(i)})-\tau)_+/c\big)$;\\
  Set $P_i\leftarrow \rho_i D_i$ (optionally with clipping and/or spectral-norm control);\\
  $v_{i-1}\leftarrow v_{i-1}'-P_i\,\nabla_{v_i}\mathcal J(\hat v_0^{(i)})$;\\
}
\caption{OT-guided Wavefield-Enhanced Preconditioned DPS (OT-WE-PDPS)}
\label{al:precond_dps}
\end{algorithm}

\section{Experiments}\label{sec:exp}
We evaluate our method on OpenFWI~\cite{deng2022openfwi}. We train the score-based prior using velocity models from the Vel Family and Fault Family datasets. All inversion results are reported on held-out test samples that are disjoint from the training set; representative test models are shown in the first row of Figure~\ref{fig:main}.
The score neural network is a U-Net~\cite{ronneberger2015u} trained via denoising score matching under the VP diffusion described in Section~\ref{sec:diffusion}. Importantly, prior training uses only samples of the model parameter field (velocity) and does not require any wavefield simulations.

We compute posterior-driven reconstructions by running a discretized reverse-time diffusion whose drift combines the learned prior score and the data-consistency term induced by $\mathcal{F}$, implemented via our \emph{OT-guided Wavefield-Enhanced Preconditioned DPS (OT-WE-PDPS)} (Algorithm~\ref{al:precond_dps}). The likelihood potential is chosen as $\Phi=\mathcal{J}$, the normalized OT misfit of Section~\ref{subsec:ot_misfit}. The reported reconstructions correspond to unseen models that were not used in training the prior score network.

Experimental performance is evaluated using both qualitative inspection and quantitative metrics. 
As a primary accuracy measure, we report the relative $\ell^2$ reconstruction error
\begin{equation*}
e_{\ell^2}:=\frac{\|v_{\mathrm{true}}-v_{\mathrm{rec}}\|_2}{\|v_{\mathrm{true}}\|_2},
\end{equation*}
where $v_{\mathrm{true}}\in\mathbb{R}^m$ denotes the ground-truth velocity model (vectorized on the computational grid) and $v_{\mathrm{rec}}\in\mathbb{R}^m$ is the reconstruction, and $\|\cdot\|_2$ is the Euclidean norm.

In addition, to quantify perceptual fidelity and structural consistency, we report the peak signal-to-noise ratio (PSNR) and the structural similarity index measure (SSIM) computed between $v_{\mathrm{rec}}$ and $v_{\mathrm{true}}$.

\textit{Baselines and implementation details.}
% We discretize Equation \ref{eq:tv} using finite differences, obtaining the TV regularization term used in the experiments as follows:
% \begin{equation*}
%     \mathrm{TV}(x)=\frac{1}{(H-1)W}\sum_{i=1}^{H-1}\sum_{j=1}^{W}\left|x_{i+1,j}-x_{i,j}\right|
%     +\frac{1}{H(W-1)}\sum_{i=1}^{H}\sum_{j=1}^{W-1}\left|x_{i,j+1}-x_{i,j}\right|
% \end{equation*}
We compare our method against three baselines:\\
(i) \textbf{Baseline 1 ($W_2+\mathrm{TV}$) \cite{engquist2022optimal}.}We consider a classical baseline that minimizes a normalized $W_2$ waveform misfit with total-variation (TV) regularization.
Define
\begin{equation*}
\mathcal J_{W_2}^{\mathrm{raw}}(v):=\frac{\displaystyle \sum_{s}\sum_{r}
W_2^2\left({\rho}_{\mathrm{syn}}(\cdot;s,r;v),{\rho}_{\mathrm{obs}}(\cdot;s,r)\right)}{\displaystyle \sum_{s}\sum_{r}\int_0^1\lvert {F}_{\mathrm{obs},s,r}^{-1}(\xi)\rvert^2 d\xi},
\end{equation*}
where $\rho_{\mathrm{syn}}(\cdot;s,r;v)$ and $\rho_{\mathrm{obs}}(\cdot;s,r)$ are probability densities obtained
from the \emph{raw} synthetic and observed traces $d_{\mathrm{syn}}(\cdot;s,r;v)$ and $d_{\mathrm{obs}}(\cdot;s,r)$, respectively,
via the fixed shift--rescale density map (cf.\ \eqref{eq:pdf_map}); in particular, no wavefield enhancement is applied here.
Moreover, $F_{\mathrm{obs},s,r}$ denotes the CDF of $\rho_{\mathrm{obs}}(\cdot;s,r)$, and $F_{\mathrm{obs},s,r}^{-1}$ is the associated quantile function.

We then minimize the TV-regularized objective
\begin{equation*}
\min_{v}\ \mathcal J_{W_2}^{\mathrm{raw}}(v)+\alpha\ \mathrm{TV}(v),
\end{equation*}
where $\alpha>0$ is the regularization parameter and $\mathrm{TV}(\cdot)$ is defined in \eqref{eq:tv}.
Using a first-order method, a (sub)gradient descent iteration reads
\begin{equation*}
v \leftarrow v-\rho_0\Big(\nabla_v \mathcal J_{W_2}^{\mathrm{raw}}(v)+\alpha\,g_{\mathrm{TV}}(v)\Big),
\qquad g_{\mathrm{TV}}(v)\in \partial\,\mathrm{TV}(v),
\end{equation*}
with a fixed step size $\rho_0>0$.

(ii) \textbf{Baseline 2 (OT-WE$+\mathrm{TV}$).}
We consider an improved version of the above method that minimizes the wavefield-enhanced normalized OT misfit $\mathcal{J}(v)$ in \eqref{eq:final_objective}
(Section~\ref{subsec:ot_misfit}), augmented with a TV penalty:
\begin{equation*}
\min_{v}\ \mathcal{J}(v)+\alpha\,\mathrm{TV}(v),
\end{equation*}
where $\alpha>0$ is the regularization parameter and $\mathrm{TV}(\cdot)$ is defined in \eqref{eq:tv}.
We solve this problem using a first-order method.
In particular, we use a diagonal preconditioned (sub)gradient step
\begin{equation*}
v \leftarrow v - P(v)\Big(\nabla_v \mathcal{J}(v)+\alpha\,g_{\mathrm{TV}}(v)\Big),
\qquad g_{\mathrm{TV}}(v)\in \partial\,\mathrm{TV}(v),
\end{equation*}
where $P(v)=\rho(v)\,D(v)$ is a positive diagonal preconditioner updated from the current iterate. Here $\rho(v)$ and $D(v)$ are defined by \eqref{eq:rho} and \eqref{eq:kappa}, respectively, after replacing $\hat v_0^{(i)}$ by $v$.

(iii) \textbf{Baseline 3 (DPS)}~\cite{chung2022diffusion}.
We use the original diffusion posterior sampling baseline with an $\ell_2$ (MSE) data-consistency potential,
\begin{equation*}
\Phi_{\ell_2}(v):= \tfrac12\|\mathcal{F}(v)-y^\delta\|_{\Sigma^{-1}}^2,
\end{equation*}
and apply the standard scalar guidance update \eqref{eq:dps_baseline_step} at each reverse-time step.
No OT-based misfit, wavefield enhancement, or preconditioned guidance is used in this baseline.

Across all methods, we use the same forward solver, acquisition geometry, and observation data. For a fair comparison, we tune the hyperparameters of each competing method as carefully as possible (within the same experimental protocol) to achieve its best attainable performance. 

\subsection{Forward Solver and Observation Data}\label{subsubsec:forward}
The forward operator $\mathcal F$ (defined in Section~\ref{sec:setting}) is evaluated by numerically solving problem \eqref{eq:solving} for the wavefield $u$ given the velocity field $v$. For a source located at $\boldsymbol{\xi}_j$, we use the source term $s(\mathbf{x},t;\boldsymbol{\xi_j})=f(t)\delta(\mathbf{x}-\boldsymbol{\xi}_j)$, where $f$ is a Ricker wavelet
\begin{equation*}
f(t)=\bigl(1-2\pi^2 f_p^2(t-t_0)^2\bigr)\exp\!\bigl(-\pi^2 f_p^2(t-t_0)^2\bigr),
\end{equation*}
with $f_p=15$ and $t_0=1.1/f_p$.
We discretize \eqref{eq:solving} on a uniform Cartesian grid of size $71\times 71$ with spacing $dx=dz=10$ and advance in time with step $dt=10^{-3}$ for $n_t=1000$ steps.
We use an 8th-order accurate finite-difference discretization in space and a second-order accurate explicit time-stepping scheme, as implemented in Deepwave~\cite{richardson_alan_2025}.
To approximate an unbounded domain, we truncate the computational domain using absorbing perfectly matched layers (PMLs) of thickness $6$ grid points on the left, right, and bottom boundaries. On the top boundary, we impose a reflecting (free-surface) boundary condition.
The observed data are generated according to \eqref{eq:data_model}, with additive Gaussian noise $\eta\sim\mathcal{N}(0,\sigma^2 I)$ and $\sigma=0.05$.

\subsection{Numerical calculation of the Wasserstein--2 distance}\label{subsec:w2_num}
For each source--receiver trace, we compute a one-dimensional Wasserstein--2 distance between the observed and synthetic signals
after enforcing nonnegativity and unit mass. In our method, the Wasserstein distance is computed using the wavefield-enhanced traces
$\tilde d_{\mathrm{obs}}(t)$ and $\tilde d_{\mathrm{syn}}(t)$ \eqref{eq:weighted_traces}.

We first shift each trace by a constant $c'>0$ to ensure positivity and then normalize to obtain probability densities:
\begin{equation*}
\tilde{\rho}_{\mathrm{obs}}(t)=
\frac{\tilde d_{\mathrm{obs}}(t)+c'}{\int_0^T \bigl(\tilde d_{\mathrm{obs}}(\tau)+c'\bigr)\,d\tau+\varepsilon'},\qquad
\tilde{\rho}_{\mathrm{syn}}(t)=
\frac{\tilde d_{\mathrm{syn}}(t)+c'}{\int_0^T \bigl(\tilde d_{\mathrm{syn}}(\tau)+c'\bigr)\,d\tau+\varepsilon'},
\end{equation*}
where the integrals are taken over the recording window $[0,T]$ and evaluated on the discrete sampling grid, $\varepsilon'=10^{-9}$ is a small constant used to avoid division by a vanishing discrete normalization factor, and we fix
$c'=1.1\bigl|\min_{t\in[0,T]} \tilde d_{\mathrm{obs}}(t)\bigr|$ in all experiments.

Let $\tilde F_{\mathrm{obs}}$ and $\tilde F_{\mathrm{syn}}$ denote the corresponding CDFs, and let
$\tilde Q_{\mathrm{obs}}=\tilde F_{\mathrm{obs}}^{-1}$ and $\tilde Q_{\mathrm{syn}}=\tilde F_{\mathrm{syn}}^{-1}$ be the associated quantile functions.
Using the one-dimensional representation
\begin{equation*}
W_2(\rho_{\mathrm{obs}},\rho_{\mathrm{syn}})=\big(\int_0^1 \bigl|\tilde Q_{\mathrm{obs}}(\xi)-\tilde Q_{\mathrm{syn}}(\xi)\bigr|^2\,d\xi\big)^{\frac{1}{2}},
\end{equation*}
we discretize $\tilde F_{\mathrm{obs}}$ and $\tilde F_{\mathrm{syn}}$ on the time grid $\{t_k\}_{k=0}^{n_t-1}$ by trapezoidal quadrature,
evaluate $\tilde Q_{\mathrm{obs}}(\xi)$ and $\tilde Q_{\mathrm{syn}}(\xi)$ via linear interpolation, and approximate the $\xi$-integral on a uniform grid
$\xi_j=j/(N_\xi-1)$ by the trapezoidal rule ($N_\xi=1000$).
The overall OT misfit is obtained by averaging the resulting $W_2$ values over all source--receiver pairs.

\subsection{Quantitative and qualitative comparisons}\label{subsec:qq_comp}
\begin{table}[th]
  \caption{Hyperparameters used in each dataset.
Here $\gamma$ is the exponent in the diagonal preconditioner $D_i$ for OT-WE-PDPS \eqref{eq:kappa}.
The symbol $\rho_0$ denotes the base guidance strength in the noise-aware schedule $\rho_i$ for OT-WE-PDPS \eqref{eq:rho},
and it also denotes the (constant) step size used in the first-order solver for the $W_2+\mathrm{TV}$ baseline.
The parameter $\alpha$ is the TV regularization weight used in the TV-regularized baselines (Baseline~1 and Baseline~2),
with $\mathrm{TV}$ defined in \eqref{eq:tv}.}
  \centering
  \small
  \setlength{\tabcolsep}{1pt} % 列间距（“一点点间隔”）

  \begin{subtable}{0.98\linewidth}
  \centering
    \begin{tabularx}{\linewidth}{
      >{\centering\arraybackslash}p{1.3cm}
      *{8}{>{\centering\arraybackslash}X}
    }
      \toprule
       & \multicolumn{4}{c}{CurveVel-A} & \multicolumn{4}{c}{CurveVel-B} \\
       \cmidrule(lr){2-5}\cmidrule(lr){6-9}
       Param & Ours & DPS & OT-WE$+\mathrm{TV}$ & $W_2+\mathrm{TV}$ &
               Ours & DPS & OT-WE$+\mathrm{TV}$ & $W_2+\mathrm{TV}$ \\
      \midrule
      $\gamma$ & 0.55 & - & 0.65 & - & 0.35 & - & 0.65 & - \\
      $\rho_0$ & 1.15 & 8 & 1.3 & 13 & 12.4 & 4 & 4 & 14 \\
      $\alpha$ & - & - & 0.05 & 0.2 & - & - & 1.0 & 0.2  \\
      \bottomrule
    \end{tabularx}
  \end{subtable}

  \vspace{0.6em}

  \begin{subtable}{0.98\linewidth}
  \centering
    \begin{tabularx}{\linewidth}{
      >{\centering\arraybackslash}p{1.3cm}
      *{8}{>{\centering\arraybackslash}X}
    }
    \toprule
      & \multicolumn{4}{c}{FlatFault-A} & \multicolumn{4}{c}{FlatFault-B} \\
      \cmidrule(lr){2-5}\cmidrule(lr){6-9}
       Param & Ours & DPS & OT-WE$+\mathrm{TV}$ & $W_2+\mathrm{TV}$ &
               Ours & DPS & OT-WE$+\mathrm{TV}$ & $W_2+\mathrm{TV}$ \\
      \midrule
      $\gamma$ & 0.35 & - & 0.65 & -  & 0.45 & - & 0.65 & - \\
      $\rho_0$& 1.2 & 0.3 & 1.1 & 5  & 1.45 & 9 & 0.6 & 2 \\
      $\alpha$ & - & - & 0.3 & 1.0 & - & - & 0.2 & 0.5 \\
      \bottomrule
    \end{tabularx}
  \end{subtable}

  \vspace{0.6em}

  \begin{subtable}{0.98\linewidth}
  \centering
    \begin{tabularx}{\linewidth}{
      >{\centering\arraybackslash}p{1.3cm}
      *{8}{>{\centering\arraybackslash}X}
    }
    \toprule
      & \multicolumn{4}{c}{CurveFault-A} & \multicolumn{4}{c}{CurveFault-B} \\
      \cmidrule(lr){2-5}\cmidrule(lr){6-9}
       Param & Ours & DPS & OT-WE$+\mathrm{TV}$ & $W_2+\mathrm{TV}$ &
               Ours & DPS & OT-WE$+\mathrm{TV}$ & $W_2+\mathrm{TV}$ \\
      \midrule
      $\gamma$ & 0.45 & - & 0.65 & - & 0.55 & - & 0.65 & - \\
      $\rho_0$ & 1.8 & 4.15 & 1.3 & 12 & 1.75 & 5 & 0.6 & 14 \\
      $\alpha$ & - & - & 0.1 & 0.2 & - & - & 0.1 & 0.5 \\
      \bottomrule
    \end{tabularx}
  \end{subtable}

  \label{tab:parameter}
\end{table}

\subsubsection{Acquisition geometry and discretization.}\label{subsubsec:operator_settings}
We consider a surface acquisition setting in which both sources and receivers are located on the free surface $z=L_z$, where $L_x=L_z=700$.
The computational domain is discretized on a $71\times 71$ Cartesian grid with spacing $d x=d z=10$~m,
where $x$ denotes the horizontal axis and $z$ the depth axis.
Grid nodes are given by $x_i=i \, dx $ and $z_j=j \, dz$ for $i,j=0,\dots,70$.
Ten point sources are placed at surface nodes
\begin{equation*}
(x_{s_k},z_{s_k})=(x_{i_k},L_z)=(i_k \, dx,L_z),\qquad
i_k\in\{0,7,14,21,28,35,42,49,56,63\},
\end{equation*}
and $70$ receivers are deployed at all surface nodes,
\begin{equation*}
(x_r,z_r)=(x_r,L_z)=(r \, d x,L_z),\qquad r=0,1,\dots,69.
\end{equation*}

\subsubsection{Test datasets.}
We evaluate six velocity models shown in Figure~\ref{fig:main}:
CurveVel-A/B, FlatFault-A/B, and CurveFault-A/B.
These models are designed to probe increasing structural complexity, from smooth layered interfaces to sharp discontinuities and their combinations.
CurveVel-A/B are layered media with clear interfaces (including curved boundaries), intended to test interface recovery.
In CurveVel-A, velocities within layers increase gradually with depth, whereas in CurveVel-B, the layer velocities are randomly distributed.
FlatFault-A/B introduces faults that shift layers and generate sharp discontinuities; compared with CurveVel-A/B exhibits more discontinuities and more severe velocity contrasts.
CurveFault-A/B combines curved stratified interfaces with fault-induced discontinuities, thereby coupling smooth geometric variation with sharp jumps and posing the most challenging scenarios among the six datasets.

\subsubsection{Settings.}
We set $k=100$ in \eqref{eq:omega_def} and $\varepsilon=10^{-4}$ in the construction of $D_i$ (\ref{eq:kappa}). When computing $\rho_i$ \eqref{eq:rho}, we set $c=0.1$ and $\tau=0$ in our method.
These settings are kept fixed across all datasets. Other hyperparameters vary across datasets and methods and are listed in Table~\ref{tab:parameter}.

For \textbf{Baseline~1 ($W_2+\mathrm{TV}$)} and \textbf{Baseline~2 (OT-WE$+\mathrm{TV}$)}, we run gradient descent for $10{,}000$ iterations and report the final iterate; step sizes and regularization weights
are listed in Table~\ref{tab:parameter}.
Across all methods, we use the same forward solver, acquisition geometry, and observation data.
Hyperparameters are tuned under a unified protocol to obtain the best attainable performance for each competing method.

\subsubsection{Results.}
\begin{figure}[tbp]
  \centering
  
  \setlength{\subfigwidth}{2.1cm}
  
  \setlength{\insetfigwidth}{6cm}

  \begin{tikzpicture}[
      subfig/.style={anchor=center, inner sep=0pt, outer sep=2pt, draw=black},
      subfig matrix/.style={
        matrix of nodes,
        column sep=0pt,    % 基础列间距
        row sep=0pt,       % 基础行间距
        nodes={subfig, minimum width=\subfigwidth, minimum height=\subfigwidth}
      },
      % 定义虚线样式
      split line/.style={dashed, gray, line width=0.5pt},
      label text/.style={anchor=east, font=\small, inner sep=5pt, align=center}
    ]

    \matrix (subfig_matrix) [subfig matrix
      ,column 1/.style={column sep=5pt}
      ,column 2/.style={column sep=5pt}
      ,column 3/.style={column sep=5pt}
      ,column 4/.style={column sep=5pt}
      ,column 5/.style={column sep=5pt}
      ,row 1/.style={row sep=8pt}
      ,row 5/.style={row sep=8pt}
    ] {
      % 第一行：前4列 + 空白列 + 后4列
      \includegraphics[width=\subfigwidth, height=\subfigwidth, keepaspectratio]{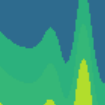} &
      \includegraphics[width=\subfigwidth, height=\subfigwidth, keepaspectratio]{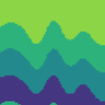} &
      \includegraphics[width=\subfigwidth, height=\subfigwidth, keepaspectratio]{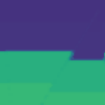} &
      \includegraphics[width=\subfigwidth, height=\subfigwidth, keepaspectratio]{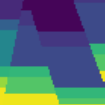} &
      \includegraphics[width=\subfigwidth, height=\subfigwidth, keepaspectratio]{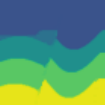} &
      \includegraphics[width=\subfigwidth, height=\subfigwidth, keepaspectratio]{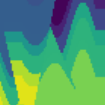}\\

      \includegraphics[width=\subfigwidth, height=\subfigwidth, keepaspectratio]{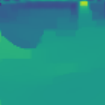} &
      \includegraphics[width=\subfigwidth, height=\subfigwidth, keepaspectratio]{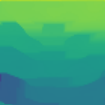} &
      \includegraphics[width=\subfigwidth, height=\subfigwidth, keepaspectratio]{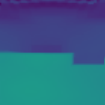} &
      \includegraphics[width=\subfigwidth, height=\subfigwidth, keepaspectratio]{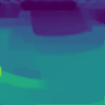} &
      \includegraphics[width=\subfigwidth, height=\subfigwidth, keepaspectratio]{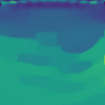} &
      \includegraphics[width=\subfigwidth, height=\subfigwidth, keepaspectratio]{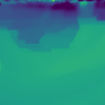}\\

      \includegraphics[width=\subfigwidth, height=\subfigwidth, keepaspectratio]{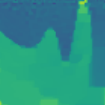} &
      \includegraphics[width=\subfigwidth, height=\subfigwidth, keepaspectratio]{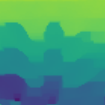} &
      \includegraphics[width=\subfigwidth, height=\subfigwidth, keepaspectratio]{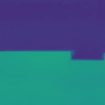} &
      \includegraphics[width=\subfigwidth, height=\subfigwidth, keepaspectratio]{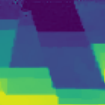} &
      \includegraphics[width=\subfigwidth, height=\subfigwidth, keepaspectratio]{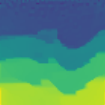} &
      \includegraphics[width=\subfigwidth, height=\subfigwidth, keepaspectratio]{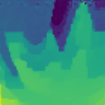}\\

      \includegraphics[width=\subfigwidth, height=\subfigwidth, keepaspectratio]{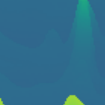} &
      \includegraphics[width=\subfigwidth, height=\subfigwidth, keepaspectratio]{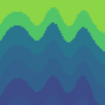} &
      \includegraphics[width=\subfigwidth, height=\subfigwidth, keepaspectratio]{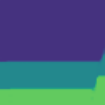} &
      \includegraphics[width=\subfigwidth, height=\subfigwidth, keepaspectratio]{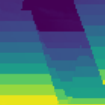} &
      \includegraphics[width=\subfigwidth, height=\subfigwidth, keepaspectratio]{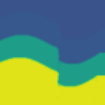} &
      \includegraphics[width=\subfigwidth, height=\subfigwidth, keepaspectratio]{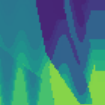}\\

      \includegraphics[width=\subfigwidth, height=\subfigwidth, keepaspectratio]{picture/main_curvevela_result_ours} &
      \includegraphics[width=\subfigwidth, height=\subfigwidth, keepaspectratio]{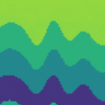} &
      \includegraphics[width=\subfigwidth, height=\subfigwidth, keepaspectratio]{picture/main_flata_result_ours} &
      \includegraphics[width=\subfigwidth, height=\subfigwidth, keepaspectratio]{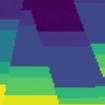} &
      \includegraphics[width=\subfigwidth, height=\subfigwidth, keepaspectratio]{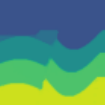} &
      \includegraphics[width=\subfigwidth, height=\subfigwidth, keepaspectratio]{picture/main_curvefaultb_result_ours}\\

      \includegraphics[width=\subfigwidth, height=\subfigwidth, keepaspectratio]{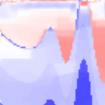} &
      \includegraphics[width=\subfigwidth, height=\subfigwidth, keepaspectratio]{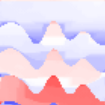} &
      \includegraphics[width=\subfigwidth, height=\subfigwidth, keepaspectratio]{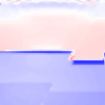} &
      \includegraphics[width=\subfigwidth, height=\subfigwidth, keepaspectratio]{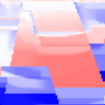} &
      \includegraphics[width=\subfigwidth, height=\subfigwidth, keepaspectratio]{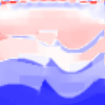} &
      \includegraphics[width=\subfigwidth, height=\subfigwidth, keepaspectratio]{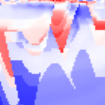}\\

      \includegraphics[width=\subfigwidth, height=\subfigwidth, keepaspectratio]{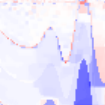} &
      \includegraphics[width=\subfigwidth, height=\subfigwidth, keepaspectratio]{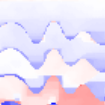} &
      \includegraphics[width=\subfigwidth, height=\subfigwidth, keepaspectratio]{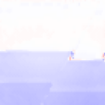} &
      \includegraphics[width=\subfigwidth, height=\subfigwidth, keepaspectratio]{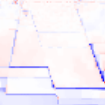} &
      \includegraphics[width=\subfigwidth, height=\subfigwidth, keepaspectratio]{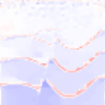} &
      \includegraphics[width=\subfigwidth, height=\subfigwidth, keepaspectratio]{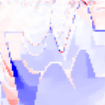}\\

      \includegraphics[width=\subfigwidth, height=\subfigwidth, keepaspectratio]{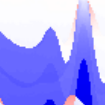} &
      \includegraphics[width=\subfigwidth, height=\subfigwidth, keepaspectratio]{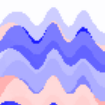} &
      \includegraphics[width=\subfigwidth, height=\subfigwidth, keepaspectratio]{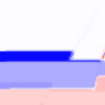} &
      \includegraphics[width=\subfigwidth, height=\subfigwidth, keepaspectratio]{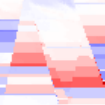} &
      \includegraphics[width=\subfigwidth, height=\subfigwidth, keepaspectratio]{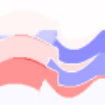} &
      \includegraphics[width=\subfigwidth, height=\subfigwidth, keepaspectratio]{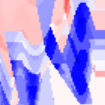}\\
            
      \includegraphics[width=\subfigwidth, height=\subfigwidth, keepaspectratio]{picture/main_curvevela_diff_ours} &
      \includegraphics[width=\subfigwidth, height=\subfigwidth, keepaspectratio]{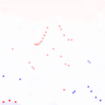} &
      \includegraphics[width=\subfigwidth, height=\subfigwidth, keepaspectratio]{picture/main_flata_diff_ours} &
      \includegraphics[width=\subfigwidth, height=\subfigwidth, keepaspectratio]{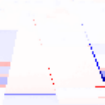} &
      \includegraphics[width=\subfigwidth, height=\subfigwidth, keepaspectratio]{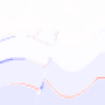} &
      \includegraphics[width=\subfigwidth, height=\subfigwidth, keepaspectratio]{picture/main_curvefaultb_diff_ours}\\
    };

    \coordinate (inset_pos1) at ($(subfig_matrix-2-6.east)!0.5!(subfig_matrix-5-6.east)+(35pt,0)$);
    \coordinate (inset_pos2) at ($(subfig_matrix-6-6.east)!0.5!(subfig_matrix-9-6.east)+(35pt,0)$);

    \node[minimum width=\insetfigwidth, minimum height=\insetfigwidth] at (inset_pos1) {
      \includegraphics[width=\insetfigwidth, height=\insetfigwidth, keepaspectratio]{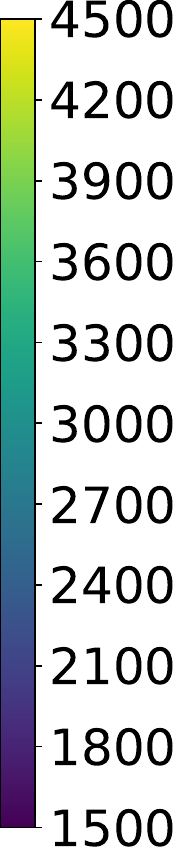}
    };
    \node[right=15pt, font=\small] at (inset_pos1) {\rotatebox{90}{velocity (m/s)}};

    \node[minimum width=\insetfigwidth, minimum height=\insetfigwidth] at (inset_pos2) {
      \includegraphics[width=\insetfigwidth, height=\insetfigwidth, keepaspectratio]{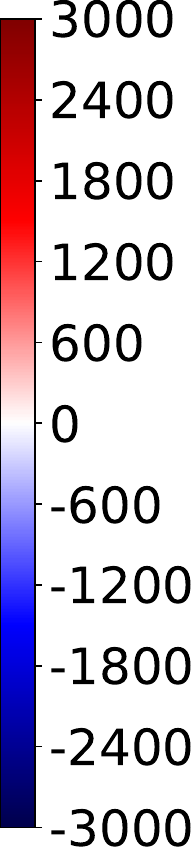}
    };
    \node[right=15pt, font=\small] at (inset_pos2) {\rotatebox{90}{velocity (m/s)}};

    % \draw[split line] 
    %   ($(subfig_matrix-1-1.west)!0.5!(subfig_matrix-2-1.west)+(-7pt,0)$) -- ($(subfig_matrix-1-6.east)!0.5!(subfig_matrix-2-6.east)+(7pt,0)$);

    \draw[split line] 
      ($(subfig_matrix-5-1.west)!0.5!(subfig_matrix-6-1.west)+(-7pt,0)$) -- ($(subfig_matrix-5-6.east)!0.5!(subfig_matrix-6-6.east)+(7pt,0)$);
    % 在第一条虚线中间添加Result文字
    
    \node[font=\small, anchor=center] at (
      $(subfig_matrix-2-6.east)!0.5!(subfig_matrix-5-6.east)+(5pt,0)$
    ) {\rotatebox{90}{Reconstruction}};

    \node[font=\small, anchor=center] at (
      $(subfig_matrix-6-6.east)!0.5!(subfig_matrix-9-6.east)+(5pt,0)$
    ) {\rotatebox{90}{Difference}};
    
    % 顶部标签
    \node[above=28pt, font=\small] at (subfig_matrix-1-1) {CurveVel-A};
    \node[above=28pt, font=\small] at (subfig_matrix-1-2) {CurveVel-B};
    \node[above=28pt, font=\small] at (subfig_matrix-1-3) {FlatFault-A};
    \node[above=28pt, font=\small] at (subfig_matrix-1-4) {FlatFault-B};
    \node[above=28pt, font=\small] at (subfig_matrix-1-5) {CurveFault-A};
    \node[above=28pt, font=\small] at (subfig_matrix-1-6) {CurveFault-B};

    \node[left=28pt, font=\small] at (subfig_matrix-1-1) {\rotatebox{90}{Truth}};
    \node[left=28pt, font=\small] at (subfig_matrix-5-1) {\rotatebox{90}{\textbf{Ours}}};
    \node[left=28pt, font=\small] at (subfig_matrix-9-1) {\rotatebox{90}{\textbf{Ours}}};
    \node[left=28pt, font=\small] at (subfig_matrix-4-1) {\rotatebox{90}{DPS}};
    \node[left=28pt, font=\small] at (subfig_matrix-8-1) {\rotatebox{90}{DPS}};
    \node[left=28pt, font=\small] at (subfig_matrix-3-1) {\rotatebox{90}{OT-WE$+\mathrm{TV}$}};
    \node[left=28pt, font=\small] at (subfig_matrix-7-1) {\rotatebox{90}{OT-WE$+\mathrm{TV}$}};
    \node[left=28pt, font=\small] at (subfig_matrix-2-1) {\rotatebox{90}{$W_2+\mathrm{TV}$}};
    \node[left=28pt, font=\small] at (subfig_matrix-6-1) {\rotatebox{90}{$W_2+\mathrm{TV}$}};

    % 底部标注（修正坐标，适配新的列结构）
    % \node[label text, anchor=north, inner sep=5pt, font=\large,align=center] at ($(subfig_matrix-2-1)!0.5!(subfig_matrix-3-1)+(-34pt, 22pt)$) {\rotatebox{90}{D   D   I   M}};

  \end{tikzpicture}
  \caption{Inversion results across different datasets. The first row displays the true velocity field $v_\mathrm{true}$. Rows 2 to 5 show the inverted velocity fields $v_\mathrm{rec}$ obtained using different algorithms. Rows 6 to 9 present the difference between the true and inverted velocity fields, $v_\mathrm{rec}-v_\mathrm{true}$. The observed wavefield has been perturbed with additive noise of intensity $\sigma = 0.05$.}
  \label{fig:main}
\end{figure}

The following observations are drawn from Table~\ref{tab:main} and Figure~\ref{fig:main}:
\begin{table}[t]  % 核心：将table*改为table，适配单栏
\caption{Quantitative comparison of inversion methods on OpenFWI velocity models.
Metrics are the relative $\ell_2$ error $e_{\ell_2}$, PSNR, and SSIM for the reconstructed velocity field (arrows indicate the preferred direction).
Methods include $W_2+\mathrm{TV}$ (Baseline 1), OT-WE$+\mathrm{TV}$ (Baseline 2), DPS (Baseline 3), and the proposed method.}
\centering
\small
\setlength{\tabcolsep}{4.2pt}
\renewcommand{\arraystretch}{1.12}
% 缩放表格至单栏宽度的98%（避免过宽）
\begin{subtable}{0.98\linewidth}
\centering
\begin{tabular}{c ccc ccc ccc}  % 添加竖线分隔：第一列后、第四列后、第七列后
\toprule
& \multicolumn{3}{c}{CurveVel-A} & \multicolumn{3}{c}{CurveVel-B} & \multicolumn{3}{c}{FlatFault-A}  \\
\cmidrule(lr){2-4}\cmidrule(lr){5-7}\cmidrule(lr){8-10}
Method & $e_{l_2}$$\downarrow$ & PSNR$\uparrow$ & SSIM$\uparrow$ & $e_{l_2}$$\downarrow$ & PSNR$\uparrow$ & SSIM$\uparrow$ & $e_{l_2}$$\downarrow$ & PSNR$\uparrow$ & SSIM$\uparrow$ \\
\midrule
$W_2+\mathrm{TV}$ & 12.10\% & 20.65 & 0.4971 & 10.72\% & 21.05 & 0.5425 & 11.84\% & 20.56 & 0.7198 \\
OT-WE$+\mathrm{TV}$ & 7.96\% & 24.28 & 0.6987 & 7.40\% & 24.27 & 0.6434 & 7.28\% & 24.79 & 0.9217 \\
DPS & 23.15\% & 15.01 & 0.5747 & 15.70\% & 17.73 & 0.4329 & 16.68\% & 17.58 & 0.6965 \\
\textbf{Ours} & \textbf{3.50\%} & \textbf{31.42} & \textbf{0.9039} & \textbf{2.04\%} & \textbf{35.45} & \textbf{0.9517} & \textbf{2.91\%} & \textbf{32.74} & \textbf{0.9646} \\
\bottomrule
\end{tabular}
\end{subtable}  
    
\vspace{0.6em}

\begin{subtable}{0.98\linewidth}
\centering
\begin{tabular}{c ccc ccc ccc}
\toprule
& \multicolumn{3}{c}{FlatFault-B} & \multicolumn{3}{c}{CurveFault-A} & \multicolumn{3}{c}{CurveFault-B}  \\
\cmidrule(lr){2-4}\cmidrule(lr){5-7}\cmidrule(lr){8-10}
Method & $e_{l_2}$$\downarrow$ & PSNR$\uparrow$ & SSIM$\uparrow$ & $e_{l_2}$$\downarrow$ & PSNR$\uparrow$ & SSIM$\uparrow$ & $e_{l_2}$$\downarrow$ & PSNR$\uparrow$ & SSIM$\uparrow$ \\
\midrule  % 仅在Method行下方添加贯穿所有列的cmidrule
$W_2+\mathrm{TV}$ & 22.82\% & 17.05 & 0.4584 & 17.63\% & 17.74 & 0.4875 & 16.34\% & 18.20 & 0.4343 \\
OT-WE$+\mathrm{TV}$ & 6.31\% & 28.21 & 0.8729 & 4.79\% & 29.06 & 0.7897 & 9.39\% & 23.01 & 0.6997 \\
DPS & 13.54\% & 21.58 & 0.5662 & 10.72\% & 22.06 & 0.6149 & 21.82\% & 15.69 & 0.2210 \\
\textbf{Ours} & \textbf{6.13\%} & \textbf{28.48} & \textbf{0.9272} & \textbf{2.41\%} & \textbf{35.03} & \textbf{0.9540} & \textbf{5.32\%} & \textbf{27.95} & \textbf{0.7956} \\
\bottomrule
\end{tabular}
\end{subtable}
\label{tab:main}
\end{table}

\begin{figure}[t]
    \centering
    
    % 第一幅子图
    \begin{subfigure}[b]{0.45\textwidth}
        \centering
        \includegraphics[width=\textwidth]{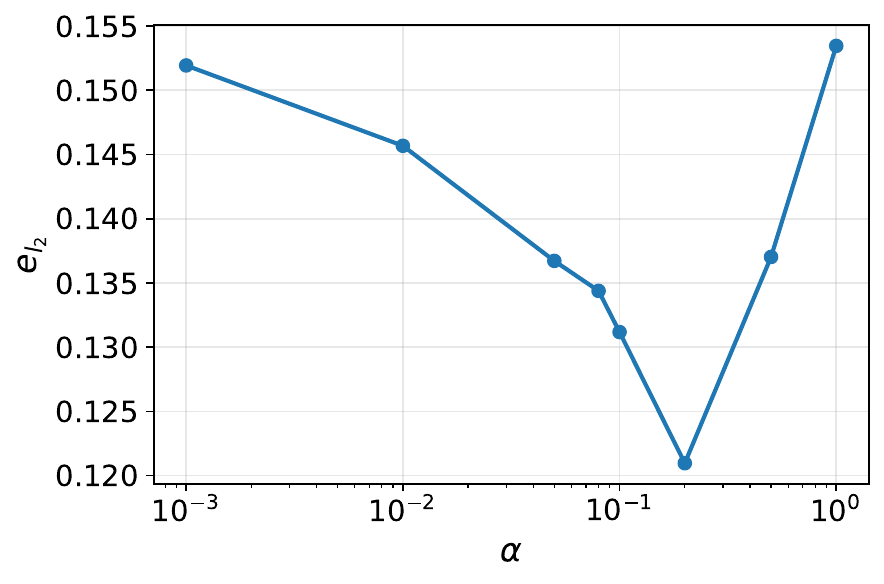} % 替换为您的图片路径
        \caption{Relative $l_2$-error of the reconstructed velocity field $v$ for the $W_2+\mathrm{TV}$ method under different values of the TV regularization parameter $\alpha$.}
        \label{fig:trad_alpha}
    \end{subfigure}
    \hfill % 添加水平间距
    % 第二幅子图
    \begin{subfigure}[b]{0.45\textwidth}
        \centering
        \includegraphics[width=\textwidth]{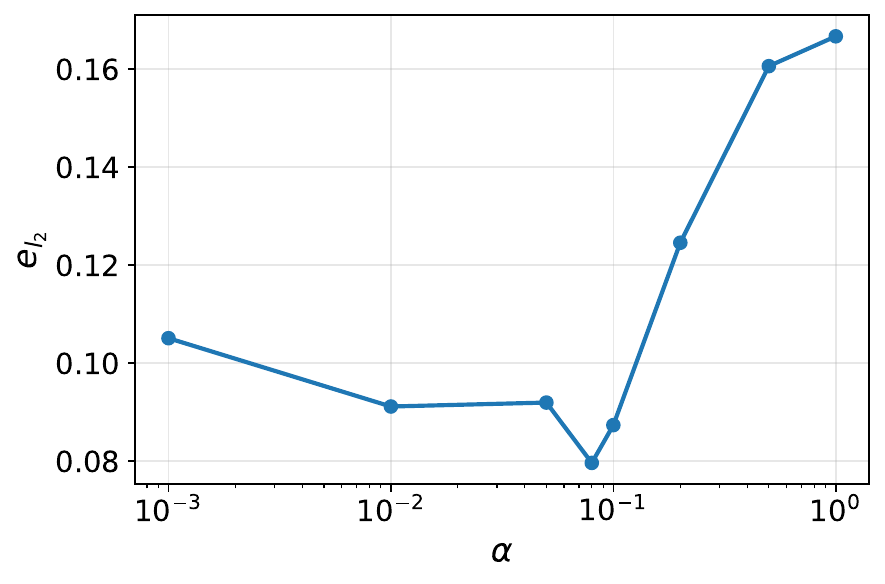} % 替换为您的图片路径
        \caption{Relative $l_2$-error of the reconstructed velocity field $v$ for the OT-WE$+\mathrm{TV}$ method under different values of the TV regularization parameter $\alpha$.}
        \label{fig:tradt_alpha}
    \end{subfigure}
    
    \caption{Relative $l_2$-error of the reconstructed velocity field $v$ versus the TV regularization parameter $\alpha$ for the two deterministic methods: (a) $W_2+\mathrm{TV}$ and (b) OT-WE$+\mathrm{TV}$.}
    \label{fig:trad}
\end{figure}

\begin{figure}[th]
    \centering
    
    % 第一幅子图
    \begin{subfigure}[b]{0.45\textwidth}
        \centering
        \includegraphics[width=\textwidth]{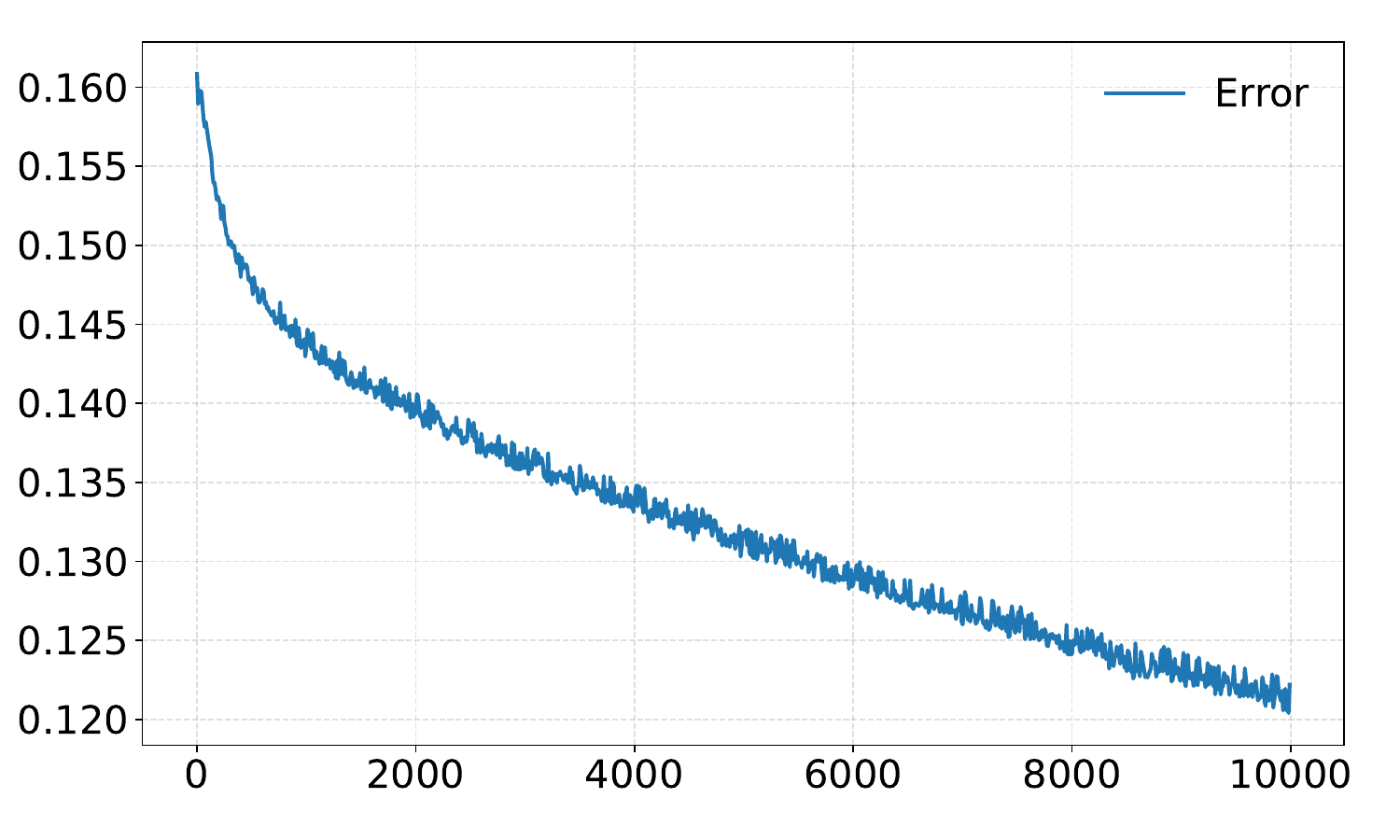} % 替换为您的图片路径
        \caption{ Relative $l_2$-error of the reconstructed velocity field $v$ versus the iteration number for the $W_2+\mathrm{TV}$ method with the optimal TV regularization parameter $\alpha=0.2$.}
        \label{fig:trad_error1}
    \end{subfigure}
    \hfill % 添加水平间距
    % 第二幅子图
    \begin{subfigure}[b]{0.45\textwidth}
        \centering
        \includegraphics[width=\textwidth]{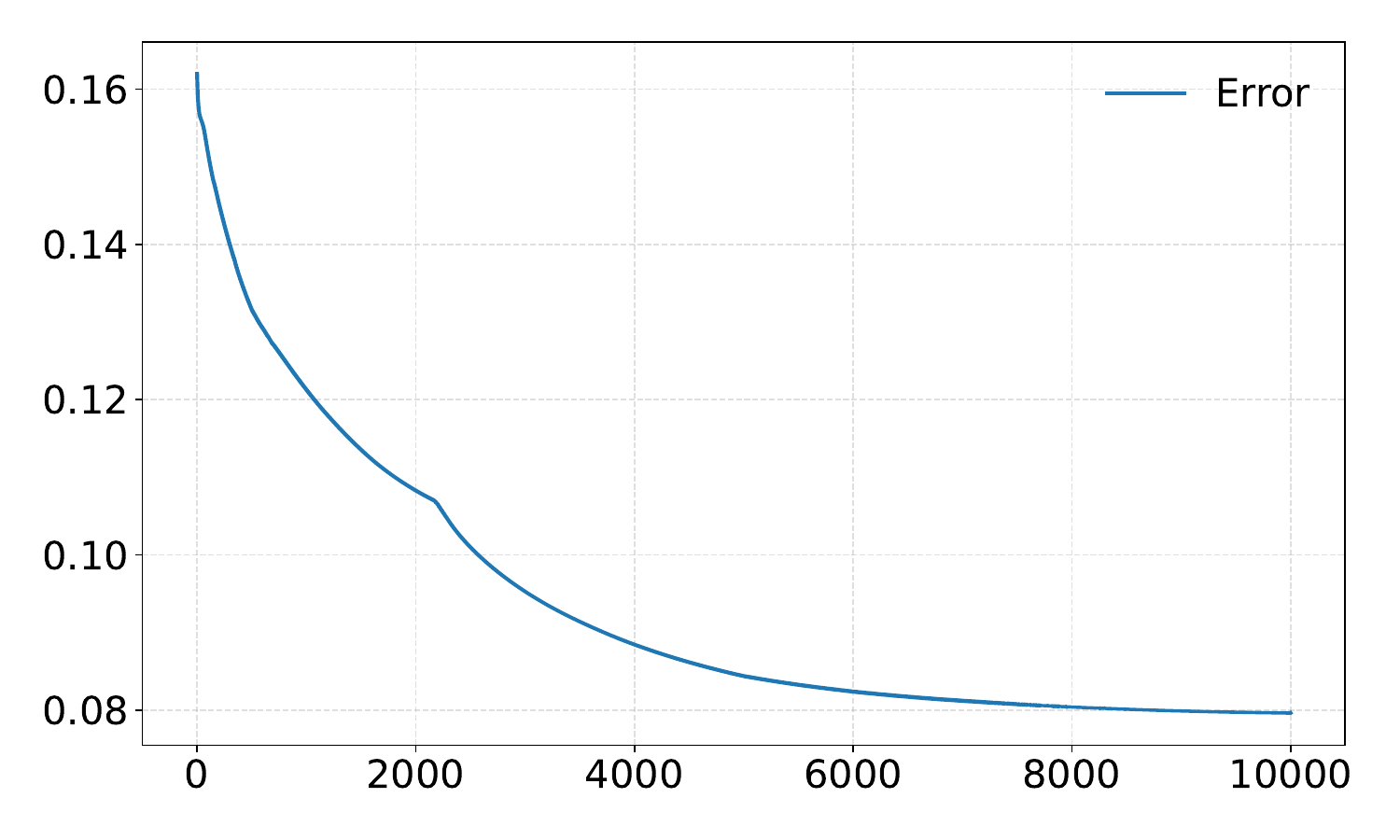} % 替换为您的图片路径
        \caption{Relative $l_2$-error of the reconstructed velocity field $v$ versus the iteration number for the OT-WE$+\mathrm{TV}$ method with the optimal TV regularization parameter $\alpha=0.08$.}
        \label{fig:trad_error2}
    \end{subfigure}
    
    \caption{Relative $l_2$-error of the reconstructed velocity field $v$ versus the iteration number for the two deterministic methods with their optimal TV regularization parameters $\alpha$: (a) $W_2+\mathrm{TV}$ and (b) OT-WE$+\mathrm{TV}$.}
    \label{fig:error_trad}
\end{figure}
\begin{enumerate}
\item Our method attains the lowest reconstruction errors and the highest image-quality metrics across the datasets,
while preserving interfaces and fault geometries in Figure~\ref{fig:main}.
In contrast, the DPS baseline often yields reconstructions with a larger mismatch to the ground truth, suggesting that the standard MSE-based wavefield guidance may be too weak to enforce data consistency under surface-only acquisition.
These results indicate that our proposed method achieves higher reconstruction quality than the DPS baseline under surface-only acquisition.

\item $W_2+\mathrm{TV}$ baseline often produces over-smoothed reconstructions, with blurred interfaces and loss of structural details,
which is reflected by lower SSIM and visually smeared layer boundaries.
Our method mitigates this over-smoothing and recovers the main structures, although errors remain in challenging regions.

\item As the model complexity increases (e.g., more layers, stronger discontinuities, and curved/faulted interfaces),
the inversion problem becomes more ill-posed and the reconstruction accuracy degrades.
For our method, residual errors concentrate near sharp velocity transitions and at larger depths,
consistent with reduced illumination and parameter sensitivity in deeper regions for surface-only acquisition.

\item The OT-WE$+\mathrm{TV}$ baseline frequently outperforms the $W_2+\mathrm{TV}$ baseline in terms of both quantitative metrics and visual quality. This observation suggests that wavefield enhancement and adaptive step sizing are effective not only within the DPS framework, but may also improve the performance of other inversion methods.
\end{enumerate}

To verify that the selected regularization parameter $\alpha$ is appropriate for the deterministic methods, we take Dataset CurveVel-A as an example and plot, in Figure~\ref{fig:trad}, the relationship between the regularization parameter and the final relative error for (i) the $W_2+\mathrm{TV}$ method and (ii) the OT-WE$+\mathrm{TV}$ method.
In addition, Figure~\ref{fig:error_trad} shows the evolution of the relative error with respect to the iteration number for these two methods. Since Figure~\ref{fig:trad_error1} indicates that the relative error of the $W_2+\mathrm{TV}$ method is still decreasing, we further continue this method to $80{,}000$ iterations under the same parameter setting. Even then, the relative error is only reduced to $10.5\%$, which is still higher than the best relative error achieved by the OT-WE$+\mathrm{TV}$ method, namely, $8\%$.

The results indicate that, for the deterministic methods, the final inversion error is sensitive to the choice of the TV regularization parameter $\alpha$, and a best-performing value can be identified within the tested parameter range. For the remaining datasets, the parameter $\alpha$ was selected by the same tuning procedure in order to achieve the best performance attainable within the tested parameter range for each deterministic method. In addition, measured in terms of iteration count, the deterministic methods require substantially more iterations to converge than our method.
\begin{figure}[th]
  \centering
  
  \setlength{\subfigwidth}{2.1cm}
  
  \setlength{\insetfigwidth}{4cm}

  \begin{tikzpicture}[
      subfig/.style={anchor=center, inner sep=0pt, outer sep=2pt, draw=black},
      subfig matrix/.style={
        matrix of nodes,
        column sep=0pt,    % 基础列间距
        row sep=0pt,       % 基础行间距
        nodes={subfig, minimum width=\subfigwidth, minimum height=\subfigwidth}
      },
      % 定义虚线样式
      no border/.style={subfig, draw=none},
      split line/.style={dashed, gray, line width=0.5pt},
      label text/.style={anchor=east, font=\small, inner sep=5pt, align=center}
    ]

    \matrix (subfig_matrix) [subfig matrix
      ,row 1/.style={row sep=18pt}
      ,row 3/.style={row sep=8pt}
    ] {
      % 第一行：前4列 + 空白列 + 后4列
      \includegraphics[width=\subfigwidth, height=\subfigwidth, keepaspectratio]{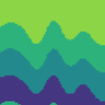}\\
      
      \includegraphics[width=\subfigwidth, height=\subfigwidth, keepaspectratio]{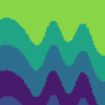} &
      \includegraphics[width=\subfigwidth, height=\subfigwidth, keepaspectratio]{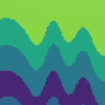} &
      \includegraphics[width=\subfigwidth, height=\subfigwidth, keepaspectratio]{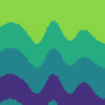} &
      \includegraphics[width=\subfigwidth, height=\subfigwidth, keepaspectratio]{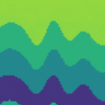} \\

      \includegraphics[width=\subfigwidth, height=\subfigwidth, keepaspectratio]{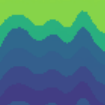} &
      \includegraphics[width=\subfigwidth, height=\subfigwidth, keepaspectratio]{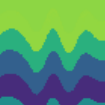} &
      \includegraphics[width=\subfigwidth, height=\subfigwidth, keepaspectratio]{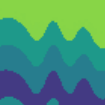} &
      \includegraphics[width=\subfigwidth, height=\subfigwidth, keepaspectratio]{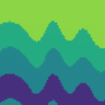} \\

      \includegraphics[width=\subfigwidth, height=\subfigwidth, keepaspectratio]{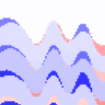} &
      \includegraphics[width=\subfigwidth, height=\subfigwidth, keepaspectratio]{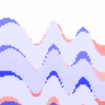} &
      \includegraphics[width=\subfigwidth, height=\subfigwidth, keepaspectratio]{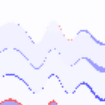} &
      \includegraphics[width=\subfigwidth, height=\subfigwidth, keepaspectratio]{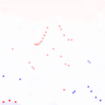} \\

      \includegraphics[width=\subfigwidth, height=\subfigwidth, keepaspectratio]{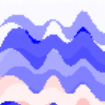} &
      \includegraphics[width=\subfigwidth, height=\subfigwidth, keepaspectratio]{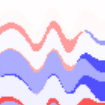} &
      \includegraphics[width=\subfigwidth, height=\subfigwidth, keepaspectratio]{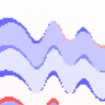} &
      \includegraphics[width=\subfigwidth, height=\subfigwidth, keepaspectratio]{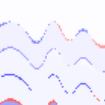} \\
    };

    \coordinate (inset_pos1) at ($(subfig_matrix-2-4.east)!0.5
    !(subfig_matrix-3-4)+(55pt,0)$);
    \coordinate (inset_pos2) at ($(subfig_matrix-4-4.east)!0.5!(subfig_matrix-5-4)+(55pt,0)$);

    \node[minimum width=\insetfigwidth, minimum height=\insetfigwidth] at (inset_pos1) {
      \includegraphics[width=\insetfigwidth, height=\insetfigwidth, keepaspectratio]{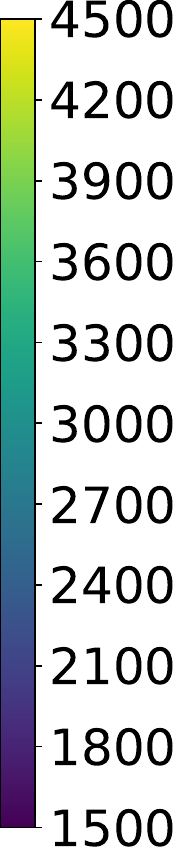}
    };
    \node[right=15pt, font=\small] at (inset_pos1) {\rotatebox{90}{velocity (m/s)}};

    \node[minimum width=\insetfigwidth, minimum height=\insetfigwidth] at (inset_pos2) {
      \includegraphics[width=\insetfigwidth, height=\insetfigwidth, keepaspectratio]{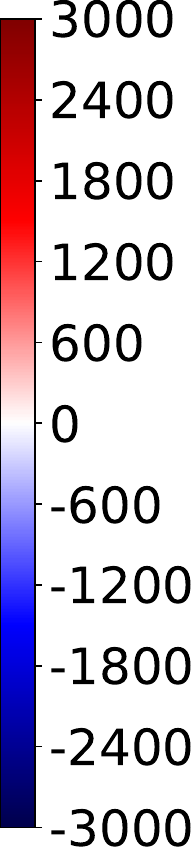}
    };
    \node[right=15pt, font=\small] at (inset_pos2) {\rotatebox{90}{velocity (m/s)}};

    \draw[split line] 
      ($(subfig_matrix-3-1.west)!0.5!(subfig_matrix-4-1.west)+(-7pt,0)$) -- ($(subfig_matrix-3-4.east)!0.5!(subfig_matrix-4-4.east)+(7pt,0)$);
    
    \node[font=\small, anchor=center] at (
      $(subfig_matrix-2-4.east)!0.5!(subfig_matrix-3-4)+(30pt,0)$
    ) {\rotatebox{90}{Reconstruction}};

    \node[font=\small, anchor=center] at (
      $(subfig_matrix-4-4.east)!0.5!(subfig_matrix-5-4)+(30pt,0)$
    ) {\rotatebox{90}{Difference}};
    
    % 顶部标签
    \node[above=30pt, font=\small] at (subfig_matrix-1-1) {\texttt{Truth}};
    \node[above=30pt, font=\small] at (subfig_matrix-2-1) {\texttt{Only}};
    \node[above=30pt, font=\small] at (subfig_matrix-2-2) {\texttt{Step}};
    \node[above=32pt, font=\small] at (subfig_matrix-2-3) {\texttt{Enh}};
    \node[above=32pt, font=\small] at (subfig_matrix-2-4) {\texttt{All}};

    \node[left=33pt, font=\small] at (subfig_matrix-2-1) {\rotatebox{90}{\texttt{$W_2$}}};
    \node[left=33pt, font=\small] at (subfig_matrix-3-1) {\rotatebox{90}{\texttt{MSE}}};
    \node[left=33pt, font=\small] at (subfig_matrix-4-1) {\rotatebox{90}{\texttt{$W_2$}}};
    \node[left=33pt, font=\small] at (subfig_matrix-5-1) {\rotatebox{90}{\texttt{MSE}}};

  \end{tikzpicture}
  \caption{Ablation study on CurveVel-B. The top panel shows the ground-truth velocity model. The second and third panels report reconstructions obtained with $W_2$-based and MSE-based guidance, respectively. From left to right, the columns correspond to \texttt{Only} (baseline DPS), \texttt{Step} (preconditioned guidance update $P_i=\rho_iD_i$), \texttt{Enh} (wavefield enhancement: amplitude-adaptive weighting and misfit-scale normalization), and \texttt{All} (\texttt{Step}+\texttt{Enh}). The bottom two panels show the corresponding error maps (reconstruction minus ground truth).}
  \label{fig:abla}
\end{figure}

\begin{table}[t]
  \caption{Ablation results on CurveVel-B: relative $\ell_2$ error ($e_{\ell_2}$; lower is better), PSNR and SSIM (higher is better). Here \texttt{Only} denotes the baseline DPS update with the specified misfit, \texttt{Enh} enables wavefield enhancement (amplitude-adaptive data weighting and misfit-scale normalization), \texttt{Step} enables the preconditioned guidance update $P_i=\rho_iD_i$, and \texttt{All} enables both \texttt{Enh} and \texttt{Step}.}
  \label{tab:ablation}
  \centering
  \small
  \setlength{\tabcolsep}{5pt}
  \renewcommand{\arraystretch}{1.08}
  \begin{tabular}{lccc}
    \toprule
    Method & $e_{\ell_2}\downarrow$ & PSNR$\uparrow$ & SSIM$\uparrow$ \\
    \midrule
    \multicolumn{4}{l}{\textbf{MSE guidance}} \\
    \texttt{Only} & 18.21\% & 16.45 & 0.3827 \\
    \texttt{Step} & 12.52\% & 19.70 & 0.5082 \\
    \texttt{Enh}  &  9.20\% & 22.37 & 0.6578 \\
    \texttt{All}  &  5.10\% & 27.50 & 0.8226 \\
    \midrule
    \multicolumn{4}{l}{\textbf{$W_2$ guidance}} \\
    \texttt{Only} & 10.81\% & 20.98 & 0.6307 \\
    \texttt{Step} &  8.65\% & 22.91 & 0.6726 \\
    \texttt{Enh}  &  5.33\% & 27.13 & 0.8244 \\
    \textbf{\texttt{All} (ours)} & \textbf{2.04\%} & \textbf{35.45} & \textbf{0.9517} \\
    \bottomrule
  \end{tabular}
\end{table}

\subsection{Ablation Study}\label{subsec:abla}
In the preceding experiments, we evaluated the reconstruction performance of our method for FWI.
We next report an ablation study to assess the contribution of its principal components.
Taking the original DPS framework as the reference, we investigate:
(i) the choice of data-misfit metric, comparing MSE with the Wasserstein-$2$ distance ($W_2$);
(ii) \emph{wavefield enhancement}, consisting of \emph{amplitude-adaptive data weighting} and \emph{misfit-scale normalization}
(Section~\ref{subsec:ot_misfit});
and (iii) the \emph{preconditioned guidance update} $P_i=\rho_iD_i$ (Section~\ref{subsec:precond_dps}), which combines stabilized guidance scheduling with spatially adaptive diagonal scaling.
Unless otherwise noted, all ablation experiments are performed on the CurveVel-B model in Figure~\ref{fig:main}.
For each setting, hyperparameters are selected according to the same tuning protocol so as to achieve the best performance attainable within that setting.

\begin{figure}[th]
    \centering
    \includegraphics[width=0.5\textwidth]{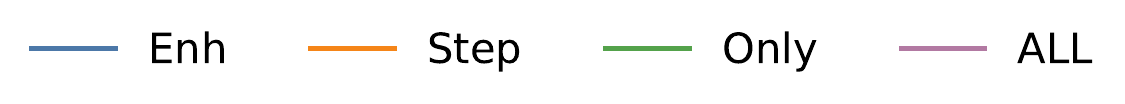}
    \vspace{1ex} % 与下方子图保持适当间距
    
    \begin{subfigure}[t]{0.3\textwidth}
        \centering
        \includegraphics[width=0.9\linewidth]{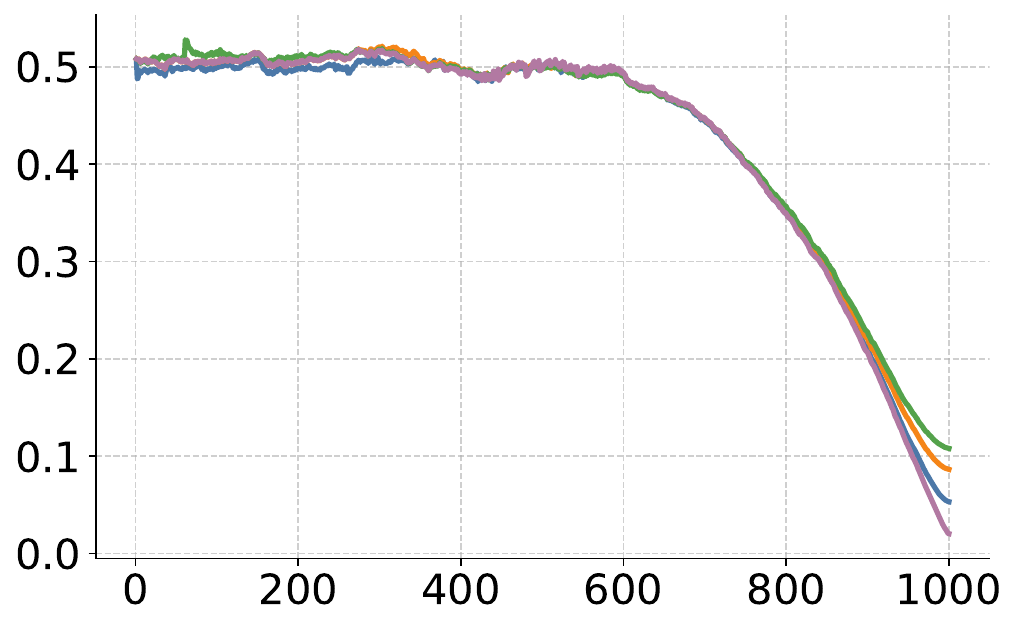}
        \caption{Relative $\ell_2$-error of $v_i$ using $W_2$}
        \label{fig:abla_w2_el2_error}
    \end{subfigure}
    \hfill
    \begin{subfigure}[t]{0.3\textwidth}
        \centering
        \includegraphics[width=0.9\linewidth]{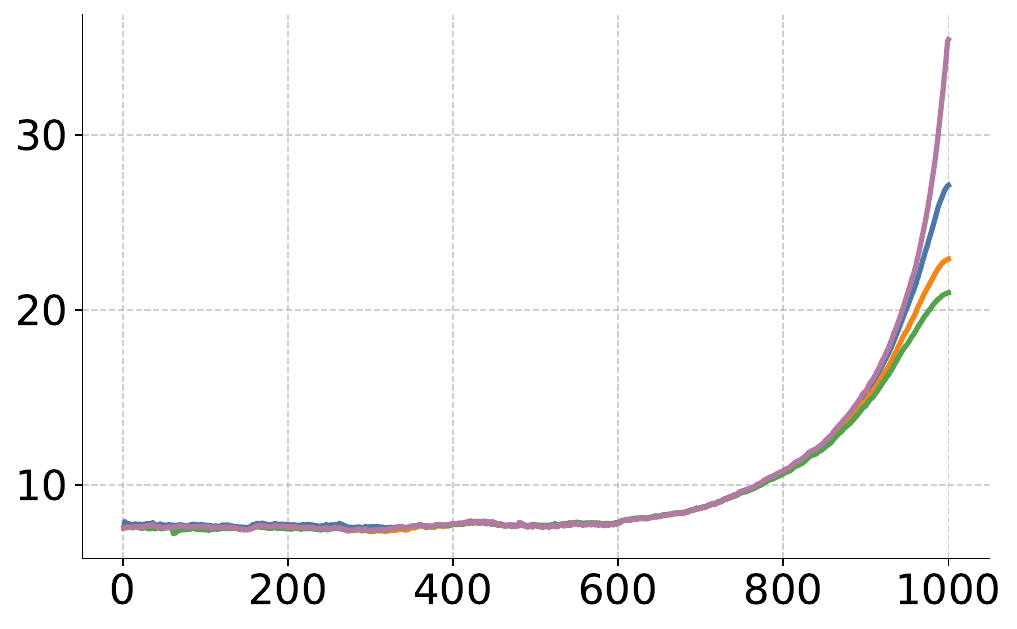}
        \caption{PSNR of $v_i$ using $W_2$}
        \label{fig:abla_w2_psnr_error}
    \end{subfigure}
    \hfill
    \begin{subfigure}[t]{0.3\textwidth}
        \centering
        \includegraphics[width=0.9\linewidth]{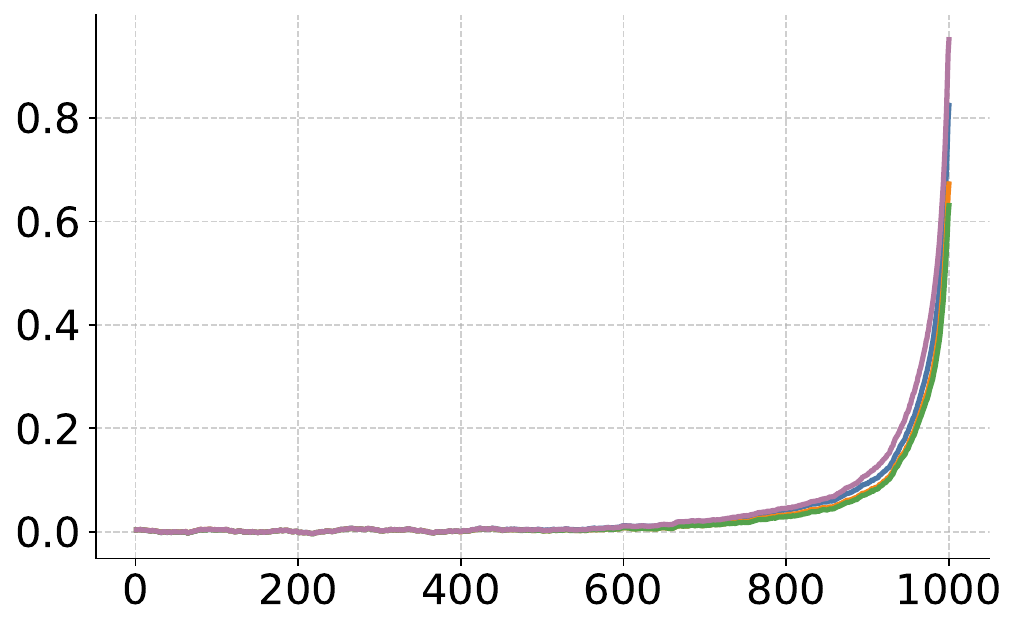}
        \caption{SSIM of $v_i$ using $W_2$}
        \label{fig:abla_w2_ssim_error}
    \end{subfigure}
    \hfill
    
    \begin{subfigure}[t]{0.3\textwidth}
        \centering
        \includegraphics[width=0.9\linewidth]{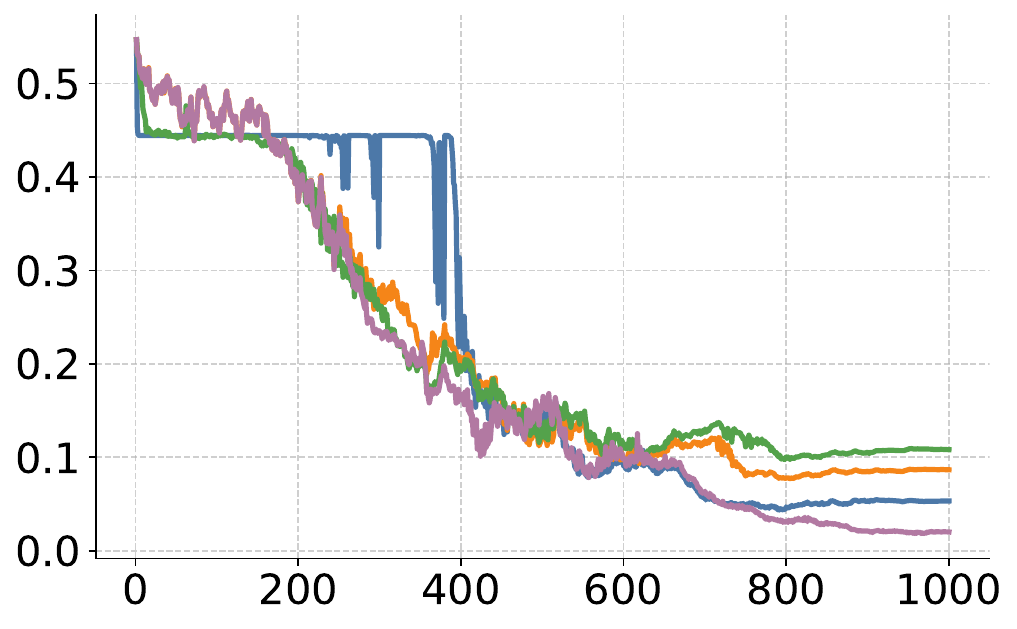}
        \caption{Relative $\ell_2$-error of $\hat v_0^{(i)}$ using $W_2$}
        \label{fig:abla_w2_el2_error0}
    \end{subfigure}
    \hfill
    \begin{subfigure}[t]{0.3\textwidth}
        \centering
        \includegraphics[width=0.9\linewidth]{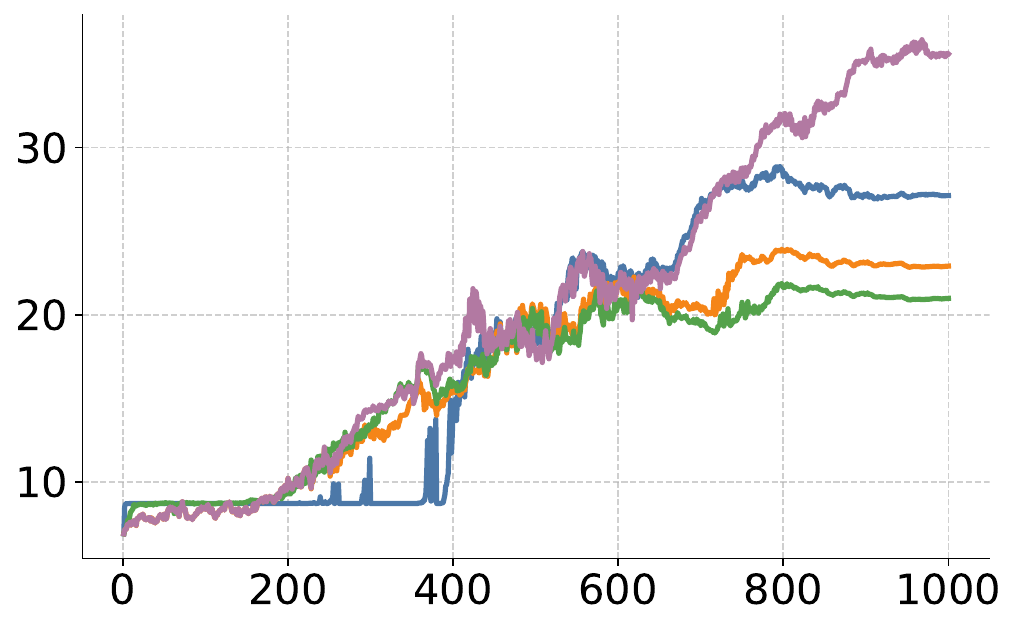}
        \caption{PSNR of $\hat v_0^{(i)}$ using $W_2$}
        \label{fig:abla_w2_psnr_error0}
    \end{subfigure}
    \hfill
    \begin{subfigure}[t]{0.3\textwidth}
        \centering
        \includegraphics[width=0.9\linewidth]{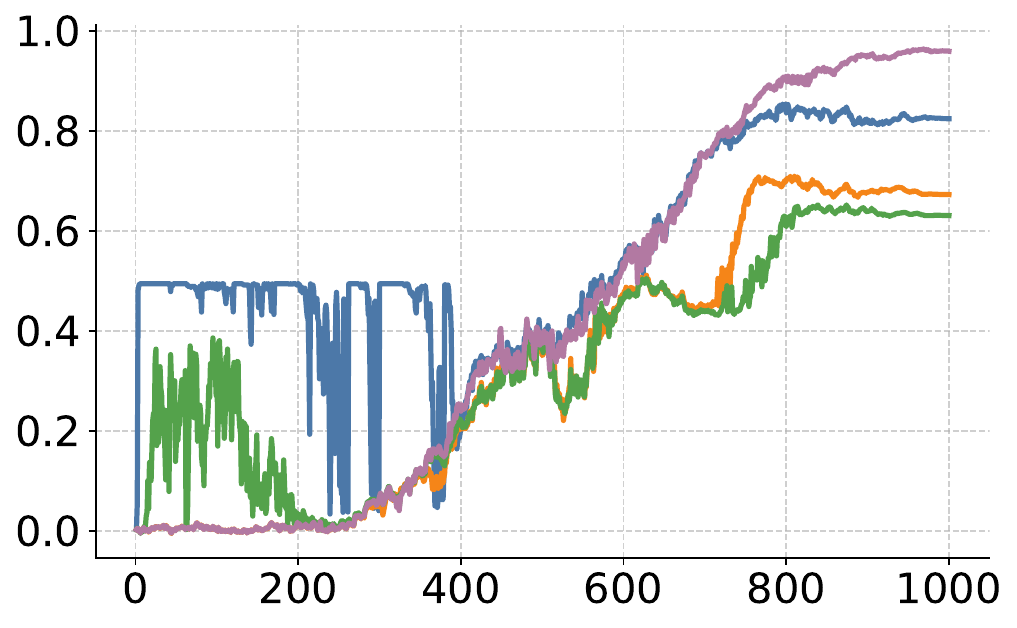}
        \caption{SSIM of $\hat v_0^{(i)}$ using $W_2$}
        \label{fig:abla_w2_ssim_error0}
    \end{subfigure}
    \hfill
    
    \begin{subfigure}[t]{0.3\textwidth}
        \centering
        \includegraphics[width=0.9\linewidth]{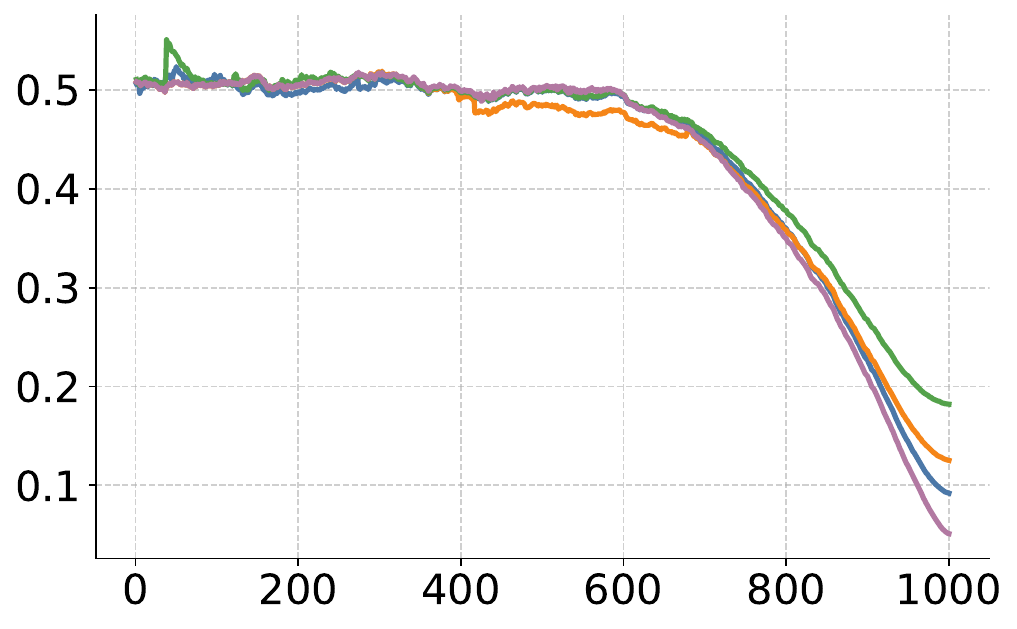}
        \caption{Relative $\ell_2$-error of $v_i$ using MSE}
        \label{fig:abla_MSE_el2_error}
    \end{subfigure}
    \hfill
    \begin{subfigure}[t]{0.3\textwidth}
        \centering
        \includegraphics[width=0.9\linewidth]{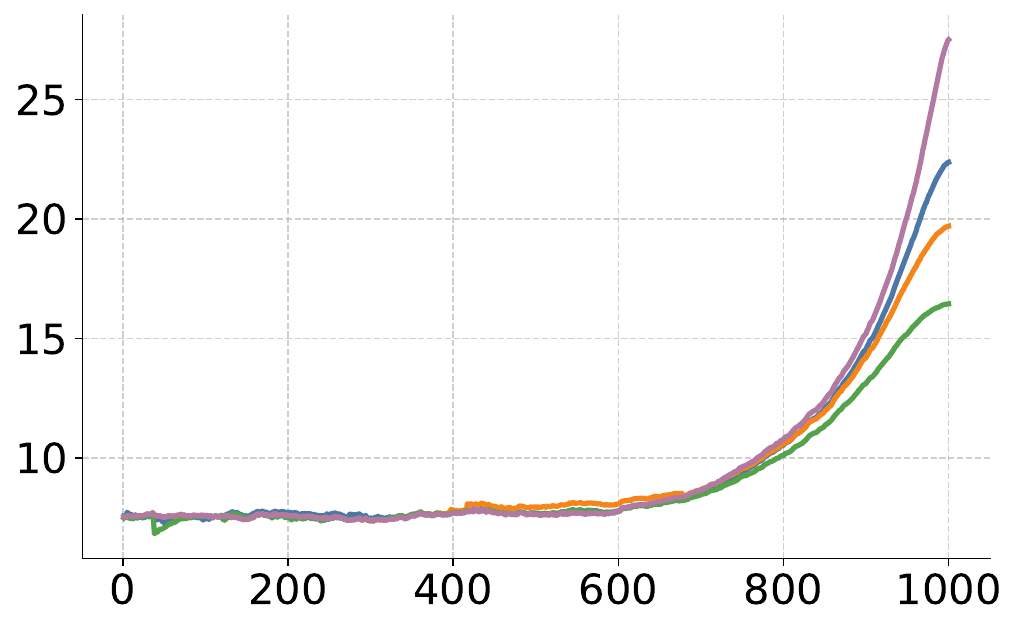}
        \caption{PSNR of $v_i$ using MSE}
        \label{fig:abla_MSE_psnr_error}
    \end{subfigure}
    \hfill
    \begin{subfigure}[t]{0.3\textwidth}
        \centering
        \includegraphics[width=0.9\linewidth]{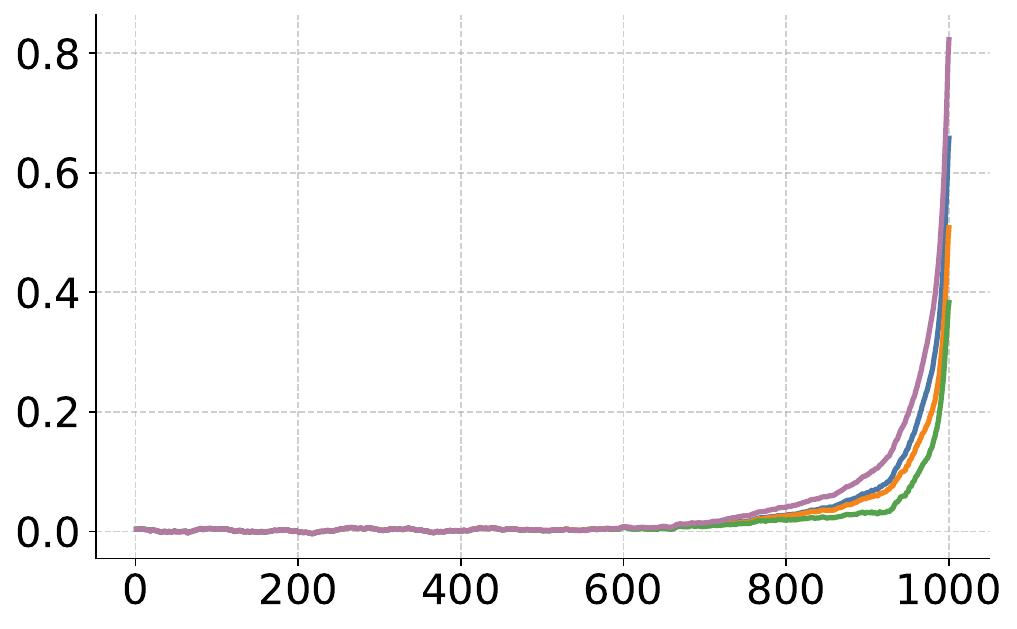}
        \caption{SSIM of $v_i$ using MSE}
        \label{fig:abla_MSE_ssim_error}
    \end{subfigure}
    \hfill

    \begin{subfigure}[t]{0.3\textwidth}
        \centering
        \includegraphics[width=0.9\linewidth]{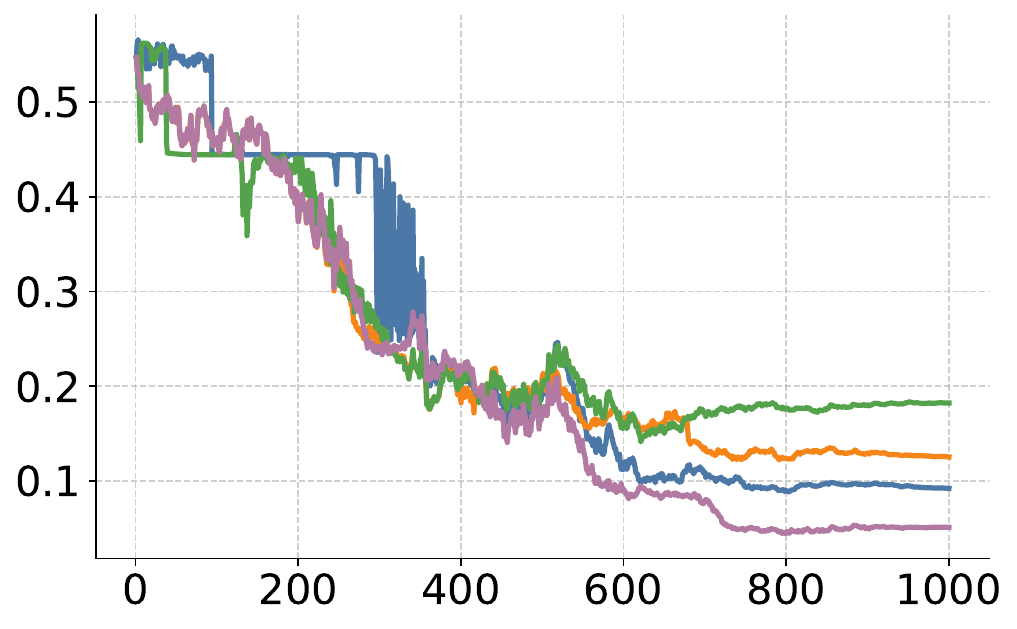}
        \caption{Relative $\ell_2$-error of $\hat v_0^{(i)}$ using MSE}
        \label{fig:abla_MSE_el2_error0}
    \end{subfigure}
    \hfill
    \begin{subfigure}[t]{0.3\textwidth}
        \centering
        \includegraphics[width=0.9\linewidth]{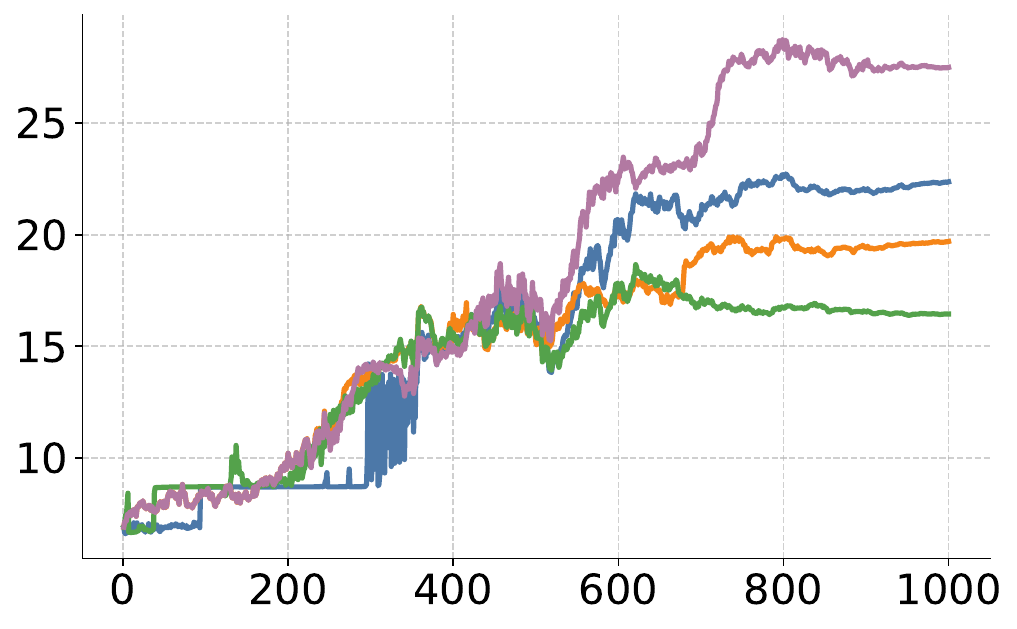}
        \caption{PSNR of $\hat v_0^{(i)}$ using MSE}
        \label{fig:abla_MSE_psnr_error0}
    \end{subfigure}
    \hfill
    \begin{subfigure}[t]{0.3\textwidth}
        \centering
        \includegraphics[width=0.9\linewidth]{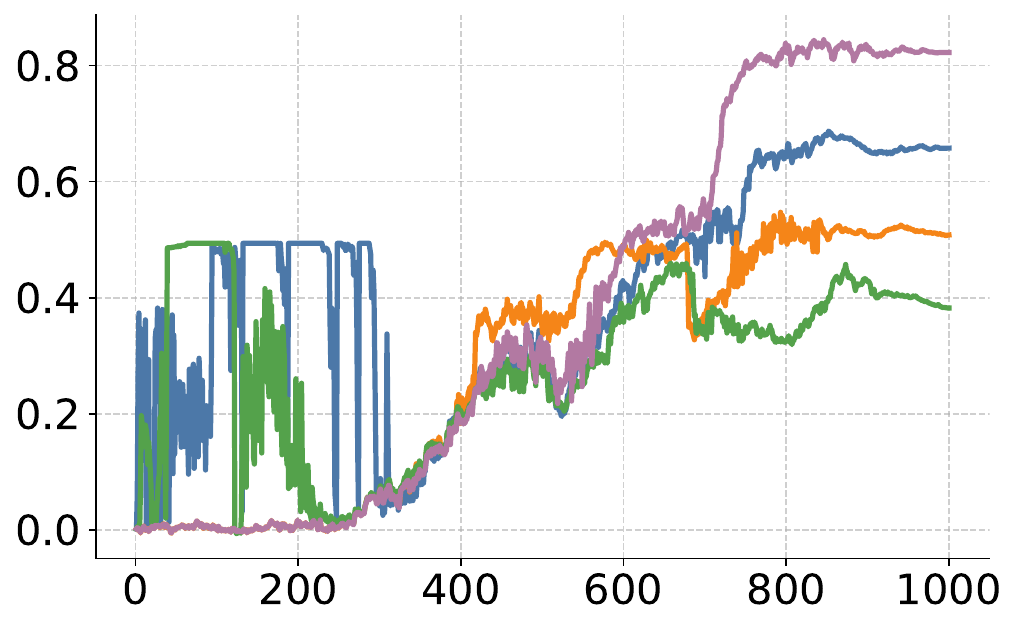}
        \caption{SSIM of $\hat v_0^{(i)}$ using MSE}
        \label{fig:abla_MSE_ssim_error0}
    \end{subfigure}
    \hfill

    \caption{Evolution of the statistical metrics over reverse-time steps for the ablation settings \texttt{Enh}, \texttt{Step}, \texttt{Only} and \texttt{All}. Panels (a)--(c) show the relative $\ell_2$-error, PSNR and SSIM of the intermediate iterate $v_i$ under $W_2$-based guidance; panels (d)--(f) show the corresponding metrics for the denoised estimate $\hat v_0^{(i)}$ under $W_2$-based guidance. Panels (g)--(i) and (j)--(l) report the analogous quantities under MSE-based guidance for $v_i$ and $\hat v_0^{(i)}$, respectively.}
    \label{fig:ablation_lines}
\end{figure}

The results are summarized in Figure~\ref{fig:abla} and Table~\ref{tab:ablation}.
Figure~\ref{fig:ablation_lines} further reports the evolution of the relative $\ell_2$ error (with respect to $v_{\mathrm{true}}$)
for the intermediate reconstructions $\hat{v}_0^{(i)}$ and $v_{i}$ over the reverse-diffusion steps $t$.
\textit{Notation.}
Throughout this ablation study, each label has the form \texttt{(misfit)-(component)}.
The prefix \texttt{$W_2$-} or \texttt{MSE-} specifies the guidance misfit used in the data-consistency term, namely the Wasserstein-$2$ distance ($W_2$) or the mean-squared error (MSE), respectively.
The suffix indicates which additional components are enabled:
\begin{itemize}
\item \texttt{Only}: baseline DPS with the chosen misfit, without any of our additional modifications.
\item \texttt{Enh}: baseline DPS plus \emph{wavefield enhancement}, i.e., amplitude-adaptive data weighting together with misfit-scale normalization (Sec.~\ref{subsec:ot_misfit}).
\item \texttt{Step}: baseline DPS plus the \emph{preconditioned guidance update} $P_i=\rho_iD_i$ (Sec.~\ref{subsec:precond_dps}), which combines the noise-aware scalar schedule $\rho_i$ and the spatial diagonal scaling $D_i$.
\item \texttt{All}: baseline DPS with both \texttt{Enh} and \texttt{Step} enabled.
\end{itemize}
For example, \texttt{$W_2$-Enh} uses $W_2$ guidance together with wavefield enhancement, while \texttt{MSE-All} enables both enhancements under MSE guidance.

\textit{Observations.}
Figures~\ref{fig:abla}--~\ref{fig:ablation_lines} and Table~\ref{tab:ablation} reveal the following trends:

(i) Combining wavefield enhancement with the preconditioned update yields the best overall performance.
In particular, \texttt{$W_2$-All} achieves $e_{\ell_2}=2.04\%$, PSNR $=35.45$, and SSIM $=0.9517$,
substantially improving over \texttt{$W_2$-Only} ($e_{\ell_2}=10.81\%$, PSNR $=20.98$, SSIM $=0.6307$).

(ii) Wavefield enhancement provides the larger gain when used alone, while the preconditioned update offers additional improvement when combined.
For example, under $W_2$ guidance, adding wavefield enhancement reduces $e_{\ell_2}$ from $10.81\%$ (\texttt{$W_2$-Only}) to $5.33\%$ (\texttt{$W_2$-Enh}),
whereas adding the preconditioned update alone reduces it to $8.65\%$ ($W_2$-Step).Enabling both components yields the best performance, achieving $e_{\ell_2}=2.04\%$ (\texttt{$W_2$-All}).

(iii) $W_2$-based guidance is consistently more robust than MSE-based guidance across all ablation settings, both in the baseline configuration (\texttt{Only}) and in the combined configuration (\texttt{All}). This is also reflected in the visual reconstructions and error maps in Figure~\ref{fig:abla}.

(iv) The metric trajectories suggest improved stability with preconditioned guidance.
As shown in Figure~\ref{fig:ablation_lines}, enabling $P_i=\rho_iD_i$ reduces transient oscillations during intermediate reverse steps
and leads to more stable convergence, with the combined setting (\texttt{All}) achieving the best final metrics.

\subsection{Robustness Study}
\begin{figure}[th]
  \centering
  % \newlength{\subfigwidth}
  \setlength{\subfigwidth}{2.1cm}
  % \newlength{\insetfigwidth}
  \setlength{\insetfigwidth}{8cm}

  \begin{tikzpicture}[
      subfig/.style={anchor=center, inner sep=0pt, outer sep=2pt, draw=black},
      subfig matrix/.style={
        matrix of nodes,
        column sep=0pt,    % 基础列间距
        row sep=0pt,       % 基础行间距
        nodes={subfig, minimum width=\subfigwidth, minimum height=\subfigwidth}
      },
      % 定义虚线样式
      split line/.style={dashed, gray, line width=0.5pt},
      label text/.style={anchor=east, font=\small, inner sep=5pt, align=center}
    ]

    \matrix (subfig_matrix) [subfig matrix
      ,column 1/.style={column sep=5pt}
      ,column 2/.style={column sep=5pt}
      ,column 3/.style={column sep=5pt}
      ,column 4/.style={column sep=5pt}
      ,column 5/.style={column sep=5pt}
      ,row 1/.style={row sep=8pt}
    ] {
      % 第一行：前4列 + 空白列 + 后4列
      \includegraphics[width=\subfigwidth, height=\subfigwidth, keepaspectratio]{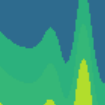} &
      \includegraphics[width=\subfigwidth, height=\subfigwidth, keepaspectratio]{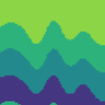} &
      \includegraphics[width=\subfigwidth, height=\subfigwidth, keepaspectratio]{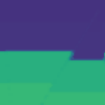} &
      \includegraphics[width=\subfigwidth, height=\subfigwidth, keepaspectratio]{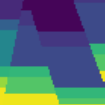} &
      \includegraphics[width=\subfigwidth, height=\subfigwidth, keepaspectratio]{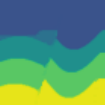} &
      \includegraphics[width=\subfigwidth, height=\subfigwidth, keepaspectratio]{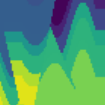}\\

      \includegraphics[width=\subfigwidth, height=\subfigwidth, keepaspectratio]{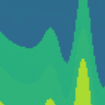} &
      \includegraphics[width=\subfigwidth, height=\subfigwidth, keepaspectratio]{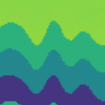} &
      \includegraphics[width=\subfigwidth, height=\subfigwidth, keepaspectratio]{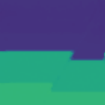} &
      \includegraphics[width=\subfigwidth, height=\subfigwidth, keepaspectratio]{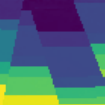} &
      \includegraphics[width=\subfigwidth, height=\subfigwidth, keepaspectratio]{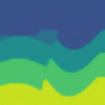} &
      \includegraphics[width=\subfigwidth, height=\subfigwidth, keepaspectratio]{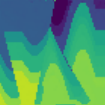}\\
      
      \includegraphics[width=\subfigwidth, height=\subfigwidth, keepaspectratio]{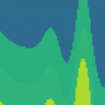} &
      \includegraphics[width=\subfigwidth, height=\subfigwidth, keepaspectratio]{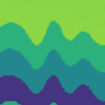} &
      \includegraphics[width=\subfigwidth, height=\subfigwidth, keepaspectratio]{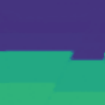} &
      \includegraphics[width=\subfigwidth, height=\subfigwidth, keepaspectratio]{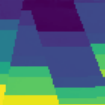} &
      \includegraphics[width=\subfigwidth, height=\subfigwidth, keepaspectratio]{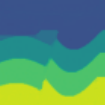} &
      \includegraphics[width=\subfigwidth, height=\subfigwidth, keepaspectratio]{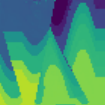}\\
      
      \includegraphics[width=\subfigwidth, height=\subfigwidth, keepaspectratio]{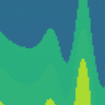} &
      \includegraphics[width=\subfigwidth, height=\subfigwidth, keepaspectratio]{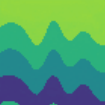} &
      \includegraphics[width=\subfigwidth, height=\subfigwidth, keepaspectratio]{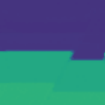} &
      \includegraphics[width=\subfigwidth, height=\subfigwidth, keepaspectratio]{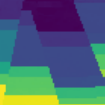} &
      \includegraphics[width=\subfigwidth, height=\subfigwidth, keepaspectratio]{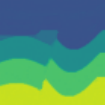} &
      \includegraphics[width=\subfigwidth, height=\subfigwidth, keepaspectratio]{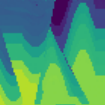}\\

      \includegraphics[width=\subfigwidth, height=\subfigwidth, keepaspectratio]{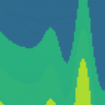} &
      \includegraphics[width=\subfigwidth, height=\subfigwidth, keepaspectratio]{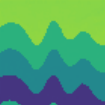} &
      \includegraphics[width=\subfigwidth, height=\subfigwidth, keepaspectratio]{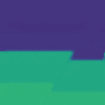} &
      \includegraphics[width=\subfigwidth, height=\subfigwidth, keepaspectratio]{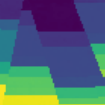} &
      \includegraphics[width=\subfigwidth, height=\subfigwidth, keepaspectratio]{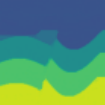} &
      \includegraphics[width=\subfigwidth, height=\subfigwidth, keepaspectratio]{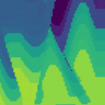}\\

      \includegraphics[width=\subfigwidth, height=\subfigwidth, keepaspectratio]{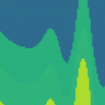} &
      \includegraphics[width=\subfigwidth, height=\subfigwidth, keepaspectratio]{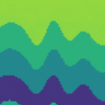} &
      \includegraphics[width=\subfigwidth, height=\subfigwidth, keepaspectratio]{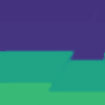} &
      \includegraphics[width=\subfigwidth, height=\subfigwidth, keepaspectratio]{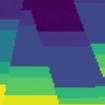} &
      \includegraphics[width=\subfigwidth, height=\subfigwidth, keepaspectratio]{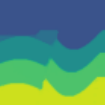} &
      \includegraphics[width=\subfigwidth, height=\subfigwidth, keepaspectratio]{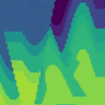}\\
    };

    \coordinate (inset_pos) at ($(subfig_matrix-1-6.east)!0.5!(subfig_matrix-6-6.east)+(25pt,0)$);
    
    % 插入图片（可替换为你的图片路径）
    \node[minimum width=\insetfigwidth, minimum height=\insetfigwidth] at (inset_pos) {
      \includegraphics[width=\insetfigwidth, height=\insetfigwidth, keepaspectratio]{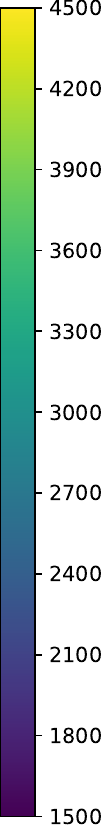}
    };
    \node[right=13pt, font=\small] at (inset_pos) {\rotatebox{90}{velocity (m/s)}};

    \draw[split line] 
      ($(subfig_matrix-1-1.west)!0.5!(subfig_matrix-2-1.west)+(-7pt,0)$) -- ($(subfig_matrix-1-6.east)!0.5!(subfig_matrix-2-6.east)+(7pt,0)$);

    % 顶部标签
    \node[above=28pt, font=\small] at (subfig_matrix-1-1) {CurveVel-A};
    \node[above=28pt, font=\small] at (subfig_matrix-1-2) {CurveVel-B};
    \node[above=28pt, font=\small] at (subfig_matrix-1-3) {FlatFault-A};
    \node[above=28pt, font=\small] at (subfig_matrix-1-4) {FlatFault-B};
    \node[above=28pt, font=\small] at (subfig_matrix-1-5) {CurveFault-A};
    \node[above=28pt, font=\small] at (subfig_matrix-1-6) {CurveFault-B};

    \node[left=28pt, font=\small] at (subfig_matrix-1-1) {\rotatebox{90}{Truth}};
    \node[left=28pt, font=\small] at (subfig_matrix-2-1) {\rotatebox{90}{$\sigma=0$}};
    \node[left=28pt, font=\small] at (subfig_matrix-3-1) {\rotatebox{90}{$\sigma=1e-4$}};
    \node[left=28pt, font=\small] at (subfig_matrix-4-1) {\rotatebox{90}{$\sigma=1e-3$}};
    \node[left=28pt, font=\small] at (subfig_matrix-5-1) {\rotatebox{90}{$\sigma=1e-2$}};
    \node[left=28pt, font=\small] at (subfig_matrix-6-1) {\rotatebox{90}{$\sigma=5e-2$}};

    % 底部标注（修正坐标，适配新的列结构）
    % \node[label text, anchor=north, inner sep=5pt, font=\large,align=center] at ($(subfig_matrix-2-1)!0.5!(subfig_matrix-3-1)+(-34pt, 22pt)$) {\rotatebox{90}{D   D   I   M}};

  \end{tikzpicture}
  \caption{Qualitative robustness to observation noise across the six datasets. Each column corresponds to one dataset. The first row shows the true velocity field, and the remaining rows show the reconstructed velocity fields under progressively increasing noise levels $\sigma\in\{0,10^{-4},10^{-3},10^{-2},5\times10^{-2}\}$. The reconstructions remain visually stable as the noise level increases, with the main interfaces and fault structures largely preserved.}
  \label{fig:robust_noise}
\end{figure}
\subsubsection{Robustness to Observation Noise}\label{subsec:robust_noise}
We evaluate the robustness of our method to observation noise by corrupting the recorded data with additive noise at increasing levels
$\sigma\in\{0,10^{-4},10^{-3},10^{-2},5\times 10^{-2}\}$.
All other hyperparameters are kept identical to those used in the previous experiments.
Qualitative reconstructions are shown in Figure~\ref{fig:robust_noise}, and quantitative results, including the relative $\ell_2$-error, PSNR, and SSIM, are reported in Table~\ref{tab:robust_noise}.

Overall, the reconstructions remain visually stable across all noise levels for all six datasets: the major interfaces and fault geometries are preserved, and the predicted velocity ranges remain consistent (Figure~\ref{fig:robust_noise}).
Quantitatively, the performance degrades only mildly as $\sigma$ increases, with nonmonotone but limited variations in $e_{\ell_2}$, PSNR, and SSIM (Table~\ref{tab:robust_noise}).
These results indicate that the proposed guidance remains effective under moderate levels of observation noise.

\begin{table}[t]
\caption{Quantitative robustness to observation noise across the six datasets. Relative $\ell_2$-error ($e_{\ell_2}$; lower is better), PSNR, and SSIM (higher is better) of the reconstructed velocity fields under different noise levels $\sigma$.}
\label{tab:robust_noise}
\centering
\small
\setlength{\tabcolsep}{4.2pt}
\renewcommand{\arraystretch}{1.12}

\begin{subtable}{0.98\linewidth}
\centering
\begin{tabular}{c ccc ccc ccc}
\toprule
& \multicolumn{3}{c}{CurveVel-A} & \multicolumn{3}{c}{CurveVel-B} & \multicolumn{3}{c}{FlatFault-A} \\
\cmidrule(lr){2-4}\cmidrule(lr){5-7}\cmidrule(lr){8-10}
$\sigma$ & $e_{\ell_2}\downarrow$ & PSNR$\uparrow$ & SSIM$\uparrow$
        & $e_{\ell_2}\downarrow$ & PSNR$\uparrow$ & SSIM$\uparrow$
        & $e_{\ell_2}\downarrow$ & PSNR$\uparrow$ & SSIM$\uparrow$ \\
\midrule
$0$          & 3.29\% & 31.96 & 0.9111 & 1.95\% & 35.83 & 0.9584 & 1.60\% & 37.94 & 0.9820 \\
$10^{-4}$    & 3.29\% & 31.96 & 0.9113 & 1.94\% & 35.90 & 0.9594 & 1.84\% & 36.74 & 0.9691 \\
$10^{-3}$    & 3.39\% & 31.70 & 0.9075 & 2.23\% & 34.69 & 0.9451 & 1.95\% & 36.20 & 0.9739 \\
$10^{-2}$    & 3.62\% & 31.12 & 0.9121 & 2.07\% & 35.34 & 0.9557 & 1.97\% & 36.15 & 0.9747 \\
$5\!\times\!10^{-2}$ & 3.50\% & 31.42 & 0.9039 & 2.04\% & 35.45 & 0.9517 & 2.91\% & 32.74 & 0.9646 \\
\bottomrule
\end{tabular}
\end{subtable}

\vspace{0.6em}

\begin{subtable}{0.98\linewidth}
\centering
\begin{tabular}{c ccc ccc ccc}
\toprule
& \multicolumn{3}{c}{FlatFault-B} & \multicolumn{3}{c}{CurveFault-A} & \multicolumn{3}{c}{CurveFault-B} \\
\cmidrule(lr){2-4}\cmidrule(lr){5-7}\cmidrule(lr){8-10}
$\sigma$ & $e_{\ell_2}\downarrow$ & PSNR$\uparrow$ & SSIM$\uparrow$
        & $e_{\ell_2}\downarrow$ & PSNR$\uparrow$ & SSIM$\uparrow$
        & $e_{\ell_2}\downarrow$ & PSNR$\uparrow$ & SSIM$\uparrow$ \\
\midrule
$0$          & 6.07\% & 28.56 & 0.9242 & 2.41\% & 35.01 & 0.9602 & 5.67\% & 27.39 & 0.7920 \\
$10^{-4}$    & 6.03\% & 28.61 & 0.9240 & 2.22\% & 35.75 & 0.9606 & 6.03\% & 26.86 & 0.7779 \\
$10^{-3}$    & 6.10\% & 28.51 & 0.9225 & 2.45\% & 34.90 & 0.9552 & 6.34\% & 26.43 & 0.7407 \\
$10^{-2}$    & 6.15\% & 28.45 & 0.9214 & 2.26\% & 35.56 & 0.9609 & 6.65\% & 26.01 & 0.7304 \\
$5\!\times\!10^{-2}$ & 6.13\% & 28.48 & 0.9272 & 2.41\% & 35.03 & 0.9540 & 5.32\% & 27.95 & 0.7956 \\
\bottomrule
\end{tabular}
\end{subtable}

\end{table}

\newcolumntype{C}[1]{>{\centering\arraybackslash}m{#1}}

% \begin{table}[]
%     \centering
%     \begin{tabular}{@{}ccccccc@{}}

%     \toprule
%     $\sigma$ & CurveVel-A & CurveVel-B & FlatFault-A & FlatFault-B & CurveFault-A & CurveFault-B \\
%     \midrule
%     0 & 0 & 1 & 2 & 3 & 4 & 5 \\
%     1e-4 & & & & & & \\
%     1e-3 & & & & & & \\
%     1e-2 & & & & & & \\
%     5e-2 & & & & & & \\
%     \bottomrule
%     \end{tabular}
%     \caption{Relative errors in robust study}
%     \label{tab:robust}
% \end{table}

\subsection{Generalization across forward operators}\label{subsec:multi_operator}

An important feature of the proposed framework is that the learned score-based prior is trained only on prior velocity models and does not depend on a specific forward operator.
Once this prior has been trained, it can be coupled with different forward models $\mathcal{F}$ for inversion without retraining the prior network.
In this subsection, we examine the robustness of the proposed method with respect to changes in the forward operator by perturbing the acquisition and physical settings and then re-running the inversion.

We consider the FlatFault-B dataset and define five forward-operator configurations.
The baseline configuration, denoted by \texttt{Init}, coincides with that used in the preceding experiments.
We then modify the acquisition and physical settings in several ways, thereby obtaining perturbed forward operators as follows:
\begin{enumerate}
\item \texttt{Freq10}: we replace the temporal source function by a Ricker wavelet of peak frequency $f_p=10\,\mathrm{Hz}$, while keeping the spatial point-source locations, all other parameters, and the acquisition geometry unchanged. Consequently, the observed and synthetic wavefields are generated with a different source peak frequency.
  \item \texttt{Depth2}: we shift the acquisition array downward by two grid intervals by setting the source and receiver depths to $z=L_z-2\,dx$ instead of $z=L_z$, while keeping their horizontal coordinates fixed. This modification changes how waves propagate from the sources to the receivers and how the subsurface structures are illuminated.
 \item \texttt{Shot17}: relative to the baseline configuration with 10 sources, we increase the number of sources to 17 and distribute them uniformly along the surface, while keeping the receiver layout unchanged. This perturbation yields denser illumination and increased data redundancy for the inversion.
\item \texttt{AllChange}: we apply \texttt{Freq10}, \texttt{Depth2}, and \texttt{Shot17} simultaneously, thereby introducing a compound perturbation of the forward operator that modifies the source peak frequency, the source and receiver depths, and the number of sources.
\end{enumerate}

\begin{figure}[th]
  \centering
  
  \setlength{\subfigwidth}{2.6cm}
  
  \setlength{\insetfigwidth}{2.7cm}

  \begin{tikzpicture}[
      subfig/.style={anchor=center, inner sep=0pt, outer sep=2pt, draw=black},
      subfig matrix/.style={
        matrix of nodes,
        column sep=0pt,    % 基础列间距
        row sep=0pt,       % 基础行间距
        nodes={subfig, minimum width=\subfigwidth, minimum height=\subfigwidth}
      },
      % 定义虚线样式
      no border/.style={subfig, draw=none},
      split line/.style={dashed, gray, line width=0.5pt},
      label text/.style={anchor=east, font=\small, inner sep=5pt, align=center}
    ]

    \matrix (subfig_matrix) [subfig matrix
      ,row 1/.style={row sep=18pt}
      ,row 2/.style={row sep=8pt}
    ] {
      % 第一行：前4列 + 空白列 + 后4列
      \includegraphics[width=\subfigwidth, height=\subfigwidth, keepaspectratio]{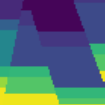}\\
      \includegraphics[width=\subfigwidth, height=\subfigwidth, keepaspectratio]{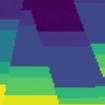} &
      \includegraphics[width=\subfigwidth, height=\subfigwidth, keepaspectratio]{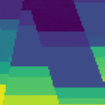} &
      \includegraphics[width=\subfigwidth, height=\subfigwidth, keepaspectratio]{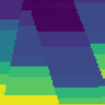} &
      \includegraphics[width=\subfigwidth, height=\subfigwidth, keepaspectratio]{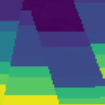} &
      \includegraphics[width=\subfigwidth, height=\subfigwidth, keepaspectratio]{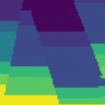}  \\

      \includegraphics[width=\subfigwidth, height=\subfigwidth, keepaspectratio]{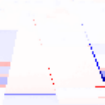} &
      \includegraphics[width=\subfigwidth, height=\subfigwidth, keepaspectratio]{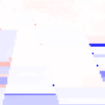} &
      \includegraphics[width=\subfigwidth, height=\subfigwidth, keepaspectratio]{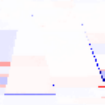} &
      \includegraphics[width=\subfigwidth, height=\subfigwidth, keepaspectratio]{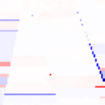} &
      \includegraphics[width=\subfigwidth, height=\subfigwidth, keepaspectratio]{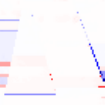}  \\
    };

    \coordinate (inset_pos1) at ($(subfig_matrix-2-5.east)+(35pt,0)$);
    \coordinate (inset_pos2) at ($(subfig_matrix-3-5.east)+(35pt,0)$);

    \node[minimum width=\insetfigwidth, minimum height=\insetfigwidth] at (inset_pos1) {
      \includegraphics[width=\insetfigwidth, height=\insetfigwidth, keepaspectratio]{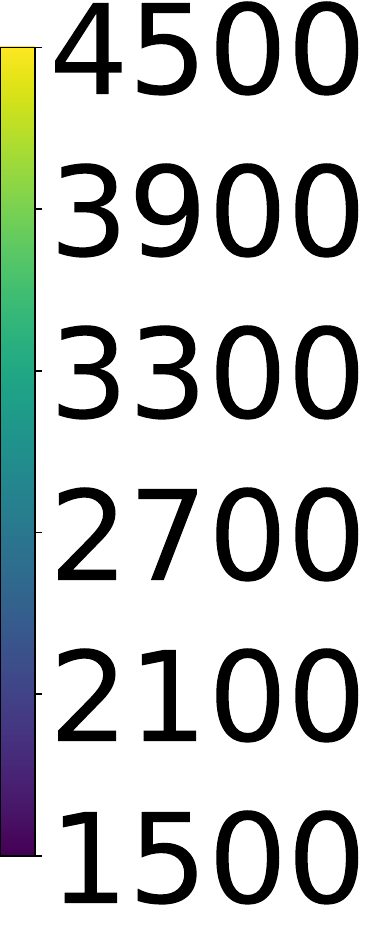}
    };
    \node[right=15pt, font=\small] at (inset_pos1) {\rotatebox{90}{velocity (m/s)}};

    \node[minimum width=\insetfigwidth, minimum height=\insetfigwidth] at (inset_pos2) {
      \includegraphics[width=\insetfigwidth, height=\insetfigwidth, keepaspectratio]{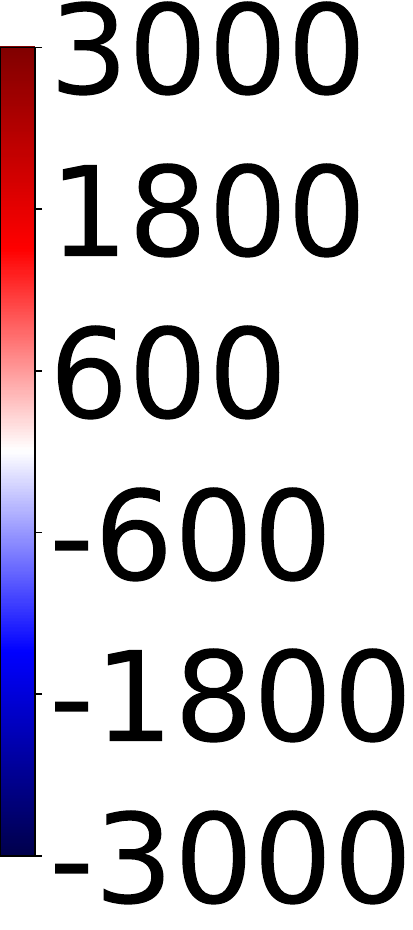}
    };
    \node[right=15pt, font=\small] at (inset_pos2) {\rotatebox{90}{velocity (m/s)}};

    \draw[split line] 
      ($(subfig_matrix-2-1.west)!0.5!(subfig_matrix-3-1.west)+(-7pt,0)$) -- ($(subfig_matrix-2-5.east)!0.5!(subfig_matrix-3-5.east)+(7pt,0)$);
    
    \node[font=\small, anchor=center] at (
      $(subfig_matrix-2-5.east)+(8pt,0)$
    ) {\rotatebox{90}{Reconstruction}};

    \node[font=\small, anchor=center] at (
      $(subfig_matrix-3-5.east)+(8pt,0)$
    ) {\rotatebox{90}{Difference}};
    
    % 顶部标签
    \node[above=38pt, font=\small] at (subfig_matrix-1-1) {\texttt{Truth}};
    \node[above=38pt, font=\small] at (subfig_matrix-2-1) {\texttt{Init}};
    \node[above=38pt, font=\small] at (subfig_matrix-2-2) {\texttt{Shot17}};
    \node[above=36pt, font=\small] at (subfig_matrix-2-3) {\texttt{Depth2}};
    \node[above=36pt, font=\small] at (subfig_matrix-2-4) {\texttt{Freq10}};
    \node[above=38pt, font=\small] at (subfig_matrix-2-5) {\texttt{AllChange}};

  \end{tikzpicture}
  \caption{Inversion results under different forward-operator configurations (FlatFault-B). From left to right: \texttt{Init} (baseline setup), \texttt{Shot17} (17 uniformly distributed sources), \texttt{Depth2} (sources and receivers shifted downward by $2\,dx$), \texttt{Freq10} (Ricker wavelet peak frequency $f_p=10,\mathrm{Hz}$), and \texttt{AllChange} (all modifications applied). The reconstructions remain visually consistent across operator changes, indicating that the proposed method generalizes well to perturbed forward operators without retraining the prior.}
  \label{fig:multi_oper}
\end{figure}
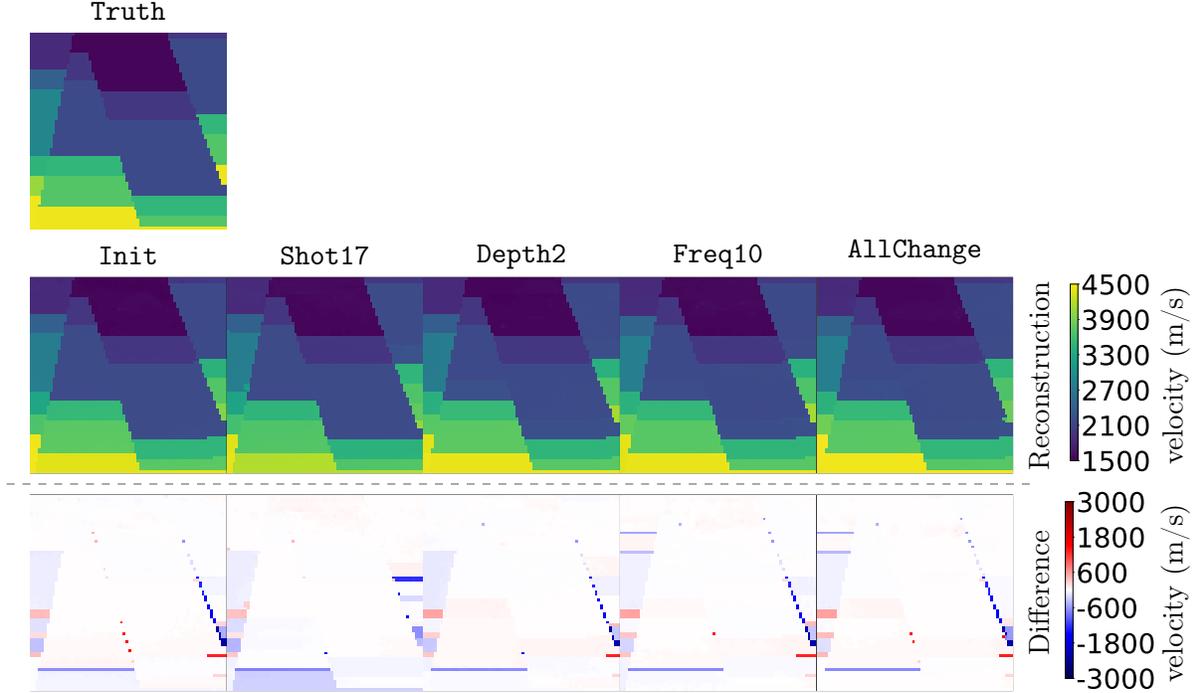

For a fair comparison across the different operator settings, the network architecture, diffusion schedule, and all remaining hyperparameters are kept identical to those used in the preceding experiments.
For each operator configuration, only the guidance-scale parameter $\rho_0$ is tuned, following the same protocol, so as to achieve the best performance attainable for our method; the selected values are listed in Table~\ref{tab:multi_oper}.
The corresponding quantitative metrics and qualitative reconstructions are reported in Table~\ref{tab:multi_oper} and Figure~\ref{fig:multi_oper}, respectively.›

\begin{table}[th]
  \caption{Reconstruction performance on FlatFault-B under different forward-operator configurations. Reported are the tuned values of $\rho_0$, together with the relative $\ell_2$-error ($e_{\ell_2}$; lower is better), PSNR, and SSIM (higher is better).}
  \label{tab:multi_oper}
  \centering
  \small
  \setlength{\tabcolsep}{6pt}
  \renewcommand{\arraystretch}{1.10}
  \begin{tabular}{lcccc}
    \toprule
    Operator & $\rho_0$ & $e_{\ell_2}\downarrow$ & PSNR$\uparrow$ & SSIM$\uparrow$ \\
    \midrule
    \texttt{Init}      & 1.45 & 6.13\% & 28.48 & 0.9272 \\
    \texttt{Shot17}    & 2.25 & 5.24\% & 29.84 & 0.9399 \\
    \texttt{Depth2}    & 1.60 & 5.37\% & 29.62 & 0.9343 \\
    \texttt{Freq10}    & 1.00 & 5.85\% & 28.87 & 0.9160 \\
    \texttt{AllChange} & 1.30 & 5.82\% & 28.92 & 0.9168 \\
    \bottomrule
  \end{tabular}
\end{table}

It can be observed that, across all forward-operator configurations, our method consistently yields high-quality reconstructions, with SSIM exceeding $0.9$ and relative errors below $6.2\%$ (Table~\ref{tab:multi_oper}). These results indicate that, once the score-based prior is trained on velocity models, it can be coupled with modified acquisition/physics settings (e.g., changes in source frequency, source/receiver depth, or the number of shots) without retraining the prior. This operator-agnostic property substantially broadens the applicability of the proposed inversion framework.

\subsection{A unified diffusion model for heterogeneous velocity fields}

\subsubsection{Performance of the mixed model on OpenFWI}
In the preceding experiments, a separate diffusion model was trained for each family of velocity fields. Although this dataset-specific strategy is effective when the structural class is known a priori, it is of limited practical value for inverse problems, where the subsurface structure is typically unknown. To address this issue, we next consider a \emph{mixed model} trained jointly on a heterogeneous dataset assembled from multiple structural classes. The resulting unified model is designed to learn a broader prior and thus avoids the need to select a class-specific model at the inference stage.

\begin{figure}[th]
  \centering
  
  \setlength{\subfigwidth}{2.6cm}
  
  \setlength{\insetfigwidth}{3.5cm}

  \begin{tikzpicture}[
      subfig/.style={anchor=center, inner sep=0pt, outer sep=2pt, draw=black},
      subfig matrix/.style={
        matrix of nodes,
        column sep=0pt,    % 基础列间距
        row sep=0pt,       % 基础行间距
        nodes={subfig, minimum width=\subfigwidth, minimum height=\subfigwidth}
      },
      % 定义虚线样式
      no border/.style={subfig, draw=none},
      split line/.style={dashed, gray, line width=0.5pt},
      label text/.style={anchor=east, font=\small, inner sep=5pt, align=center}
    ]

    \matrix (subfig_matrix) [subfig matrix
      ,row 1/.style={row sep=8pt}
      ,row 3/.style={row sep=8pt}
      ,column 1/.style={column sep=5pt}
      ,column 2/.style={column sep=5pt}
    ] {
      % 第一行：前4列 + 空白列 + 后4列
      \includegraphics[width=\subfigwidth, height=\subfigwidth, keepaspectratio]{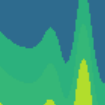} & 
      \includegraphics[width=\subfigwidth, height=\subfigwidth, keepaspectratio]{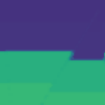} & 
      \includegraphics[width=\subfigwidth, height=\subfigwidth, keepaspectratio]{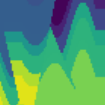} \\
      
      \includegraphics[width=\subfigwidth, height=\subfigwidth, keepaspectratio]{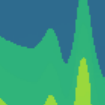} &
      \includegraphics[width=\subfigwidth, height=\subfigwidth, keepaspectratio]{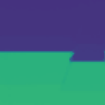} &
      \includegraphics[width=\subfigwidth, height=\subfigwidth, keepaspectratio]{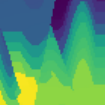} \\

      \includegraphics[width=\subfigwidth, height=\subfigwidth, keepaspectratio]{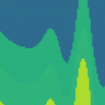} &
      \includegraphics[width=\subfigwidth, height=\subfigwidth, keepaspectratio]{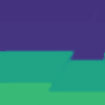} &
      \includegraphics[width=\subfigwidth, height=\subfigwidth, keepaspectratio]{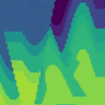} \\

      \includegraphics[width=\subfigwidth, height=\subfigwidth, keepaspectratio]{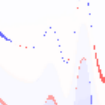} &
      \includegraphics[width=\subfigwidth, height=\subfigwidth, keepaspectratio]{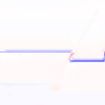} &
      \includegraphics[width=\subfigwidth, height=\subfigwidth, keepaspectratio]{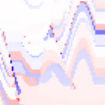} \\

      \includegraphics[width=\subfigwidth, height=\subfigwidth, keepaspectratio]{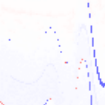} &
      \includegraphics[width=\subfigwidth, height=\subfigwidth, keepaspectratio]{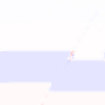} &
      \includegraphics[width=\subfigwidth, height=\subfigwidth, keepaspectratio]{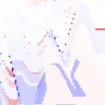} \\
    };

    % \coordinate (inset_pos0) at ($(subfig_matrix-1-3.east)+(35pt,0)$);
    \coordinate (inset_pos1) at ($(subfig_matrix-2-3.east)!0.5!(subfig_matrix-3-3.east)+(40pt,0)$);
    \coordinate (inset_pos2) at ($(subfig_matrix-4-3.east)!0.5!(subfig_matrix-5-3.east)+(40pt,0)$);

    % \node[minimum width=\insetfigwidth, minimum height=\insetfigwidth] at (inset_pos0) {
    %   \includegraphics[width=\insetfigwidth, height=\insetfigwidth, keepaspectratio]{picture/marmousi_colorbar_1545.pdf}
    % };
    % \node[right=15pt, font=\small] at (inset_pos1) {\rotatebox{90}{velocity (m/s)}};

    \node[minimum width=\insetfigwidth, minimum height=\insetfigwidth] at (inset_pos1) {
      \includegraphics[width=\insetfigwidth, height=\insetfigwidth, keepaspectratio]{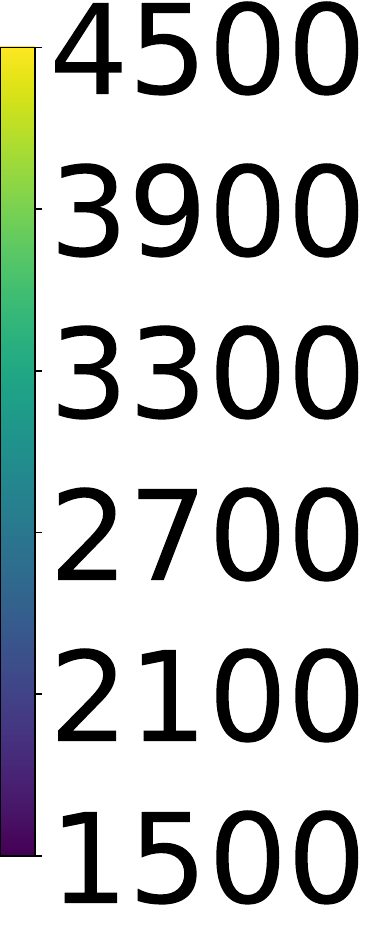}
    };
    \node[right=18pt, font=\small] at (inset_pos1) {\rotatebox{90}{velocity (m/s)}};

    \node[minimum width=\insetfigwidth, minimum height=\insetfigwidth] at (inset_pos2) {
      \includegraphics[width=\insetfigwidth, height=\insetfigwidth, keepaspectratio]{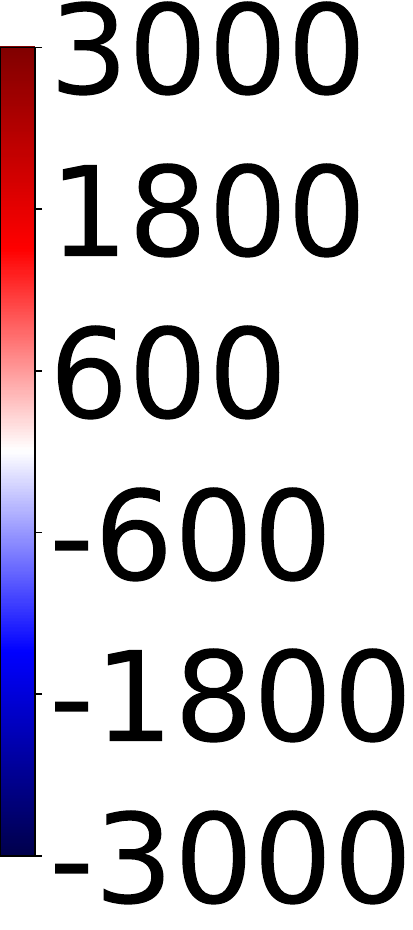}
    };
    \node[right=18pt, font=\small] at (inset_pos2) {\rotatebox{90}{velocity (m/s)}};

    \draw[split line] 
      ($(subfig_matrix-3-1.west)!0.5!(subfig_matrix-4-1.west)+(-7pt,0)$) -- ($(subfig_matrix-3-3.east)!0.5!(subfig_matrix-4-3.east)+(7pt,0)$);

    \draw[split line] 
      ($(subfig_matrix-1-1.west)!0.5!(subfig_matrix-2-1.west)+(-7pt,0)$) -- ($(subfig_matrix-1-3.east)!0.5!(subfig_matrix-2-3.east)+(7pt,0)$);

    % \node[font=\small, anchor=center] at (
    %   $(subfig_matrix-1-3.east)+(8pt,0)$
    % ) {\rotatebox{90}{Truth}};
    
    \node[font=\small, anchor=center] at (
      $(subfig_matrix-2-3.east)!0.5!(subfig_matrix-3-3.east)+(8pt,0)$
    ) {\rotatebox{90}{Reconstruction}};

    \node[font=\small, anchor=center] at (
      $(subfig_matrix-4-3.east)!0.5!(subfig_matrix-5-3.east)+(8pt,0)$
    ) {\rotatebox{90}{Difference}};
    
    % 顶部标签
    \node[above=38pt, font=\small] at (subfig_matrix-1-1) {\texttt{CurveVel-A}};
    \node[above=38pt, font=\small] at (subfig_matrix-1-2) {\texttt{FlatFault-A}};
    \node[above=38pt, font=\small] at (subfig_matrix-1-3) {\texttt{CurveFault-B}};

    \node[left=38pt, font=\small] at (subfig_matrix-1-1) {\rotatebox{90}{Truth}};
    \node[left=38pt, font=\small] at (subfig_matrix-2-1) {\rotatebox{90}{mixed model}};
    \node[left=36pt, font=\small] at (subfig_matrix-3-1) {\rotatebox{90}{separate model}};
    \node[left=38pt, font=\small] at (subfig_matrix-4-1) {\rotatebox{90}{mixed model}};
    \node[left=36pt, font=\small] at (subfig_matrix-5-1) {\rotatebox{90}{separate model}};

  \end{tikzpicture}
  \caption{Comparison of inversion results produced by the mixed and separate models for CurveVel-A, FlatFault-A, and CurveFault-B. The first row shows the true velocity models $v_{\mathrm{true}}$. The second and third rows display the reconstructed velocity models $v_{\mathrm{rec}}$ obtained by the mixed model and the separate model, respectively. The fourth and fifth rows show the reconstruction errors, $v_{\mathrm{rec}}-v_{\mathrm{true}}$, corresponding to the mixed and separate models, respectively.}
  \label{fig:mix}
\end{figure}

Specifically, we randomly selected 7,500 samples from each of CurveVel-A, CurveVel-B, FlatFault-A, FlatFault-B, CurveFault-A, and CurveFault-B, and merged them into a single training set. The mixed model was then evaluated on several test datasets, among which CurveVel-A, FlatFault-A, and CurveFault-B are reported here as representative examples. We adopted the same operator setting as in Section~\ref{subsubsec:operator_settings} and carefully tuned the remaining hyperparameters. More precisely, we set
$k=100$, $c=0.1$, $\tau=0$, $\varepsilon=10^{-4}$, and $\gamma=0.55$,
and chose $\rho_0=0.95$, $0.85$, and $1.4$ for CurveVel-A, FlatFault-A, and CurveFault-B, respectively. For each dataset, the performance of the mixed model is compared with that of the corresponding separate model trained exclusively on the same dataset. The quantitative and qualitative results are reported in Table~\ref{tab:mix} and Figure~\ref{fig:mix}, respectively.

It can be observed that, when the score-based neural network is trained on the mixed dataset rather than on a single dataset, the inversion performance degrades slightly across all three representative examples. Nevertheless, the mixed model still yields satisfactory reconstructions and remains clearly competitive in practice. These results indicate that, when strong prior knowledge of the target velocity field is available, a dataset-specific model trained on structurally homogeneous samples is generally preferable. By contrast, when prior information is limited or the structural type is unknown, the mixed model provides a practical and flexible alternative.

A plausible reason for the performance gap is that the heterogeneous training distribution makes it more challenging for a single score network to capture the full range of structural patterns with equal fidelity. As a consequence, the learned prior may become less accurate for individual classes than that obtained from dataset-specific training. Moreover, owing to computational constraints, the mixed model was not trained on the full collection of available samples, which may also contribute to the observed performance gap.

\begin{table}[th]
  \caption{Quantitative comparison of inversion results on OpenFWI velocity models using the separate and mixed models. Metrics include the relative $\ell_2$ error $e_{\ell_2}$, PSNR, and SSIM for the reconstructed velocity model (arrows indicate the preferred direction).}
  \label{tab:mix}
  \centering
  \small
  \setlength{\tabcolsep}{6pt}
  \renewcommand{\arraystretch}{1.10}
  \begin{tabular}{lccc}
    \toprule
    Dataset & $e_{\ell_2}\downarrow$ & PSNR$\uparrow$ & SSIM$\uparrow$ \\
    \midrule
    CurveVel-A (separate model)   & 3.50\% & 31.42 & 0.9039 \\
    CurveVel-A (mixed model)      & 3.87\% & 30.54 & 0.8906 \\
    \midrule
    FlatFault-A (separate model)  & 2.91\% & 32.74 & 0.9646 \\
    FlatFault-A (mixed model)     & 3.28\% & 31.72 & 0.9179 \\
    \midrule
    CurveFault-B (separate model) & 5.32\% & 27.95 & 0.7956 \\
    CurveFault-B (mixed model)    & 5.91\% & 27.04 & 0.7079 \\
    \bottomrule
  \end{tabular}
\end{table}

\subsubsection{Performance of the mixed model on the Marmousi2 model}
An important practical advantage of the mixed model is its ability to provide a unified prior for structures not explicitly represented by any single training class. To investigate this generalization capability, we applied the mixed model to the inversion of the Marmousi2 dataset \cite{martin2006marmousi2}. We emphasize that the training data contained no samples from Marmousi2. In the numerical experiments, local velocity models of different sizes were cropped from Marmousi2 and normalized to form dimensionless $71\times 71$ data matrices. For all local velocity models, we fixed
$k=100$, $c=0.1$, $\tau=0$, $\varepsilon=10^{-4}$, and $\gamma=0.55$,
while $\rho_0$ was chosen separately for each example.

\begin{figure}[th]
  \centering
  
  \setlength{\subfigwidth}{2.6cm}
  
  \setlength{\insetfigwidth}{2.7cm}

  \begin{tikzpicture}[
      subfig/.style={anchor=center, inner sep=0pt, outer sep=2pt, draw=black},
      subfig matrix/.style={
        matrix of nodes,
        column sep=0pt,    % 基础列间距
        row sep=0pt,       % 基础行间距
        nodes={subfig, minimum width=\subfigwidth, minimum height=\subfigwidth}
      },
      % 定义虚线样式
      no border/.style={subfig, draw=none},
      split line/.style={dashed, gray, line width=0.5pt},
      label text/.style={anchor=east, font=\small, inner sep=5pt, align=center}
    ]

    \matrix (subfig_matrix) [subfig matrix
      ,row 1/.style={row sep=8pt}
      ,row 2/.style={row sep=8pt}
      ,column 1/.style={column sep=5pt}
      ,column 2/.style={column sep=5pt}
      ,column 3/.style={column sep=5pt}
    ] {
      % 第一行：前4列 + 空白列 + 后4列
      \includegraphics[width=\subfigwidth, height=\subfigwidth, keepaspectratio]{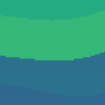} & 
      \includegraphics[width=\subfigwidth, height=\subfigwidth, keepaspectratio]{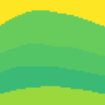} & 
      \includegraphics[width=\subfigwidth, height=\subfigwidth, keepaspectratio]{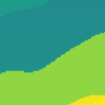} & 
      \includegraphics[width=\subfigwidth, height=\subfigwidth, keepaspectratio]{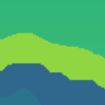} \\
      
      \includegraphics[width=\subfigwidth, height=\subfigwidth, keepaspectratio]{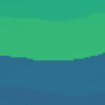} &
      \includegraphics[width=\subfigwidth, height=\subfigwidth, keepaspectratio]{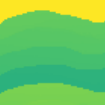} &
      \includegraphics[width=\subfigwidth, height=\subfigwidth, keepaspectratio]{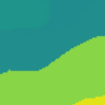} &
      \includegraphics[width=\subfigwidth, height=\subfigwidth, keepaspectratio]{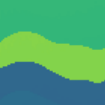}  \\

      \includegraphics[width=\subfigwidth, height=\subfigwidth, keepaspectratio]{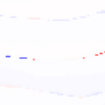} &
      \includegraphics[width=\subfigwidth, height=\subfigwidth, keepaspectratio]{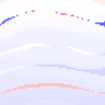} &
      \includegraphics[width=\subfigwidth, height=\subfigwidth, keepaspectratio]{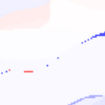} &
      \includegraphics[width=\subfigwidth, height=\subfigwidth, keepaspectratio]{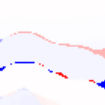}  \\
    };

    % \coordinate (inset_pos0) at ($(subfig_matrix-1-4.east)+(35pt,0)$);
    \coordinate (inset_pos1) at ($(subfig_matrix-2-4.east)+(35pt,0)$);
    \coordinate (inset_pos2) at ($(subfig_matrix-3-4.east)+(35pt,0)$);

    % \node[minimum width=\insetfigwidth, minimum height=\insetfigwidth] at (inset_pos0) {
    %   \includegraphics[width=\insetfigwidth, height=\insetfigwidth, keepaspectratio]{picture/marmousi_colorbar_1545.pdf}
    % };
    % \node[right=15pt, font=\small] at (inset_pos0) {\rotatebox{90}{velocity (m/s)}};

    \node[minimum width=\insetfigwidth, minimum height=\insetfigwidth] at (inset_pos1) {
      \includegraphics[width=\insetfigwidth, height=\insetfigwidth, keepaspectratio]{picture/marmousi_colorbar_1545.pdf}
    };
    \node[right=15pt, font=\small] at (inset_pos1) {\rotatebox{90}{velocity (m/s)}};

    \node[minimum width=\insetfigwidth, minimum height=\insetfigwidth] at (inset_pos2) {
      \includegraphics[width=\insetfigwidth, height=\insetfigwidth, keepaspectratio]{picture/marmousi_colorbar_3030.pdf}
    };
    \node[right=15pt, font=\small] at (inset_pos2) {\rotatebox{90}{velocity (m/s)}};

    \draw[split line] 
      ($(subfig_matrix-2-1.west)!0.5!(subfig_matrix-3-1.west)+(-7pt,0)$) -- ($(subfig_matrix-2-4.east)!0.5!(subfig_matrix-3-4.east)+(7pt,0)$);

    \draw[split line] 
      ($(subfig_matrix-1-1.west)!0.5!(subfig_matrix-2-1.west)+(-7pt,0)$) -- ($(subfig_matrix-1-4.east)!0.5!(subfig_matrix-2-4.east)+(7pt,0)$);

    \node[font=\small, anchor=center] at (
      $(subfig_matrix-1-4.east)+(8pt,0)$
    ) {\rotatebox{90}{Truth}};
    
    \node[font=\small, anchor=center] at (
      $(subfig_matrix-2-4.east)+(8pt,0)$
    ) {\rotatebox{90}{Reconstruction}};

    \node[font=\small, anchor=center] at (
      $(subfig_matrix-3-4.east)+(8pt,0)$
    ) {\rotatebox{90}{Difference}};
    
    % 顶部标签
    \node[above=38pt, font=\small] at (subfig_matrix-1-1) {Marmousi2-P1};
    \node[above=38pt, font=\small] at (subfig_matrix-1-2) {Marmousi2-P2 };
    \node[above=38pt, font=\small] at (subfig_matrix-1-3) {Marmousi2-P3};
    \node[above=38pt, font=\small] at (subfig_matrix-1-4) {Marmousi2-P4};

  \end{tikzpicture}
  \caption{Inversion results for the four local velocity models extracted from Marmousi2, denoted by Marmousi2-P1--P4. The first row shows the true velocity models $v_{\mathrm{true}}$. The second row presents the reconstructed velocity models $v_{\mathrm{rec}}$ obtained by the proposed method with the mixed model. The third row shows the corresponding reconstruction errors $v_{\mathrm{rec}}-v_{\mathrm{true}}$.}
  \label{fig:marmousi}
\end{figure}

We considered four local velocity models cropped from Marmousi2 at different physical scales. Their physical sizes, grid spacings, and the corresponding values of $\rho_0$ are summarized in Table~\ref{tab:marmousi_local_settings}. For all cases, we used exactly the same trained mixed model without any fine-tuning or distillation. We followed the same experimental setup as in the depth2 experiment in the multi-operation generalization section. The quantitative metrics and qualitative reconstructions are reported in Table~\ref{tab:marmousi_mix} and Figure~\ref{fig:marmousi}, respectively.

\begin{table}[th]
\centering
\caption{Experimental settings for the four local velocity models extracted from Marmousi2.}
\label{tab:marmousi_local_settings}
\small
\setlength{\tabcolsep}{8pt}
\renewcommand{\arraystretch}{1.10}
\begin{tabular}{ccccc}
\toprule
Local model & Physical size & $dx$ (m) & $dz$ (m) & $\rho_0$ \\
\midrule
Marmousi2-P1 & $700\,\mathrm{m}\times 175\,\mathrm{m}$  & 10              & 2.5  & 3.9  \\
Marmousi2-P2 & $700\,\mathrm{m}\times 263\,\mathrm{m}$  & 10              & 3.75 & 17.9 \\
Marmousi2-P3 & $1000\,\mathrm{m}\times 175\,\mathrm{m}$ & $\frac{100}{7}$ & 2.5  & 6.1  \\
Marmousi2-P4 & $1500\,\mathrm{m}\times 175\,\mathrm{m}$ & $\frac{150}{7}$ & 2.5  & 5.0  \\
\bottomrule
\end{tabular}
\end{table}

\begin{table}[th]
\centering
\caption{Quantitative comparison of inversion results obtained by the mixed model for the four local velocity models extracted from Marmousi2. Metrics include the relative $\ell_2$ error $e_{\ell_2}$, PSNR, and SSIM for the reconstructed velocity model (arrows indicate the preferred direction).}
\label{tab:marmousi_mix}
\small
\setlength{\tabcolsep}{6pt}
\renewcommand{\arraystretch}{1.10}
\begin{tabular}{lccc}
\toprule
Local model & $e_{\ell_2}\downarrow$ & PSNR$\uparrow$ & SSIM$\uparrow$ \\
\midrule
Marmousi2-P1 ($700\,\mathrm{m}\times 175\,\mathrm{m}$)  & 1.80\% & 35.95 & 0.9206 \\
Marmousi2-P2 ($700\,\mathrm{m}\times 263\,\mathrm{m}$)  & 3.57\% & 30.49 & 0.7553 \\
Marmousi2-P3 ($1000\,\mathrm{m}\times 175\,\mathrm{m}$) & 2.81\% & 33.21 & 0.8975 \\
Marmousi2-P4 ($1500\,\mathrm{m}\times 175\,\mathrm{m}$) & 5.31\% & 26.96 & 0.8536 \\
\bottomrule
\end{tabular}
\end{table}

The local velocity models extracted from Marmousi2 differ substantially from the OpenFWI datasets in structural features. Despite this mismatch, the proposed method using the mixed model still produces satisfactory reconstructions for local velocity models of different sizes. This is supported both qualitatively by Figure~\ref{fig:marmousi}, where the reconstructed velocity models preserve the main subsurface structures, and quantitatively by the relative $\ell_2$ error, PSNR, and SSIM reported in Table~\ref{tab:marmousi_mix}, which indicate good agreement between the reconstructed and true velocity models. Since no Marmousi2-specific samples were included in the training data, these results provide evidence of the cross-dataset generalization capability of the proposed approach. In this sense, the method offers a promising data-driven framework for mitigating the difficulty arising from limited training data in full waveform inversion.

\section{Conclusion}\label{sec:conclusion}
In this paper, we proposed a physics-guided diffusion framework for seismic velocity reconstruction in full waveform inversion. The method combines an optimal-transport (OT) data-consistency potential with a stabilized variable-metric guidance rule in the reverse diffusion process. On the data-consistency side, the OT potential incorporates bounded amplitude-adaptive reweighting together with a one-dimensional Wasserstein discrepancy evaluated through quantile functions. This construction reduces the dominance of large-amplitude early arrivals and alleviates the sensitivity of pointwise $\ell_2$ objectives to time and phase misalignment. On the guidance side, the scalar DPS correction is replaced by a diagonal preconditioner $P_i=\rho_iD_i$, where $\rho_i$ provides a step-dependent scalar schedule and $D_i$ supplies spatially adaptive scaling. This yields a more balanced and stable guidance mechanism for the heterogeneous sensitivities arising in FWI.

The numerical results on the OpenFWI benchmarks show that, under matched computational budgets, the proposed method improves reconstruction quality relative to both the standard DPS approach and the deterministic baselines. In particular, the reconstructions produced by the proposed method exhibit lower relative $\ell_2$ errors and fewer visible artifacts across the tested datasets.

Several directions remain for future work. First, the computational cost of physics-based guidance is still substantial, and it would be valuable to develop more efficient implementations and acceleration strategies, for example through multi-fidelity surrogate strategies. Second, an important extension is to consider more realistic forward models, including elastic, anisotropic, and attenuative media, as well as three-dimensional acquisition geometries.

\bigskip
\noindent{\bf Acknowledgments}\\
Zheng Ma is supported by NSFC Grant No. 12531016 and  Beijing Institute of Applied
Physics and Computational Mathematics funding HX02023-6. Xiong-Bin Yan is supported by NSFC Grant No. 12401563. Additionally, we also thank Shanghai Institute for Mathematics and Interdisciplinary Sciences (SIMIS) for their financial support. This research was funded by SIMIS under grant number SIMIS-ID-2025-ST. The authors are grateful for the resources and facilities provided by SIMIS, which were essential for the completion of this work.

\clearpage

\bibliographystyle{elsarticle-num} % 数字引用风格
\bibliography{references}          % references.bib 的文件名（不带 .bib）

\end{document}